\documentclass[11pt]{article}

\usepackage{amsmath,amsthm,amsfonts,amssymb,mathtools}  
\usepackage{textgreek,upgreek,bbm}  %bm creates heavy math
 
\usepackage{enumitem}
\usepackage{microtype}
\usepackage{booktabs}
\usepackage{graphicx}
 
 \theoremstyle{plain}
 \usepackage{titlesec}
 
 \usepackage[usenames,dvipsnames]{xcolor}
\usepackage[colorlinks = true, 
            linkcolor = RoyalPurple,
            urlcolor  = CadetBlue,
            citecolor = blue,
            anchorcolor = blue]{hyperref}

\usepackage{geometry}%
\def\urls#1{{\footnotesize\url{#1}}}

\makeatletter
\newcommand\gobblepars{%
    \@ifnextchar\par%
 {\expandafter\gobblepars\@gobble}%
{}}
\makeatother

\def\wham#1{\smallbreak\pagebreak[3]%
	\noindent\textup{\textbf{#1}}\ \ \gobblepars}

\def\whamb{\wham{$\bullet$}}

\def\mindex#1{\index[main]{#1}}

\DeclareFontFamily{U}{mathx}{\hyphenchar\font45}
\DeclareFontShape{U}{mathx}{m}{n}{<-> mathx10}{}
\DeclareSymbolFont{mathx}{U}{mathx}{m}{n}
\DeclareMathAccent{\widebar}{0}{mathx}{"73}

\def\temp{\upvarrho}
 
\def\spn{\operatorname{span}}
\def\hapfail{\widehat{p}_{\textsf{\tiny{fail}}}}

\newcommand{\Gammafn}{\operatorname{Gamma}}

\def\Obj{\Upgamma}
\def\optObj{\Obj^{\star}}

\def\sflog{\textsf{\tiny log}}
\newcommand{\logprec}{\preceq_{\sflog}}
\newcommand{\logasyeq}{ \asymp_{\sflog}}
\newcommand{\logsucc}{\succeq_{\sflog}}

\DeclareFontFamily{U}{bbold}{}
\DeclareFontShape{U}{bbold}{m}{n}
   {  <-5.5> bbold5 <5.5-6.5> bbold6 <6.5-7.5> bbold7
      <7.5-8.5> bbold8 <8.5-9.5> bbold9 <9.5-11> bbold10
      <11-15> bbold12 <15-> bbold17 }{}
\DeclareSymbolFont{bbold}{U}{bbold}{m}{n}
\DeclareSymbolFontAlphabet{\mathbbold}{bbold}

\makeatletter
\newcommand{\bigblock}[2][\LARGE]{%
  \mkern4mu
  \vcenter{\kern.3ex\hbox{#1$\m@th\mathstrut\mathbbold{#2}$}}%
  \mkern1mu
}
\makeatother

\newcommand{\ind}{{\mathchoice{\bigblock[\large]{1}}%
{\mathlarger{\mathbbold{1}}}%    %%   Seems to remain small! 
{\mathbbold{1}}%
{\mathbbold{1}}}}

\renewcommand{\ind}{\mathbbold{1}}

\def\Lip{L}  % Lipschitz constant.   
 
\def\Obj{\Upgamma}  %      No good: \Uplambda , \Upomega , L

\def\fee{\upphi}

\def\uQfee{\underline{Q}_{\raise 1.25pt\hbox{\tiny$\fee$}} }

\newcommand{\bbblot}{\raise1pt\hbox{\vrule height .4ex width .4ex depth .05ex}}
\long\def\defbox#1{\framebox[.9\hsize][c]{\parbox{.85\hsize}{%
\parindent=0pt
\baselineskip=12pt plus .1pt      % STYLE 
\parskip=6pt plus 1.5pt minus 1pt % CHANGES
 #1}}}

\long\def\beginbox#1\endbox{\subsection*{}%
\hbox{\hspace{.05\hsize}\defbox{\medskip#1\bigskip}}%
\subsection*{}}

\def\endbox{}

 \def\archival#1{} %Notes I'd like to save

\def\FRAC#1#2#3{\genfrac{}{}{}{#1}{#2}{#3}}

\def\ddt{{\mathchoice{\FRAC{1}{d}{dt}}%
{\FRAC{1}{d}{dt}}%
{\FRAC{3}{d}{dt}}%
{\FRAC{3}{d}{dt}}}}

\def\ddtp{{\mathchoice{\FRAC{1}{d^{\hbox to 2pt{\rm\tiny +\hss}}}{dt}}%
{\FRAC{1}{d^{\hbox to 2pt{\rm\tiny +\hss}}}{dt}}%
{\FRAC{3}{d^{\hbox to 2pt{\rm\tiny +\hss}}}{dt}}%
{\FRAC{3}{d^{\hbox to 2pt{\rm\tiny +\hss}}}{dt}}}}

\def\ddyp{{\mathchoice{\FRAC{1}{d^{\hbox to 2pt{\rm\tiny +\hss}}}{dy}}%
{\FRAC{1}{d^{\hbox to 2pt{\rm\tiny +\hss}}}{dy}}%
{\FRAC{3}{d^{\hbox to 2pt{\rm\tiny +\hss}}}{dy}}%
{\FRAC{3}{d^{\hbox to 2pt{\rm\tiny +\hss}}}{dy}}}}

\def\darrow{\buildrel{\rm d}\over\longrightarrow}

\newsavebox{\junk}
\savebox{\junk}[1.6mm]{\hbox{$|\!|\!|$}}

\def\det{{\mathop{\rm det}}}
\def\limsup{\mathop{\rm lim{\,}sup}}
\def\liminf{\mathop{\rm lim{\,}inf}}
\def\argmin{\mathop{\rm arg{\,}min}}

\def\bfmath#1{{\mathchoice{\mbox{\boldmath$#1$}}%
{\mbox{\boldmath$#1$}}%
{\mbox{\boldmath$\scriptstyle#1$}}%
{\mbox{\boldmath$\scriptscriptstyle#1$}}}}

\def\bfmY{\bfmath{Y}}

\def\bfmhhaY{\bfmath{\hhaY}} %\widehat{\widehat{Y}}}}
\def\bfmhhaY{\hbox to 0pt{$\widehat{\bfmY}$\hss}\widehat{\phantom{\raise 1.25pt\hbox{$\bfmY$}}}}

\newlength{\dhatheight}

\def\haX{\widehat X}

\def\clB{\mathcal{B}}
\def\clC{\mathcal{C}}
\def\clE{\mathcal{E}}
\def\clG{\mathcal{G}}
\def\clH{\mathcal{H}}    % Danger -- also used for  history

\def\clL{\mathcal{L}}
\def\clM{\mathcal{M}}

\def\clS{\mathcal{S}}

\def\clL{\mathcal{L}}

\def\clC{\mathcal{C}}

\newcommand{\eqdef}{\ensuremath{\stackrel{\hbox{\sf\tiny def}}{=}}}

\def\Prob{{\sf P}}

\def\Expect{{\sf E}}

 \def\epsy{\varepsilon}

\def\formtmp#1#2{{\vskip12pt\noindent\fboxsep=0pt\colorbox{#1}{\vbox{\vskip3pt\hbox to \textwidth{\hskip3pt\vbox{\raggedright\noindent\textbf{#2\vphantom{Qy}}}\hfill}\vspace*{3pt}}}\par\vskip2pt%
\noindent\kern0pt}}

\titleformat\subparagraph[runin]
                 {\normalfont\normalsize\bfseries}
                 {mypar}
                 {0pt}
                 {}{}
\titlespacing\subparagraph{0pt}%%             left margin spacing
                {.1ex minus 0.2ex}%% spacing above the \subparagraph
                {.75em}%%             horizontal offset from title (horizontal because the `runin` shape was used.)

 \definecolor{pcodecolor}{gray}{0.87}
 \definecolor{shadecolor}{gray}{.95}%  %.85 

\def\zat{0pt}
\newdimen\svparindent
\newenvironment{newshaded}{%
  \MakeFramed
{\FrameRestore}}%
{\endMakeFramed}
\AtBeginDocument{\renewenvironment{svgraybox}%
{\fboxsep=12pt\relax
 \begin{newshaded}\vspace*{7pt}%
 \list{}{\leftmargin=12pt\rightmargin=12pt\topsep=\zat\relax}%
 \expandafter\item\parindent=\svparindent
 \hskip-\listparindent}%
{\gobblepars\vspace*{7pt}\endlist\end{newshaded}}}%

{\endlist\end{newshaded}\pagebreak[3]%
}%

\newtheoremstyle{thm}{12pt}{15pt}%
     {\itshape}%         Body font
     {}%         Indent amount (empty = no indent, \parindent = para indent)
     {\bfseries}% Thm head font
     {}%        Punctuation after thm head
     {0pt}%{8pt}%     Space after thm head (\newline = linebreak)
     {\thmname{#1}\thmnumber{ #2.}\thmnote{ \textbf{(#3)}}\quad}% Thm head spec

{\end{list}}

\def\ass(#1:#2){(#1\ref{#1:#2})}

\def\ritem#1{
\item[{\sf \ass(\current_model:#1)}]
}

\newenvironment{recall-ass}[1]{% 
\begin{description}
\def\current_model{#1}}{
\end{description}
}

\def\sq{\hbox{\rlap{$\sqcap$}$\sqcup$}}
\def\qed{\ifmmode\sq\else{\unskip\nobreak\hfil
\penalty50\hskip1em\null\nobreak\hfil\sq
\parfillskip=0pt\finalhyphendemerits=0\endgraf}\fi}

\newcommand{\blot}{\vrule height 1.1ex width .9ex depth -.1ex }
\def\qedb{\ifmmode\blot\else{\vspace{-.2cm}\unskip\nobreak\hfil
\penalty50\hskip1em\null\nobreak\hfil\blot
\parfillskip=0pt\finalhyphendemerits=0\endgraf}\fi}

\def\PF{\wham{Proof:}}

\newcommand{\clThm}{%
  \ifvmode\noindent\fi
  \unskip
  \nobreak\hfill\penalty50
  \hskip1em\hbox{}\nobreak\hfill
  \hbox{\footnotesize$\square$}%
  \par
  \gobblepars
}

\newtheoremstyle{example}{15pt}{20pt}%
     {}%         Body font
     {}%         Indent amount (empty = no indent, \parindent = para indent)
     {\bfseries}% Thm head font
     {}%        Punctuation after thm head
     {1pt}%{8pt}%     Space after thm head (\newline = linebreak)
     {\thmname{#1}\thmnumber{ #2.}~\thmnote{\textit{\textbf{#3}}}%
     \\[.15cm]\unskip\nobreak}% Thm head spec

\newcounter{rmnum}

\newcounter{anum}

\newlist{EX}{enumerate}{1}
\setlist[EX]{font=\bfseries,label=\thechapter.\arabic*,%
topsep=0pt,%
partopsep=0pt,%
listparindent = 0pt,%
itemsep=3pt,%
parsep=3pt,%
labelwidth=0pt,%
labelsep=8pt,%
itemindent=25pt,%
leftmargin=0pt,% 
rightmargin=0pt,%
 }

\newcommand{\field}[1]{\mathbb{#1}}

\def\Re{\field{R}}

\def\argmin{\mathop{\rm arg\, min}}

\def\epsy{\varepsilon}

\def\haY{\widehat{Y}}

\def\hhaY{\hbox to 0pt{$\haY$\hss}\widehat{\phantom{\raise 1.25pt\hbox{Y}}}}

\def\hax{\widehat x}

\def\haY{\widehat Y}

\newcommand{\sfrac}[2]{#1/#2}

\def\Lip{\ell_\nabla}

\graphicspath{{./Figures/}}

\theoremstyle{definition}
 
\newtheorem{theorem}{Theorem}[section]

\newtheorem{proposition}[theorem]{Proposition}
\newtheorem{lemma}[theorem]{Lemma}

\newtheorem{assumption}{Assumption}[section]

\theoremstyle{remark}

\def\spn{\operatorname{span}}

\def\hapfail{\widehat{p}_{\textsf{\tiny{fail}}}}

 \newcommand{\Dbar}{E^\bullet} 

 \def\wGD {w_{\textsf{GD}}}

 \def\TGD {T_{\textsf{GD}}}
 \newcommand{\Tf}{T_{\textsf{f}}}

\def\optObj{\Obj^{\star}}

\def\gap{\textsf{\textup{r}}}
 
\def\RoA{\mathcal{S}}

 \def\wGD {w_{\textsf{GD}}}

\def\clA{\mathcal{A}}
\def\clH{\mathcal{H}}

\providecommand{\ind}[1]{\mathbbm{1}\{#1\}}
\newlength{\noteWidth}
\long\def\notes#1{\ifinner
	{\tiny #1}
\else
\marginpar{\parbox[t]{\noteWidth}{\raggedright\tiny #1}}
\fi}

\def\notes#1{}

 \def\nmax{n_{\rm max}}
\def\temp{\upvarrho}
\def\tempLG{\temp_{\textsf{\tiny max}}}

\def\tempSimA{\temp_{\textsf{\tiny SA}}}

\def\LG#1{\textsf{LG}_{#1}}
\def\LA#1{\textsf{L}_{#1}}

\def\LGFE#1{\textsf{LG}_{#1}^{\textsf{\tiny FE}}}
\def\LGRM#1{\textsf{LG}_{#1}^{\textsf{\tiny RM}}}

\def\LGSFE#1{\textsf{LGS}_{#1}^{\textsf{\tiny FE}}}
\def\LGSRM#1{\textsf{LGS}_{#1}^{\textsf{\tiny RM}}}

\def\SimA{\textsf{SimA}}

\newcommand{\sigmaRMswitch}{\sigma^0_{\textsf{\tiny RM}}}

\newcommand{\CompLGRM}{\clC^{\textsf{\tiny RM}}}
 \newcommand{\CompSimA}{\clC^{\textsf{\tiny SA}}}

\newcommand{\WLGRM}{W^{\textsf{\tiny RM}}}
\newcommand{\WSimA}{W^{\textsf{\tiny SA}}}

\def\nGD{n_{\rm\tiny GD}}
\def\uptauGD{\uptau_{\rm\tiny GD}}

\usepackage[nameinlink,noabbrev]{cleveref}

\Crefname{corollary}{Corollary}{Corollaries}
\Crefname{eqnarray}{eq.}{eqs.}
\Crefname{equation}{eq.}{eqs.}

\Crefname{figure}{Fig.}{Figs.}
\Crefname{tabular}{Tab.}{Tabs.}
\Crefname{table}{Tab.}{Tabs.}
\Crefname{lemma}{Lemma}{Lemmas}
\Crefname{proposition}{Prop.}{Propositions}
\Crefname{theorem}{Thm.}{Thms.}
\Crefname{definition}{Def.}{Defs.} 
\Crefname{section}{Section}{Sections}
\Crefname{assumption}{Assumption}{Assumptions}
\Crefname{exmp}{Example}{Examples}
\Crefname{exercise}{Exercise}{Exercises}

\title{Fast PAC Global Optimization via Restarted Langevin: 
\\
Exploration, Exploitation, and Degenerate Cooling
}
\author{Ioannis Kontoyiannis and Sean Meyn%
\thanks{I.K.\ is with the
Statistical Laboratory, Centre for Mathematical Sciences, 
University of Cambridge, Wilberforce Road, Cambridge CB3 0WB, 
UK. Email: yiannis@maths.cam.ac.uk. 
S.M.\ is with the  Department of Electrical and Computer Engineering,  University of Florida, Gainesville,  FL, USA.   Email: meyn@ufl.edu. 
\\
  Financial support:
  I.K. was supported in part
by the EPSRC-funded INFORMED-AI project EP/Y028732/1.  S.M.\ was supported by ARO W911NF2410389 and 
the  NSF award CCF-2306023 }
}

\begin{document}
\maketitle

\begin{abstract}

This paper addresses a fundamental question in non-convex optimization:
\textit{How should a stochastic optimizer allocate computation between global exploration and local exploitation?}

We study this question in continuous time, using Langevin dynamics for global exploration and deterministic gradient flow for local exploitation.  The simplest algorithm in the resulting class is the best-state Langevin--gradient method: several Langevin trajectories explore the objective landscape, the best state encountered is retained, and a single gradient trajectory then refines this state to high terminal accuracy.

Our main objective is to characterize the computational work required to achieve a prescribed accuracy with prescribed confidence.  We develop low-temperature approximations for this PAC work--accuracy tradeoff, both for the best-state method and for related Langevin schemes.  Global exploration is governed by   energy barriers, through spectral-gap and Eyring--Kramers asymptotics.  The resulting approximations expose the competing effects of temperature, computational budget, and the number of exploratory trajectories, and quantify their dependence on the confidence parameter.  In contrast, terminal accuracy is largely separated from global exploration: once a suitable region of attraction has been reached, gradient flow requires only $O(\log(1/\epsy))$ additional work.

We also investigate corrections that disappear at logarithmic low-temperature precision but can become important with increasing dimension.  Numerical experiments on challenging multimodal objectives serve both to compare algorithms and to test the quantitative predictions of the asymptotic theory.  They identify regimes in which the low-temperature approximations provide useful guidance, as well as regimes exhibiting substantial temperature sensitivity and potentially significant high-dimensional limitations.

\end{abstract}

 \clearpage
 
 \setcounter{tocdepth}{2}
 \tableofcontents

 \clearpage

\section{Introduction}
\label{s:intro}

A central difficulty in non-convex optimization is the allocation of computation between exploration and exploitation. Random perturbations can move an iterate between attraction regions, but continued exploration prevents accurate local convergence. Deterministic gradient methods have the opposite strengths: once initialized in a favorable region they can converge rapidly, but they provide no mechanism for escaping a poor basin.

Denoting the objective $\Obj\colon\Re^d\to\Re$,  in evaluating an algorithm for global optimization  we adopt terminology from statistical learning theory.  
For prescribed accuracy $\epsy>0$ and confidence level $1-\delta$,
we say that an algorithm satisfies the $(\epsy,\delta)$-PAC  [Probably Approximately Correct]
performance requirement if its output $\haX$ satisfies
\begin{equation}
    \Prob\{\Obj(\haX)-\optObj>\epsy\}\leq\delta.
\label{e:performance-objective}
\end{equation}
Given any algorithm producing a sequence of approximations  $\{ \haX_n\}$, the associated PAC sample complexity is denoted
\begin{equation}
    n^\star(\epsy,\delta)
    \eqdef
    \min\left\{
        n:
        \Prob\{\Obj(\haX_n)-\optObj>\epsy\}\leq\delta
    \right\}.
\label{e:nstarIntro}
\end{equation}

For example, stochastic gradient descent (SGD) with convex objective can achieve $O(1/n)$ rate of convergence for the excess cost $\{ \Obj(\theta_n)-\optObj \}$.    Sharp bounds, such as  developed in \cite{bacmou13,shebeldummounausam24}, lead to   
\begin{equation}
        n^\star(\epsy,\delta)
        = O( L /\epsy),
        \qquad
        L=\log(1/\delta).
\label{e:SGDcomplexity}
\end{equation}

For the preferred approaches surveyed in this paper, analysis reveals a strong separation between global exploration and terminal accuracy. 
The stochastic phase need only locate a suitable region of attraction; once this occurs, gradient flow supplies accuracy $\epsy>0$
with an additional cost only logarithmic in $1/\epsy$, resulting in significant improvement over \eqref{e:SGDcomplexity}.

Hence design for successful exploration is the main algorithmic goal of this paper.

\subsection{Assumptions}

For the purposes of analysis,  algorithms are formulated in continuous time,  as variants of stochastic-gradient dynamics. 
For  fixed temperature $0<\temp<\infty$, the Langevin diffusion is the stochastic differential equation, 
\begin{equation}
    dX_t=-\nabla \Obj(X_t)\,dt+\sqrt{2\temp }\,dB_t 
    \label{eq:langevin}
\end{equation}
On setting the temperature to zero 
 we obtain the gradient flow, 
\begin{equation}
    \ddt x_t=-\nabla \Obj(x_t) \,,
\label{e:GF}
\end{equation}

\begin{subequations}

A favored approach   is simply described as follows:

\begin{quote}
\wham{Best-state switching \rm ($\LGRM{\temp,N}$).}
Run   $N$ Langevin trajectories until the deterministic exploration
time $\Tf^0$, and define
\begin{equation}
    (i^\star,\sigmaRMswitch)
    \in
    \argmin  \{  
        \Obj(X_t^i) :  1\leq i\leq N\, , \ 0\leq t\leq \Tf^0 \} 
    \label{e:PerformanceStop}
\end{equation}
After this computation, 
exploration is halted and the temperature is set to zero:
\begin{equation}
 \frac{d}{dt} X_t^{i^\star} = - \nabla\Obj(X_t^{i^\star}) \,, \qquad     \sigmaRMswitch < t \le \Tf
\label{e:GF_LGRM}
\end{equation}
\qed
\end{quote}
The meta-parameters $(\temp, N, \Tf^0 , \Tf)$ are chosen in advance---theory in this paper provides guidance for selection. 
The ``\textsf{RM}'' in $\LGRM{\temp,N}$
 stands for \textit{running min}.

\end{subequations}

The assumptions imposed below introduce two key quantities with the following interpretations:
\whamb  $\Dbar$:  the  \textit{dominant energy barrier},  controlling how hard it is to explore globally.
\whamb $\eta$:  the  \textit{attraction margin},  an optimality gap ensuring  local exploitation is successful.
\\
Implementation does not require a~priori knowledge of these quantities.  

\smallskip

For positive temperature, the diffusion has a unique invariant probability measure,
\begin{align}
    \pi_\temp(dx) & =Z_\temp^{-1} \exp\Big\{- \frac{1}{\temp} \Obj(x) \Big\}\,dx.
\label{e:pi_temp}
\\
    Z_{\temp}
&=
    \int_{\Re^d} \exp\Big\{- \frac{1}{\temp} \Obj(x) \Big\}\,dx 
\label{e:Ztemp}
\end{align}
\Cref{a:jacquot} ensures that the normalizing constant $Z_{\temp}$  is finite for all  $\temp>0$.

\begin{assumption}[Primary assumptions on the objective]
\label{a:jacquot}
The objective $\Obj\colon\Re^d\to\Re$ satisfies the following:

\wham{(i)}  It is twice continuously differentiable ($C^2$).

\wham{(ii)}
It is coercive and satisfies 
\[
      \lim_{\|x\|\to\infty}
    \frac{\Obj(x)}{\log(1+\|x\|)}
    =
    \lim_{\|x\|\to\infty}\|\nabla\Obj(x)\|=\infty,
    \qquad
    \sup_{x\in\Re^d}\|\nabla^2\Obj(x)\|<\infty .
\]

\wham{(iii)} 
There exists a closed and bounded set $\RoA\subset\Re^d$ such that
$
    \RoA=\overline{\RoA^\circ}$, and moreover  

\whamb
The  attraction margin $\eta$ is strictly positive, where   
\begin{equation}
    \eta
    \eqdef
    \inf_{x\in\RoA^c}
    [\Obj(x)-\optObj]
\label{e:etaRoA}
\end{equation}

\whamb     
There are positive  constants  $\upmu_{\RoA}$, 
  $c_{\RoA}$ such that for any solution to \eqref{e:GF} with $x_0\in\RoA$,
    \begin{equation}
 \Obj(x_t)-\optObj  \le c_{\RoA} \exp \{ - \upmu_\RoA t \} \,, \quad t\ge  0
\label{e:BexpGF}
\end{equation}
\clThm
 \end{assumption}

In the optimization approaches introduced in this paper, the set $\RoA$ is regarded as a \textit{target set}: ideally, the gradient descent phase begins at the first time this set is reached.  
We do not require existence of a test for membership in the proposed algorithms, but  the value of $\eta$,
$\upmu_{\RoA}$, and
  $c_{\RoA}$ appear in  complexity bounds.   
  
  \wham{Level sets as target sets}   
  It is reasonable to take the sublevel set  $\RoA =
    \RoA_\eta$, with
  \begin{equation} 
    \RoA_\eta
    =
    \{x:\Obj(x) \le \optObj+\eta\} \,,
\label{e:RoAsublevel}
\end{equation}
for which the attraction margin \eqref{e:etaRoA} is exactly $\eta$.
The gradient flow satisfies
\[
 \frac{d}{dt}\Obj(x_t) = -\|\nabla\Obj(x_t)\|^2\le0 \, ,
 \]
 so that
 the set $\RoA_\eta$ is forward invariant: $x_t \in \RoA_\eta $ for all $t\ge0$ provided $x_0 \in \RoA_\eta $. 
 
Suppose that for sufficiently small
$\eta>0$ we can find $\upmu_\eta>0$ such that
$$ 
\|\nabla\Obj(x)\|^2 \ge \upmu_\eta [\Obj(x) -  \optObj] \,, \quad x\in \RoA_\eta   
$$
We then obtain  the required assumptions for \Cref{a:jacquot}~(iii). In particular, 
\eqref{e:BexpGF} holds with $c_{\RoA}  = \eta$ and $\upmu_\RoA = \upmu_\eta$.

\wham{Notation surrounding spectra and energy landscape}

We let   $\gap_{\!\temp}$ denote the spectral gap of the differential generator in $L_2$.   
Under the assumptions imposed in this paper it is strictly positive for positive temperature.    
Moreover,  the differential generator  has a discrete and real spectrum, assumed ordered so that  $\lambda_{i+1}\le \lambda_i$ for $i\ge 0$. 
We have $\lambda_1=0$ and $\lambda_2 = -     \gap_{\!\temp}<0$.

 The spectral gap is approximated for low temperature based on  the  dominant energy barrier
 \begin{equation}
    \Dbar
    \eqdef 
    \sup_{x,y\in\Re^d}
    \bigl\{
        H(x,y)-\Obj(x)-\Obj(y)+\optObj
    \bigr\}.
\label{e:EnergyBarrier}
\end{equation}
where, for $x,y\in\Re^d$,   the \textit{communication height} is defined by
\[
    H(x,y)
    \eqdef 
  \inf _\gamma    \max_{0\leq s\leq1} \Obj(\gamma(s)),
\]
in which the infimum is over continuous paths $\gamma$ 
satisfying
 $ \gamma(0)=x$  and $  \ \gamma(1) =y$.   
 An illustration is provided in \Cref{s:RastriginAppendix}.

We use the standard term \textit{energy landscape} for the geometric
structure induced by $\Obj$, including its minima, saddle points,
regions of attraction, and communication heights; see, for example,
\cite{huimeysch04a,MenzSchlichting2014}.

%  Fluff
%We will see that quantitative information about the objective landscape is needed to
%predict \emph{how fast} the algorithm succeeds, but not necessarily to
%guarantee eventual success.
 
 \subsection{Complexity}

The complexity measure used throughout is the total time required to
meet the PAC criterion~\eqref{e:performance-objective}; this quantity
is denoted by $\clC$.   This is  used as a proxy for
computational effort and is not intended to account for
implementation-level costs.

 In the comparisons based on simulation surveyed in \Cref{s:num} we take into account the realities of implementation. 
In particular, if $\nabla\Obj$ is globally Lipschitz with constant
$\Lip$, the natural stability scale for an explicit Euler discretization
is $\Lip^{-1}$.  Once a common discretization rule and numerical accuracy
are fixed, the number of gradient and noise updates is proportional,
up to method-independent and discretization-dependent factors, to the
aggregate simulated time.

Approximations of $\clC$ for the various methods are developed in
\Cref{s:low}.  These are low-temperature logarithmic asymptotics; for
the parallel schemes the relevant regime includes $N/L\to0$, where
$L=\log(1/\delta)$.  The propositions in \Cref{s:low} give precise
statements of the limiting regimes underlying the summary below.

 In \Cref{s:highd}
 we explain how dimension may enter through constants that are hidden in logarithmic approximations.

\Cref{t:4LGtempN} concerns $\LGRM{\temp,N}$,  where for fixed $N$ the dependency on $\delta,\epsy$   is made explicit
\begin{equation}
    \CompLGRM  \approx    N\exp\left(\frac{ LE}{\eta N} \right)      +     \TGD  (\epsy)
\label{e:CompLGlow_tempN}
\end{equation}
where   $E>\Dbar$ is a slack parameter in the
low-temperature bounds,  and  
\begin{equation}
        \TGD = 
    \TGD  (\epsy)
    \eqdef
    \frac{1}{\upmu_{\RoA}}\log\left(\frac{c_{\RoA}}{\epsy}\right) 
\label{e:TGD}
\end{equation}

Hence the first term in \eqref{e:CompLGlow_tempN} is the price of exploration; the second is the price of exploitation. Restart changes the exploration exponent through $N$, while terminal accuracy appears only through the gradient term.
The bound  
follows from a choice of exploration temperature and horizon
that balances two competing low-temperature effects: concentration near
the attraction region, governed by $\eta$, and  exploration 
governed by the energy barrier $\Dbar$.  
The required landscape  information is essentially the same as that entering classical simulated
annealing, together with a lower bound on the attraction margin $\eta$
in \Cref{a:jacquot}.

The unconstrained minimizer of the exploration term in
\eqref{e:CompLGlow_tempN} would take $N$ proportional to
$L=\log(1/\delta)$.  This lies outside the low-temperature regime
required by the approximation, 
which assumes $N/L\to0$.  It is shown in \Cref{t:BetterBddsThroughRestart}  
that choosing $N$ to grow appropriately sub-linearly with $L$ leads to 
\begin{equation}
     \CompLGRM
   =E  L^{1+o(1)}
    +O\Big(\log\frac1\epsy\Big).
 \label{e:CompLG}
\end{equation}
Thus the low-temperature theory yields
a separation of contributions of confidence from
terminal accuracy,  with 
 exploration complexity
arbitrarily close, on a logarithmic scale, to linear dependence on
$\log(1/\delta)$

%Consider for illustration the Rastrigin objective, which is the focus of \Cref{s:rast}.    
%It is shown in \Cref{t:Rastrigin}  that   $  \Dbar \approx 19.26$ for any dimension $d$, and $\eta\approx 0.5$.  Hence on choosing  $E = 20 > \Dbar$ and
%   $N = \lfloor   E/\eta\rfloor$ in 
%\eqref{e:CompLGlow_tempN}, 
% \[
% \CompLG \approx 
%40    \log\left(    \frac{1}{\delta}  \right)
%  \frac{1}{\delta}      +    \TGD  (\epsy)
%   \]

The contrast with classical simulated annealing is sharp.  There the
temperature remains coupled to terminal accuracy, leading on the same
logarithmic scale to an exponentially unfavorable dependence on
$1/\epsy$.   \Cref{t:SimAnn}, based on bounds from \cite{TangZhou2023},    establishes upper and lower bounds on simulation time to conclude
\begin{equation}
    \CompSimA \approx \exp\Big(\frac{LE}{\epsy}\Big),
\label{e:CompSimA}
\end{equation}
where $E>\Dbar$ is a parameter appearing  in the standard cooling schedule.   

\begin{figure}[h!]
	\centering
	\includegraphics[width=\hsize]{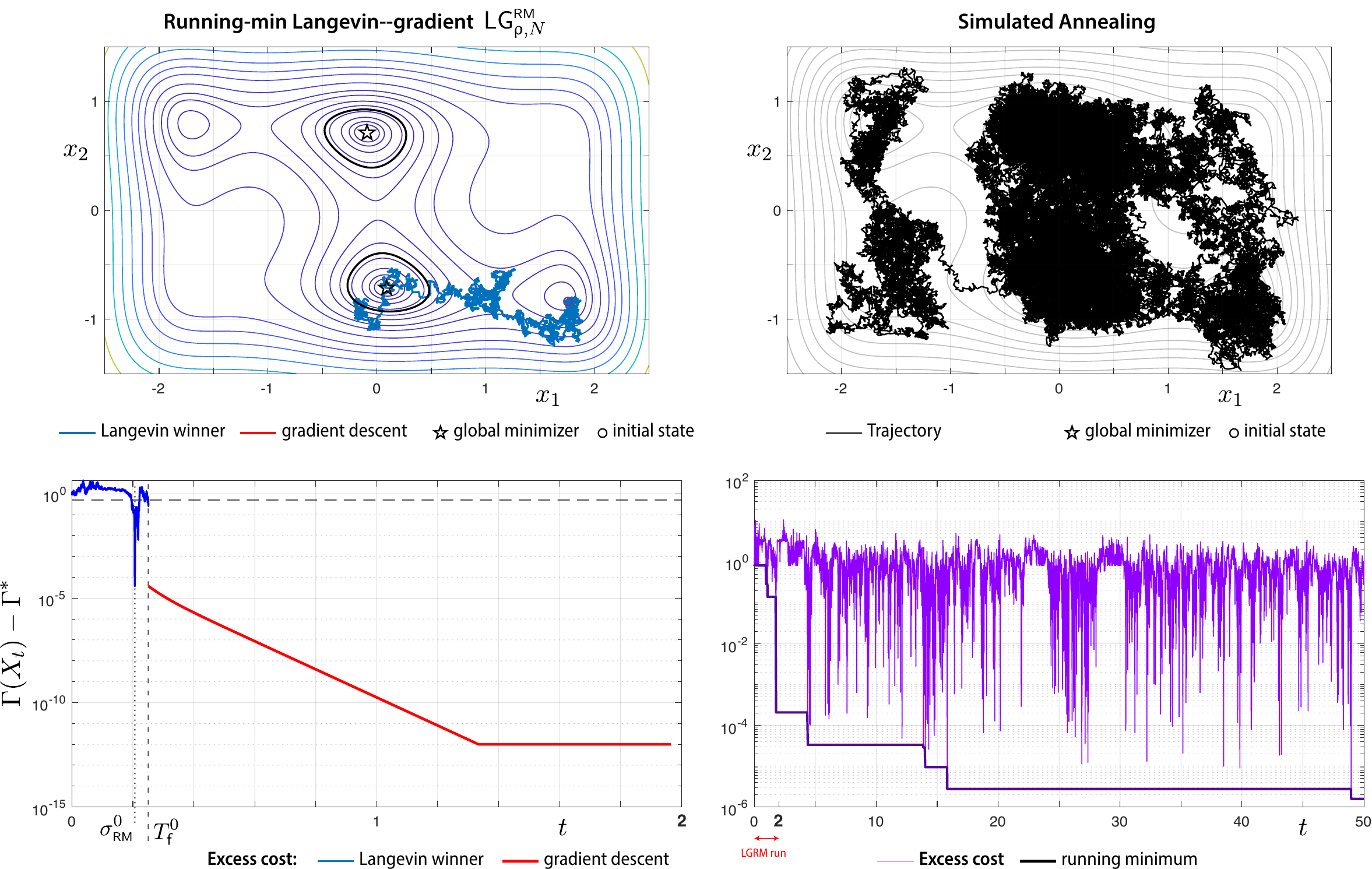}  
\caption{
Representative sample-path comparison for the modified six-hump camel
objective.  
The value $\temp=1$ was used in $\LGRM{\temp,N}$  to obtain good visibility of the trajectory--a larger value is much more effective in this example, especially combined with $N>1$.}
		\label{f:LG24temp1vsSimA}
\end{figure}

\Cref{f:LG24temp1vsSimA}
 provides a comparison of $\LGRM{\temp,N}$  and simulated annealing for an objective on $\Re^2$.   
Details on this experiment and others are found in \Cref{s:num}.

\subsection{Contributions}

The main contributions are as follows.

\wham{PAC exploration--exploitation complexity.}
We formulate the allocation between stochastic global exploration and
deterministic local exploitation through an $(\epsy,\delta)$-PAC work
criterion.    For Langevin--gradient schemes, the resulting approximation \eqref{e:CompLGlow_tempN}
separates the cost of reaching a suitable region of attraction
from the cost of terminal accuracy.   This conclusion is a corollary to \Cref{t:4LGtempN}, and related bounds for other approaches is the topic of \Cref{s:low}.

\wham{Quantifying the value of independent exploration.}
Independent restart is classical in global optimization
\cite{mus97,SchutteHaftkaFregly2007}.
For Langevin exploration,  \Cref{t:3LtempN,t:4LGtempN}  quantify how the number $N$
of independent trajectories alters the confidence--work tradeoff and
permits exploration at a correspondingly higher temperature.
For the Langevin--gradient scheme, choosing $N$ to increase sufficiently
slowly with $L=\log(1/\delta)$ while maintaining $N/L\to0$ gives
exploration complexity $L^{1+o(1)}$.

\wham{Best-state Langevin--gradient optimization.}

The best-state variant studied in \Cref{s:LGmethods} retains the
lowest-cost state encountered over all exploratory trajectories and
times and initializes a single gradient trajectory from this state.
Consequently a favorable transient visit is not lost, and no online
test for membership in the target region is required; see
\Cref{s:prior-knowledge}.  Its work advantage over continuing all $N$
trajectories through the gradient stage is made explicit in
\eqref{e:CompLGlow_tempN} and \Cref{t:4LGtempN}.

\wham{Temperature design from landscape geometry.} 

The low-temperature analysis in \Cref{s:low} combines Gibbs
concentration with spectral-gap and Eyring--Kramers asymptotics.
The resulting bounds distinguish the attraction margin $\eta$,
which controls concentration near states from which local optimization
is effective, from the dominant barrier $\Dbar$, which controls the
time scale of   exploration.  \Cref{t:4LGtempN} turns these
two geometric quantities into explicit scaling rules for temperature,
run length, confidence, and the number of exploratory trajectories.

 \wham{High-dimensional restart and temperature adaptation.}
Classical order-statistic arguments imply a severe dimensional penalty
for best-of-$N$ independent search; in our setting this appears as the
$N^{-2/d}$ rate for fixed-temperature Gibbs samples--see \Cref{t:RunMax_dimension,t:aStatsRunningMin}.
We show in \Cref{t:running-min-temp-scaling} that allowing the Gibbs temperature to depend on the
requested accuracy changes the $\epsy^{-d/2}$ sample-complexity
dependence to $\epsy^{-1}$ for each fixed $d>2$, although with a
dimension-dependent prefactor that grows rapidly with $d$.
This also illustrates why logarithmic low-temperature asymptotics can
conceal substantial high-dimensional costs.

%\wham{Beyond logarithmic low-temperature asymptotics.}
%We identify dimension-dependent prefactors that disappear on the
%logarithmic scale and show that favorable dependence on
%$\epsy$ and $\delta$ need not imply dimension-robust computational
%performance.

\wham{Quantitative assessment of the asymptotic theory.} 

The experiments in \Cref{s:num} are designed   to test the quantitative predictions used in their
design.  The six-hump camel, Rastrigin, and Griewank studies examine
Gibbs concentration, spectral-gap scaling, temperature selection,
and PAC work, and identify both regimes in which low-temperature
theory provides useful guidance and regimes in which its predictions
are quantitatively unreliable.

\subsection{Related work}
\label{s:related}

The work most closely related to ours is that of Dong and Tong
\cite{DongTong2021}.  They also exploit the complementary roles of
Langevin dynamics and gradient descent in non-convex optimization:
the former provides global exploration, while the latter gives rapid
local convergence.  Their principal algorithm couples the two through
an exchange mechanism, and they establish high-probability linear
convergence under a local strong-convexity condition near a unique
global minimizer.

In this paper and prior work,   the gradient flow phase \eqref{e:GF_LGRM} is known as \textit{quenching}.  

Our objectives are complementary.  We take the exploration--exploitation
decomposition itself as the object of study and ask how temperature,
restart, exploration time, and confidence should be allocated when
performance is measured by PAC computational work.  Low-temperature
analysis leads to an explicit separation between the   cost
of reaching a suitable region of attraction and the subsequent cost of
local optimization.  This exposes the distinct roles of the dominant
energy barrier $\Dbar$, the attraction margin $\eta$, the number $N$
of exploratory trajectories, and the confidence parameter $\delta$.
We also investigate dimension-dependent terms hidden by logarithmic
asymptotics and test quantitatively the range over which these
approximations are predictive.
 
 The analytical approach here is also different from
\cite{DongTong2021}.  Their principal high-probability bound is
obtained through the geometric drift criterion  (V4) of  \cite{meytwe93b,dowmeytwe95a}
to bound the return times to a bounded region.  This is combined with a log-Sobolev inequality to control the
exploration time.   Our analysis instead begins with the multiplicative drift condition
(DV3) of \cite{konmey05a,konmey17a}.  Under \Cref{a:jacquot}, this condition holds and implies discreteness of
the spectrum of the Langevin generator.  Reversibility then gives
$L_2(\pi_\temp)$ contraction at the exact spectral-gap rate
$\gap_{\!\temp}$, while regularization of the diffusion semigroup converts
this $L_2$ estimate into a pointwise bound on transition
probabilities; see \Cref{t:L2}.  Low-temperature spectral
asymptotics subsequently connect this rate to the dominant energy
barrier $\Dbar$.  Thus the spectral gap is not introduced as an
unspecified functional-inequality constant, but is tied directly to
the geometry of the objective landscape.  
 
A second line of work treats the exploration mechanism itself as a
design variable.  Gao, Xu, and Zhou \cite{GaoXuZhou2022} formulate
temperature selection for Langevin diffusion as a stochastic control
problem, leading to state-dependent randomized temperature policies.
Intermittent diffusion \cite{ChowYangZhou2013} alternates stochastic
perturbation and quenching, while more recent
replica-exchange methods constrain high-temperature exploration to
avoid costly excursions into the tails \cite{ZhengEtAl2024}.
These approaches reinforce a point that is central here: more
exploration is not uniformly better.  Temperature controls a genuine
tradeoff between rapid movement across the energy landscape and
concentration near regions from which local optimization is effective.
Our aim is to quantify this tradeoff directly in terms of PAC work.

Independent restart also has a long history in stochastic global
optimization; see, for example, \cite{mus97} and subsequent work
on stopping and restarting stochastic searches.  Classical multi-start
methods compare repeated independent searches with a single prolonged
search, often followed by local optimization.  The role of restart here
is more specific: independent Langevin trajectories alter the
temperature--confidence tradeoff through the product structure of the
failure probability.  Combined with a separate deterministic
exploitation phase, this leads to the complexity decomposition studied
in this paper.

The low-temperature analysis builds on the classical theory of
simulated annealing
\cite{GemanHwang1986,ChiangHwangSheu1987,HolleyKusuokaStroock1989}.
In the finite-state setting, Hajek \cite{haj88} identified the depth of
the deepest nonglobal local minimum as the critical quantity governing
logarithmic cooling.  More recent refinements, including
\cite{TangZhou2023,tanwuzho24}, provide spectral and low-temperature
estimates that we adapt here.  Classical simulated annealing uses these
ideas to determine a cooling schedule that eventually resolves terminal
accuracy.  We instead use the same energy-landscape information to ask
how long stochastic exploration should continue before it is preferable
to hand the problem to a deterministic local method.

\section{Langevin--gradient methods}
\label{s:LGmethods}

\begin{subequations}

Common to each of the Langevin--gradient methods  is a deterministic time $\Tf >0$ at which   the estimate of an optimizer is defined through minimization.  
For analysis we take 
 \begin{equation}
    \haX
    \in
    \argmin \{   \Obj(X_{\Tf}^i)  :  1\le i\le N\}  
    \label{e:bestOfN}
\end{equation}
With the understanding that any implementation of these methods will be based on a discrete-time recursive algorithm, the running minimum over all observed values is a practical alternative:  
 \begin{equation}
    \haX
    \in
    \argmin \{   \Obj(X_t^i)  :  1\le i\le N\,, \  0\le t\le \Tf\}  .   
    \label{e:bestOfN+T}
\end{equation}

\end{subequations}

\subsection{Exploration and exploitation}

The exploration component of these approaches is achieved through the realization of $N$ independent diffusions up to a terminal time $\Tf^0$.   
This data is then used to obtain an initial condition $\hax^0$ to define the exploitation phase defined by a gradient-flow.

Several variations are considered,  differing in the form of the diffusion and the mechanism to determine $\hax^0$.
We begin with the most basic approach, in which $N$ gradient flows are obtained in the exploitation phase.

\wham{Parallel Langevin--gradient \rm ($\LG{\temp,N}$).}
Run $N$ independent Langevin diffusions at temperature $\temp$ until
the switching time $\Tf^0$, then set $\temp=0$ and continue under
gradient flow.
 
The temperature is thus piecewise constant: for a specified switching time
$\Tf^0 = \Tf^0(N)$,
\[
    \temp(t)
    =
    \begin{cases}
        \temp>0, & 0\leq t<\Tf^0,\\
        0,      & \Tf^0\leq t\leq \Tf .
    \end{cases}
\]
Hence all $N$ trajectories follow the gradient flow for $t\ge \Tf^0$,
and we view  $0\le t<\Tf^0$ as an interval for exploration.

For analysis, the choices of $\Tf^0$ and $\Tf$ are based on the
property of the set $\RoA$ in \eqref{e:BexpGF}.  This bound implies that
\[
    \Obj(x_t)-\optObj\leq\epsy,
    \qquad x_0\in\RoA \,, \quad     t\geq
        \TGD \,
\]
where $\TGD $ is given in \eqref{e:TGD}, 
giving
$    \TGD =O(\log(1/\epsy))
$, and by definition
$    \Tf=\Tf^0+\TGD  $.

The switching time $\Tf^0$  is chosen so that
\[
    \Prob\{E_{\RoA}\}\geq1-\delta,
    \qquad 
    \textit{where} \ \ 
    E_{\RoA}=\bigcup_{i=1}^N E_{\RoA}^i,
    \qquad
    E_{\RoA}^i=\{X^i_{\Tf^0}\in\RoA\}.
\]
On the event $E_{\RoA}^i$ we have $
    \Obj(X^i_{\Tf})-\optObj\leq\epsy$,  and hence on the union of these events  
we achieve the desired bound,
\begin{equation}
    \Obj(\haX_{\Tf})-\optObj\leq\epsy \quad \textit{on the event
 $E_{\RoA}$.}
\label{e:LGdesired}
\end{equation}

Assuming that $N$ and   $\Tf>0$ are chosen so that   the 
PAC criterion~\eqref{e:performance-objective} is achieved,   the complexity in units of total time,  as defined in \Cref{s:intro}, is  
\begin{equation}
    \clC 
    =  N\Tf = 
    N [\Tf^0
    +
    \TGD  ]
    \label{e:CompLG0}
\end{equation}

\wham{First-entry switching \rm ($\LGFE{\temp,N}$).}
 
In this approach   Langevin
exploration is stopped as soon as the first trajectory reaches the target region.
This approach assumes that membership in the target region can be detected online, which is a reasonable assumption of the objective value $\optObj$ is known a-priori.

For trajectory $i$, define
\begin{equation}
    \tau_{\RoA}^i
    =
    \min\{t\geq0:X_t^i\in\RoA\} 
\label{e:BasinStop}
\end{equation}
At the first index
$i^\star$ for which $\tau_{\RoA}^{i^\star}\le \Tf^0$ is observed, exploration is halted and the temperature is set to zero to obtain,
\[
 \frac{d}{dt} X_t^{i^\star} = - \nabla\Obj(X_t^{i^\star}) \,, \qquad  t>      \tau_{\RoA}^{ i^\star}
\]
Failure to reach $\RoA$ by $\Tf^0$ for any $i$ 
is regarded as failure of the exploration phase.   

This version can reduce   simulation work relative to the deterministic switching rule,
since unsuccessful Langevin trajectories need not be continued once
one successful trajectory has been identified.

This 
 first-entry variant has an additional qualitative property that
does not require a prescribed confidence level or switching horizon.
Taking $N=1$ and removing the upper bound on the switching time gives
the following.

\begin{proposition}[Almost-sure exponential convergence]
\label[proposition]{t:LG-hitting}
Suppose that \Cref{a:jacquot} holds and fix $\temp>0$.   
Consider $\LGFE{\temp,1}$ with no prescribed
upper bound on the exploration time. 
That is,  run the Langevin diffusion 
until
\[
    \tau_{\RoA}
    \eqdef
    \min\{t\geq0:X_t\in\RoA\},
\]
and for $t>     \tau_{\RoA}$   the temperature is set to zero.

Then the resulting process
satisfies, for  $ t\geq\tau_{\RoA}$,
\begin{equation}
    \Obj(X_t)-\optObj
    \leq
\clG e^{- \upmu_{\RoA} t} \,, \qquad \clG =    \eta     e^{\upmu_{\RoA} \tau_{\RoA} }
    \label{e:semi-exponential-convergence}
\end{equation} 
In particular,
$
    \Obj(X_t)\longrightarrow\optObj$ a.s.,  
with exponential convergence following the random switching time.
Moreover,  for every initial condition $X_0 = x$, 
 the hitting time
$\tau_{\RoA}$ is finite a.s.,   it 
  has finite moments of every order,
and on reducing $\upmu_{\RoA} >0$ if necessary we have $\Expect_x[  \clG  ]<\infty $.
    \clThm
\end{proposition}

We have not pursued this approach in this paper because we don't yet have tools to choose an appropriate   value of $\temp>0$.  
When  $m_\RoA ( \mu \mid x) <\infty \eqdef \Expect_x[ e^{\mu  \tau_{\RoA} }  ]$ is finite then applying Markov's inequality we obtain a bound on the time required to achieve the PAC performance requirement      \eqref{e:performance-objective}.    In future research we will consider the optimization of $m_\RoA ( \mu \mid x) $   over $\temp$.

\wham{Best-state switching \rm ($\LGRM{\temp,N}$).}

The first-entry rule requires an online test for membership in
$\RoA$.  A simpler implementation is available when objective values
can be evaluated but membership in $\RoA$ cannot be tested.
This is why the best-state  switching approach $\LGRM{\temp,N}$  is highlighted in the introduction.

This variant is especially natural when the target region is the 
sublevel set $  \RoA_\eta
    =  \{x:\Obj(x)\leq \optObj+\eta\}$ introduced in \eqref{e:RoAsublevel}:    
    if any exploratory trajectory reaches $\RoA_\eta$, then the
state having the smallest objective value among all states observed
during exploration must also belong to $\RoA_\eta$.    
Consequently the
low-temperature bound on $\Tf^0$ derived below remains valid without
modification, although it may be conservative.

There is also a computational advantage.  The original
$\LG{\temp,N}$ implementation continues all $N$ trajectories
under gradient flow, whereas
$\LGRM{\temp,N}$ requires only one deterministic
trajectory after time $\Tf^0$.
Its total time is therefore
\begin{equation}
    \CompLGRM
    =
    N \Tf^0
    +
    \TGD  ,
    \label{e:CompLGRM}
\end{equation}
so that the gradient-flow contribution to the total time is
independent of $N$.

 Low temperature approximations obtained in \Cref{s:low} suggest that the temperature should scale linearly with $N$,  such as
\begin{equation}
       \tempLG
    =
    N \frac{\eta}{L} 
\label{e:temp_nom}
\end{equation}
  \Cref{t:4LGtempN} then provides an estimate for a switching time that will ensure that the PAC criterion~\eqref{e:performance-objective} is met: 
\begin{equation}
    \Tf^0  =      
       \exp\left\{  \frac{EL}{N\eta}   \right\}   =     \exp\left\{  \frac{E}{\Dbar}  \frac{1}{\tempLG}  \right\} 
    \label{e:T0LG}
\end{equation}
  where $E>\Dbar$ is fixed but arbitrary. 
These approximations require   small temperature, which imposes   constraints on $N$.    

%Hence the low-temperature theory is intended primarily to provide intuition for algorithm design.  

\subsection{A menu of temperatures}
\label{s:coupons}

In several examples we see in \Cref{s:num}  that the performance of  $\LGRM{\temp,N}$ is highly sensitive to the temperature parameter.
There is no reason to fix temperature across the $N$ trajectories in the exploration phase, and the low temperature complexity bounds are easily extended to this heterogeneous setting.

Other approaches are based on the 
 ``state-dependent Langevin diffusion'',
\begin{equation}
    dX_t
    =
    -\nabla \Obj(X_t)\,dt
    +
    \sqrt{2 \temp_{\rm eff}(\Obj(X_t)) }  \,dB_t \,,  
\label{e:LGSdiffusion}
\end{equation}
in which the state dependent temperature is defined by
\begin{equation}
    \temp_{\rm eff}(x)
    =
    \temp\,[z+\Gamma(x)_+]\,,   \quad \textit{with $z>0$.}
\label{e:LGStemp_eff}
\end{equation}

The previous approaches based on replacing the Langevin diffusion with \eqref{e:LGSdiffusion}
are  denoted $
\LGSFE{\temp,N}$ and 
$\LGSRM{\temp,N}$.
Theory surrounding low temperature asymptotics is weaker than for the standard Langevin diffusion, but enough structure remains to provide guidelines for design;  see recent work 
\cite{lelievrepavliotisrobinsantetstoltz26} on reversible diffusions that preserve the Gibbs distribution as the invariant measure.

\subsection{Complexity comparisons}
\label{s:comparisonMethods}

%\subsection{Four approaches}
%\label{s:fourL}
 
 The first two of the methods are standard, while the third is simply    $\LG{\temp,N}$ without the exploitation phase.   

 \notes{Eventually we must revise, explaining why we do not include recent approaches}

 \wham{1.} \textbf{Simulated annealing \rm (SimA).} 
A single Langevin diffusion is run with time-varying   temperature
  that vanishes according to the cooling schedule,
\begin{equation}
    \temp(t)\sim \frac{E}{\log t}    \,, \qquad \textit{with $E>\Dbar$.} 
    \label{e:coolingA}
\end{equation}

\wham{2.}
\textbf{Fixed-temperature Langevin \rm  ($\mathrm{L}_{\temp}$).}
A single Langevin diffusion is run at a constant temperature $\temp$.
\notes{See Dong.  They suggest a better bound is possible}

\wham{3.}
\textbf{Parallel Langevin \rm  ($\mathrm{L}_{\temp,N}$).}
Run $N$ independent copies of $\mathrm{L}_{\temp}$ and return the
terminal state having smallest objective value.

\wham{Overview of approximations}

 We fix  $E>\Dbar$ in the comparisons that follow.   
 
 In simulated annealing we choose the cooling schedule \eqref{e:coolingA} with this value,  where in numerical experiments we took the specific form   
\begin{equation}
    \tempSimA(t)
    =
    \frac{E}{\log(e+t)}
\label{e:cooling}
\end{equation}
Wiith $    L=\log(1/\delta)$  we fix the following  temperature parameters,
 \begin{equation}
    \temp_2
    =
    \frac{\Dbar}{E}\frac{\epsy}{L},
    \qquad
    \temp_3
    =
    N
    \frac{\Dbar}{E}\frac{\epsy}{L},
\label{e:23temp}
\end{equation}
and recall $\tempLG$ is common to any of the Langevin--gradient methods considered.  

The approximations of total simulation times, summarized in the table below, all  require  \Cref{a:jacquot}.  These approximations  are on a logarithmic scale with vanishing temperature,  as explained in \Cref{s:low}.     The proofs require fixed  $E>\Dbar$, $\epsy>0$, and the approximations become tight as $N/L \to 0$ so that $\temp\downarrow0$.

\[
\begin{array}{l|l|l}
\text{Method} & \text{Total simulation time $\clC$}
              & \text{Details}
\\[1ex]
\hline
1.\  \SimA
&
\exp(LE/\epsy)
&
\textrm{\Cref{t:SimAnn}}
\\[1ex]
2.\  \LA{\temp}
&
\exp(LE/\epsy)
&
\textrm{\Cref{t:2Ltemp}}
\\[1ex]
3.\  \LA{\temp,N}
&
N\exp(LE/(\epsy N))
&
\textrm{\Cref{t:3LtempN}}
\\[1ex]
4.\  \LG{\temp,N}
&
N\exp(LE/(\eta N))
+
O\!\left(N\log(1/\epsy)\right)
&
\textrm{\Cref{t:4LGtempN}}
\\[1ex]
5. \  \LGRM{\temp,N}
&
N\exp(LE/(\eta N))
+
O\!\left(\log(1/\epsy)\right)
&
\textrm{\Cref{t:4LGtempN}}
\end{array}
\]

Approaches~1 and~2 have essentially the same exponential dependence on the confidence and accuracy requirements.

The reduction obtained for $\LA{\temp,N}$  follows from the fact that   increasing $N$ allows each trajectory
to be run at a warmer temperature, which reduces the required single-run simulation time $T_3$.   
The improvement is most evident with appropriate choice of $N$.

 \wham{Restart and temperature}

On choosing    $N$ as a function of $\delta$, while ensuring that the low temperature constraint is met,  we conclude that 
the use of a large number of independent runs reduces complexity dramatically when 
compared to either simulated annealing or fixed-temperature Langevin.    

The proposition that follows provides two examples.   Observe that the approximation \eqref{e:betterC3}  is very close to the complexity of SGD with \textit{convex} objective (recall  \eqref{e:SGDcomplexity}).

\begin{proposition}
\label[proposition]{t:BetterBddsThroughRestart}  
Using Langevin methods with restart, the value of $N$ can be chosen to grow appropriately with $L$ to achieve  
\eqref{e:CompLG} for   $  \LGRM{\temp,N}$, and with   $ \LA{\temp,N}$ the approximation
\begin{equation}
 \clC^{\LA{\temp,N}} \approx  E L^{1+o(1)}   \frac{1}{\epsy} 
\label{e:betterC3}
\end{equation} 
\end{proposition}

\PF
To obtain \eqref{e:CompLG} we choose     $N = \lceil EL/a_L \rceil$,  with $a_L\to\infty$ as $L\to\infty$, and satisfying  $ a_L=o(\log L)$  (for example $a_L = \log\log(L)$, defined for $\delta \le  \exp(- e)$).  On choosing $\temp(L) = a \tempLG$, for fixed $a>0$ and $ \tempLG$ defined in  \eqref{e:temp_nom},
we conclude that the low temperature assumption is met: 
 $\temp(L) \to 0$ as $\delta\downarrow 0$, and moreover we obtain \eqref{e:CompLG}:
 \begin{equation}
\begin{aligned}
    \CompLGRM
    &\approx
    \frac{EL}{a_L}
    \exp \Big\{\frac{1}{\eta}a_L\Big\}
    +\TGD(\epsy) \\
    &=E  L^{1+o(1)}
    +O\Big(\log\frac1\epsy\Big).
\end{aligned}
\label{e:CompLGb}
\end{equation}
 
 \notes{return to this argument since I think $a>1$ will be needed}

In the case of  $ \LA{\temp,N}$  we take  
$    N_L
    =
    \left\lfloor
         {E L}/(\epsy a_L)
    \right\rfloor $, maintaining the  same constraints on $a_L$.    
    We obtain  $N_L/L  \sim 1/(\epsy a_L)\to 0$ as $\delta\downarrow 0$, and   
    the low-temperature approximations remain applicable: 
\[
    \temp_3
    =
    \frac{E^\bullet}{E}
    \frac{N_L\epsy}{L}   \sim   E^\bullet \frac{1}{a_L} 
    \longrightarrow0.
\]

Moreover,  $ 
      {LE}/(\epsy N_L) 
    \sim
      a_L$ 
and hence from the table, the total simulation time   satisfies
\begin{equation*}
    \clC^{\LA{\temp,N}}
\approx
N_L\exp(LE/(\epsy N_L)) 
    =
    \frac{EL}{\epsy a_L}
    \exp\{  a_L+o(a_L)\}     =
    \frac{EL}{\epsy  }
    \exp\{  a_L+o(a_L)\}.
%\label{e:C3-slow}
\end{equation*}
Under the assumption  $ a_L=o(\log L)$  we have
$
    \exp\{  a_L+o(a_L)\}
    =
    L^{o(1)}
$,  establishing \eqref{e:betterC3}.  
\qed

%%%%%%%%%%%%%%%%%%%%%%%%%%%%%%%%%%%%%%%%%%%%%%

\subsection{Prior knowledge required for implementation}
\label{s:prior-knowledge}

An important question is how much prior information about the objective
landscape is required by the different schemes.  We distinguish between
quantities used to obtain sufficient conditions for a finite-time PAC
guarantee and information that must be available online while the
algorithm is running.  Structural assumptions used in the analysis are
not, by themselves, prior information required by the algorithm.

 \wham{Simulated annealing.}
Classical simulated annealing already requires nontrivial global
information.  In the logarithmic cooling schedule \eqref{e:cooling},
the constant $E$ must exceed the critical depth associated with
simulated annealing; in many settings this coincides with the dominant
energy barrier $\Dbar$ used here.  Thus exact knowledge of the energy
landscape is not required, but a sufficiently large upper bound on this
global energy scale is already needed in classical simulated annealing.

 \begin{table}[t]
\centering
\small
\begin{tabular}{p{0.20\textwidth}|p{0.36\textwidth}|p{0.35\textwidth}}
\textbf{Method}
&\textbf{Prior information used in sufficient PAC conditions}&
\textbf{Information needed online}
\\ \hline

$\SimA$
&
A sufficiently large upper bound on the dominant energy barrier
$\Dbar$ to select the cooling constant.
&
Evaluation of $\Obj$ and $\nabla\Obj$, together with the prescribed
cooling schedule.
\\[1ex]

$\LA{\temp}$
&
An upper bound on $\Dbar$ for conservative selection of temperature
and simulation time.
&
Evaluation of $\Obj$ and $\nabla\Obj$.
\\[1ex]

$\LA{\temp,N}$
&
As for $\LA{\temp}$, together with the chosen number $N$ of
independent trajectories.
&
Evaluation of $\Obj$ and $\nabla\Obj$.   
\\[1ex]

$\LG{\temp,N}$
&
An upper bound on $\Dbar$, a lower bound on the attraction margin
$\eta$, and conservative local convergence constants
$c_{\RoA},\upmu_{\RoA}$.
&
Evaluation of $\Obj$ and $\nabla\Obj$; the deterministic switching
time is fixed in advance, so no online test for membership in
$\RoA$ is required.
\\[1ex]

$\LGFE{\temp,N}$
&
For a finite-time PAC guarantee, quantitative information controlling
the hitting time of the target region, together with local convergence
constants.
&
Evaluation of $\Obj$ and $\nabla\Obj$, together with an online test
for entry into the switching region.
\\[1ex]

$\LGRM{\temp,N}$
&
An upper bound on $\Dbar$, a lower bound on the attraction margin
$\eta$, and conservative local convergence constants.
&
Evaluation of $\Obj$ and $\nabla\Obj$, together with storage of the
state having lowest observed objective value.  No online test for
membership in $\RoA$ is required.
\end{tabular}

\caption{Prior information required by the Langevin-based schemes.
The second column lists quantities used in sufficient conditions for a
finite-time PAC guarantee; the third lists information required to run the algorithm.   }
\label{t:prior-knowledge}
\end{table}

 \wham{Fixed-temperature approaches.}   These require   no additional
information.  In particular, in the pathwise-minimum variant $\mathrm{L}_{\temp,N}^{\mathrm{(a)}}$, the approach proceeds by retaining 
the lowest objective value observed
during the simulations.   The upper bound $E$ is only required to establish  a bound on the needed simulation time to meet the PAC requirement.

 \wham{Langevin--gradient approaches.} 
 These  require   additional   information:
the energy margin $\eta$, target set $\RoA$, and constants $c_{\RoA}$, $\upmu_{\RoA}$ appearing in  \Cref{a:jacquot}.    
Upper bounds on $\Dbar$, $c_{\RoA}$, and  lower bounds on  $\eta$,  $\upmu_{\RoA}$    are sufficient to obtain bounds on the required total simulation times.

The amount of this information that must be available \emph{online}
depends strongly on the variant.  The first-entry scheme
$\LGFE{\temp,N}$ requires an online test for
membership in the switching region $\RoA$.  When the region is a sublevel
set, such a test may require knowledge of both $\eta$ \textit{and} $\optObj$. 
%At the opposite extreme,
%$\LGFE{\temp,1}$ with an unbounded exploration
%horizon requires no information about $\Dbar$ to obtain almost-sure
%global convergence at any fixed positive temperature.  What is lost
%in this case is not convergence but a useful finite-time PAC
%complexity estimate.  Obtaining such an estimate requires quantitative
%information about the hitting time of the target region, which
%motivates topics for future research discussed in the
%conclusion.

By contrast,
$\LGRM{\temp,N}$ requires no such membership
test.  It records only objective values during exploration, retains
the best state observed, and initializes gradient flow from that
state.  The region $\RoA$ and the constant $\eta$ enter the
\emph{analysis} to obtain sufficient conditions for the PAC guarantee, but need not be
identified by the algorithm while it is running.

The distinction between the information used to establish sufficient
conditions for a finite-time PAC guarantee and the information needed
to implement the algorithm is  summarized in  \Cref{t:prior-knowledge}.
Observe that   the improvement
from restart is not derived from greater knowledge
of the global landscape.  In particular,
$\LA{\temp,N}$ and its pathwise variant $\LGRM{\temp,N}$ require essentially the same
global information as classical simulated annealing.

%%%%%%%%%%%%%%%%%%%%%%%%%%%%%
\subsection{Challenges in high dimensions}
\label{s:highd}

The basic phenomenon considered in this subsection is classical in
random search and multi-start optimization.  If independent trials have
a small probability of entering a near-optimal region, then the
performance of the best of $N$ trials is governed by an order statistic;
in high dimension the probability mass of such regions may decrease
very rapidly.  Closely related consequences for random search and
multi-start methods are developed in
\cite{PepelyshevZhigljavskyZilinskas2018,
NoonanZhigljavsky2024,MasuyamaDanUmetani2025}.

%see also \cite{mus97,SchutteHaftkaFregly2007} for the role of restart in global optimization.

Our purpose here is to examine this phenomenon for Langevin restart,
where the distribution of each exploratory state is approximately
Gibbs and, crucially, its temperature is itself a design parameter.
This leads to a distinction between fixed-temperature restart and
accuracy-dependent temperature selection that is not visible from the
classical order-statistic argument alone.

We consider here an abstraction of the restart approaches.

Let  $\{ X_i : i\ge 1 \}$  be an i.i.d.\ sequence, and denote  $G_i = \Obj(X_i)  - \optObj$.
It is assumed that the common mean of $G_i$ is finite and the common distribution   has density $f$ and  CDF $F$ satisfying    $F(x)>0$ for all $x>0$.     Denote
\[
M_n \eqdef  \min\{G_1,\dots,G_n\}.
\]
We have by independence,  
\begin{equation}
\begin{aligned}
\Prob(M_n > \epsy) & = (1 - F(\epsy))^n  \,, \quad \epsy>0
\\
\Expect[M_n] & = \int_0^\infty (1 - F(t))^n\,dt.
\end{aligned}
\label{e:MnQandMean}
\end{equation}
In particular, the exceedance probability $\Prob(M_n > \epsy) $ vanishes at a geometric rate as $n\to\infty$.

  Adopting the notation of   \eqref{e:bestOfN}, write
\[
    \haX_n    \in    \argmin \left\{ \Obj( X_k) :  1\le k\le n   \right \}     
\]
from which we obtain the representation
\[
     \Prob\{\Obj(\haX_n)-\optObj>\epsy\} =\Prob\{M_n >\epsy \}  
\] 
Consequently, this provides an abstraction of the parallel Langevin
$\mathrm{L}_{\temp,N}$ approach when the $X_i$ have common Gibbs
distribution.  Based on the exceedance probability in
\eqref{e:MnQandMean}, it is straightforward to identify the
$(\epsy,\delta)$-PAC sample complexity for this idealized model:
\begin{equation}
    n^\star(\epsy,\delta)
    \eqdef
    \min\left\{
        n:
        \Prob\{\Obj(\haX_n)-\optObj>\epsy\}\leq\delta
    \right\}
    =
    \min\left\{
        n:
        \Prob\{M_n>\epsy\}\leq\delta
    \right\}.
\label{e:nstarIdeal}
\end{equation}

While the rate of convergence of the exceedance probability is geometric, the specific exponent in this convergence rate is 
 determined by   properties of $F$ near the origin.  If
  $F(x) \sim c x^\alpha $  as $x \downarrow 0$ we obtain,  for any $r>0$,
  \begin{equation}
  \begin{aligned} 
    \Prob\{\Obj(\haX_n)-\optObj> r n^{-1/\alpha}  \} 
&  =
  \Prob\{M_n> r n^{-1/\alpha}\} 
  \\
  &=
 [1-F(r n^{-1/\alpha})]^n  \to \exp(-cr^\alpha)\,, \quad n\to\infty
 \end{aligned}
\label{e:BadMinTails}
\end{equation}
A large value of $\alpha$ can be expected in large dimensions,  implying that convergence may be very slow.
This is seen in the   ideal quadratic objective   $\Obj(x) = \| x\|^2$ for which we obtain  $\alpha =d/2$ whenever the distribution of $X_i$ has a continuous density that is positive at the origin, as   is the case for the associated Gibbs distribution (which is Gaussian for this objective).

Below is a general statement for an objective satisfying  \Cref{a:jacquot} and a few additional assumptions.  
In the following we use $\Gammafn(\cdot)$ to denote the Gamma function, to distinguish
it from the objective $\Obj$. 
 
\begin{proposition}
\label[proposition]{t:RunMax_dimension}
Suppose that \Cref{a:jacquot} holds.  Assume moreover that there is  a unique minimizer $x^*$ satisfying $\Obj(x^*) =0$ and $H=\nabla^2\Obj(x^\star)>0$.   

Let $\{ X_i : i\ge 1 \}$ be i.i.d.,  whose distribution is Gibbs with parameter $\temp>0$,  and let $G_i = \Obj(X_i)$.   Then the common CDF of $G_i$ is
\[
    F_\temp(x)
 =
    \frac{1}{Z_\temp}
    \int_{A_x}  
    \exp\left\{ 
        - \temp^{-1}  \Obj(y)
    \right\}dy
\]  
where $A_x = \{y: \Obj(y)\leq x\}$.   Consequently, for  each fixed $\temp>0$ there is $c_\temp$ such that 
\[
    F_\temp(x)
    \sim
    c_\temp x^{d/2},
    \qquad x\downarrow0.
\]
Moreover,
\[
    c_\temp
    \sim
    \frac{1}
    {\Gammafn(d/2+1)\temp^{d/2}},
    \qquad \temp\downarrow0.
\]
\clThm
\end{proposition}

 The proposition combined with \eqref{e:BadMinTails}
 implies a curse of dimensionality:  while $M_n$ converges to zero as $n\to\infty$, the rate of convergence is 
 $O( n^{-2/d} ) $.   The value of low temperature is only in the increased value of $c$ 
in \eqref{e:BadMinTails}:  under the assumptions of \Cref{t:RunMax_dimension} we have for large $n$,  
\[
  \Prob\{n^{2/d}M_n>r \}  \approx   \exp(- c_\temp r^{d/2}) 
 \]
 In the following we find that the mean of $M_n$ also vanishes at rate  $O( n^{-2/d} ) $.
 
 %with constant   proportional to $c_\temp^{-2/d}$. 

\begin{proposition}[Limit law for the running minimum]
\label[proposition]{t:aStatsRunningMin}
Suppose that $F(x) \sim c x^\alpha$ as $x \downarrow 0$. 
Then, 
\wham{(i)}
the stochastic process $\{ \clM_n = n^{1/\alpha} M_n : n\ge 1\}$ converges in distribution to a random variable $\clM$
with CDF $F_\clM (t) = 1 -  \exp(-c t^\alpha)$:
for any $t > 0$,
\begin{equation}
 \lim_{n\to\infty} 
\Prob(\clM_n > t) =\Prob(\clM > t) = 1-F_\clM(t)
\label{e:LimitLaw_b}
\end{equation} 
\wham{(ii)} If the common mean of $\{G_k \}$ is finite then  
\[
 \lim_{n\to\infty}  n^{1/\alpha}  \Expect[M_n]  =  \Expect[\clM] = \int_0^\infty \exp( -c t^\alpha)\, dt
 = \Gammafn(1/\alpha)  \frac{1}{\alpha} \frac{1}{ c^{1/\alpha} }
\]
 
\wham{(iii)}  Suppose that $\Expect[\exp(\theta_0 G_k )] <\infty$ for some $\theta_0>0$ and any $k$.
Then for $\theta<\theta_0$, 
\[
\Expect[\exp(\theta M_n) ] = 1 + \theta \Expect[\clM]  n^{-1/\alpha}  + o(n^{-1/\alpha} )
\]
 \qed
 \end{proposition}

 \wham{Reconciling two approximations}

To achieve the error bound $    \Prob\{\Obj(\haX_n)-\optObj>\epsy\}\leq\delta$,  the limit  \eqref{e:BadMinTails}
tells us that for large $n$ we require the constraints    
\[
\begin{aligned}
\exp(-cr^\alpha) \le  \delta 
& \Longleftrightarrow  r \ge (L/c)^{1/\alpha}
\\
r n^{-1/\alpha}\le \epsy &  \Longleftrightarrow  n  \ge (r/ \epsy )^\alpha 
\end{aligned}
\]
where $L=\log(1/\delta)$.    Applying \Cref{t:RunMax_dimension} to substitute $\alpha =d/2$, 
  the required value of $n$ is roughly
\begin{equation}
n(\delta,\epsy) = \frac{L}{c} \frac{1}{\epsy^{d/2}} \,,
\label{e:BADn}
\end{equation}
which appears to contradict the approximation $ \clC^{\LA{\temp,N}} \approx L /{\epsy} 
$  given in \eqref{e:betterC3}.    
This apparent contradiction is resolved if the distribution of $X_i$
is allowed to depend on the accuracy requirement $\epsy$.
 
 We denote    \notes{use a different symbol than $\kappa$ perhaps}
\begin{equation}
K_d(\kappa)
    \eqdef
    \Gammafn\left(\frac d2+1\right)\kappa^{d/2},
\label{e:Kd_kappa}
\end{equation}    
and recall $    \RoA_\eta
    =
    \{x:\Obj(x) \le \optObj+\eta\} $,  $    L=\log(1/\delta)$.
The proof of \Cref{t:running-min-temp-scaling} is postponed to the Appendix.

\begin{proposition}[Accuracy-dependent temperature]
\label[proposition]{t:running-min-temp-scaling}
Suppose that   \Cref{a:jacquot} holds 
 with
$d>2$, and moreover,
\whamb
 $\Obj$ has a unique global minimizer
$x^\star$   that it is   non-degenerate:
$
    H
    \eqdef
    \nabla^2\Obj(x^\star)>0$.    
    
    The objective is normalized so that  $\Obj(x^\star)=\optObj=0$.
    
\whamb
The common distribution of $X_i$ is $\pi_{\temp_{\epsy}}$,  the Gibbs distribution with temperature $\temp_{\epsy}$, where for fixed 
 $\kappa>0$  
\begin{equation}
    \temp_{\epsy}
    =
    \kappa\,\epsy^{\,1-2/d}
    \label{e:rho-epsilon-scaling}
\end{equation}
Then, 
\begin{align}
    \pi_{\temp_{\epsy}}
    \{A_\epsy\}
 &   \sim
    \frac{\epsy}{     K_d(\kappa) } \,, 
    \qquad
    \epsy\downarrow0.
    \label{e:Gibbs-target-temp-scaled}
\\
    n^\star(\epsy,\delta)
 &   \sim
    K_d(\kappa)\frac{L}{\epsy},
    \label{e:n-epsilon-delta}
\end{align}
where in     \eqref{e:n-epsilon-delta}
 the asymptotic relation means that
\[
    \frac{n^\star(\epsy,\delta)}
    {K_d(\kappa)L/\epsy}
    \longrightarrow 1
\]
for any joint limit of $\epsy,\delta$ satisfying
$    \epsy\downarrow0$ and $L/\epsy \to \infty$.   
\clThm
\end{proposition}

Taking $\kappa=1$ gives the particularly simple choice
$
    \temp_{\epsy}
    =
    \epsy^{\,1-2/d}$,
for which
\begin{equation}
    n^\star(\epsy,\delta)
    \sim
    K_d
    \frac{L}{\epsy}.
    \label{e:n-epsilon-simple}
\end{equation}
with $K_d\eqdef K_d (\kappa)$ with $\kappa=1$.   
Thus, for every fixed $d>2$, adaptation of the temperature replaces
the fixed-temperature dependence $\epsy^{-d/2}$
seen in \eqref{e:BADn}
 by the much more
favorable dependence $\epsy^{-1}$, consistent with \eqref{e:betterC3}.

Dimension enters strongly through
the prefactor.  Stirling's approximation gives
\[
    K_d
    =
    \Gammafn\left(\frac d2+1\right)
    \sim
    \sqrt{\pi d}
    \left(\frac{d}{2e}\right)^{d/2}.
\]
This does not constitute a uniform-in-$d$ complexity estimate:
\Cref{t:running-min-temp-scaling} is an asymptotic result for each fixed
dimension.  It nevertheless illustrates how dimension-dependent
prefactors can be hidden by logarithmic approximations.

There is one important limitation to this calculation.  The quantity
$n^\star$ counts independent Gibbs samples and assigns no cost to
producing them.  If temperature could be reduced without cost, then
there would be no meaningful optimization over $\temp$: sufficiently
small temperature would make a single Gibbs sample successful with
arbitrarily high probability.

For Langevin diffusion, however, low temperature generally increases the time required to generate an
approximately stationary sample.  Estimation of the required simulation time is based on mixing-time approximations surveyed \Cref{s:low}.  
The conclusions of this subsection are   therefore complementary to the
low-temperature theory developed in that section.

\section{Numerical Results}
\label{s:num}

In implementation, the various approaches become recursive algorithms following approximation through an Euler scheme or some other discrete-time approximation to the diffusion and gradient flow.  

\subsection{Building blocks for numerical experiments}   

Before we begin a survey of numerical experiments we require guidelines for parameter selection.
Given our goal is to achieve a PAC constraint with low computational effort, we begin with a way of fairly comparing this quantity for the very different algorithms.  

\subsubsection{Computational work}
To compare the algorithms on a common computational scale, we use a
calibrated work model.  Separate MATLAB microbenchmarks were used to
estimate the cost of the primitive operations appearing in the
algorithms---for example, evaluating the objective and gradient, and
generating a Gaussian random vector.  ``Work'' is the sum of these
measured primitive costs over the operations performed by an algorithm.
It is expressed in seconds, but should not be confused with the total running time $\clC$ used in the previous section.  
wall-clock running time of the vectorized MATLAB implementation.

See \Cref{s:WorkModel} for a table containing explanations of all terms and their  numerical values. 
In particular, $    c_{\Obj,\nabla}$ is the cost per iteration of simultaneously computing the objective value and its gradient at a selected point.

Once the stepsize $\uptau$ has been specified,  simulation of the diffusion up to time $\Tf$ corresponds to $\nmax$ iterations, with
$\nmax \eqdef     \left\lceil  {\Tf}/{\uptau}\right\rceil$

The computational work for simulated annealing after $n$ Euler iterations  is
\[
    \WSimA_n
    =
    c_{\Obj,\nabla}
    +n\bigl(c_{\Obj,\nabla}+c_Z+c_{\rm cmp}\bigr) \, .
\]

In the Langevin-gradient method we have $\Tf = \Tf^0 + \TGD  $.     Computation of work is more involved in part because   there may be two stepsizes:  $\uptau$ for the diffusion, and $\uptauGD$ for the gradient descent phase.   
In this case the number of iterations in the respective phases are
 \begin{equation}
 n^0_{\max}
    =
    \left\lceil  {\Tf^0}/{\uptau}\right\rceil
    \qquad
    n_{\rm GD}
    \eqdef
    \left\lceil
        \frac{\TGD  }{\uptau_{\rm GD}}
    \right\rceil.
    \label{eq:n0nGD}
\end{equation}
Ignoring the  computational work in the recursion used   to select the best state among
the $N$ trajectories, the total computational work for $\LGRM{\temp,N}$  has the form
\begin{equation}
    \WLGRM
    =
    N w^0+\wGD ,
    \label{e:WorkTempScale}
\end{equation}
where $ \wGD =c_\nabla \nGD $ and
$
    w^0
    =
    c_{\Obj,\nabla}(K_a+1)
    +(c_Z+c_{\rm cmp}) n_{\max}
$.

\subsubsection{Choice of parameters}

Although $\LGRM{\temp,N}$ does not require the parameter $\eta>0$ for implementation (see main assumption in paper),  to test the theory in numerical experiments we consider $\tempLG$ in \eqref{e:temp_nom}
 as a nominal temperature value and use     \eqref{e:T0LG} to define the nominal exploration time with $E= 2\Dbar$, so that 
\begin{equation}
    \tempLG   =    N \frac{\eta}{L} 
\,, 
    \qquad 
    \Tf^0  =      
       \exp\left\{  \frac{EL}{N\eta}   \right\}   =     \exp\left\{  2 \frac{1}{\tempLG}  \right\} 
\label{e:nomTempAndTim}
\end{equation}
     This choice of $E$ gives accurate predictions in some of the
low-dimensional experiments.  However, the simulated-annealing
condition appears only to require that $E$ exceed the relevant
critical depth; it does not prescribe the factor $2$.  For example,
the equally admissible choice $E=4\Dbar$ would divide the predicted
temperature by two and would give a much less satisfactory prediction.

In testing sensitivity of performance to the value of this switching time we introduce a scaling parameter $a>0$ and define 
\begin{equation}
    \Tf^0(a) =a\Tf^0    \,,\qquad a>0 
\label{e:T4a}
\end{equation}
In the definition of the random switching time   $\sigma_4^0$ we  minimize $\Obj(X_t^i) $ over $t\le   \Tf^0(a) $ for each $1\le i\le N$.

The intended use of the theory is not to prescribe a universal
restart count  $N$ or temperature $\temp$.  Rather, an estimate of the local accuracy
scale $\eta$   informs the amount of parallel exploration $N$ as well as the value of $\temp$.   
 Numerical experiments test the sensitivity of the
algorithm to these theoretically motivated design choices.

\subsection{Six hump camel function}  
 
This is the objective function on $\Re^2$,  
\begin{equation}
\Obj(x_1,x_2) = \left(4-2.1x_1^2+\tfrac{1}{3} x_1^4 \right)x_1^2 +x_1x_2 +\left(-4+4x_2^2\right)x_2^2 \,,
\label{e:sixHumps}
\end{equation}
with two global minima at approximately  $(0.0898,-0.7126)$,  $ (-0.0898,0.7126)$,  giving $\optObj \approx-1.0316$.   
The dominant energy
barrier is approximately $ \Dbar=1.032$.

This objective function violates \Cref{a:jacquot}   since $\nabla \Obj(x)$ is not globally Lipschitz continous.    The theory of this paper applies on modifying the objective outside of a large compact set so it remains $C^2$ and the remaining required assumptions hold.

Common parameters in   simulations are
$\epsy=10^{-5}$,  $\delta=0.05$,  $\uptau=10^{-4}$ and $\uptauGD=10^{-3}$.

 \subsubsection{Comparison with simulated annealing}

Two approaches are compared here:  Simulated annealing using the cooling schedule \eqref{e:cooling} with 
$
E=2     >\Dbar$,  and a version of $\LGRM{\temp,N}$  with parameters 
$\temp = 1 $ and a small value of $\Tf^0$, 
 chosen so that  the use of ``best-state'' in this method is clearly illustrated.   

 \Cref{f:LG24temp1vsSimA} provides a comparison of the two approaches,  in which the lower panels show  the excess
cost as a function of time.  Observe that the time horizon is $t\le 2$ in the case of $\LGRM{\temp,N}$ and $t\le 50$ for simulated annealing.

The upper-left panel shows
the exploratory trajectory that produced the retained best state for $\LGRM{\temp,N}$,
together with the subsequent gradient-descent trajectory.  The retained
state is observed at time $\sigmaRMswitch$, while exploration continues until
the prescribed horizon $\Tf$; gradient descent is then initialized from
the stored state.  The exponential convergence for $t>\sigmaRMswitch$ is evident in the lower left panel.

\begin{figure}[h!]
	\centering
	\includegraphics[width=0.5\hsize]{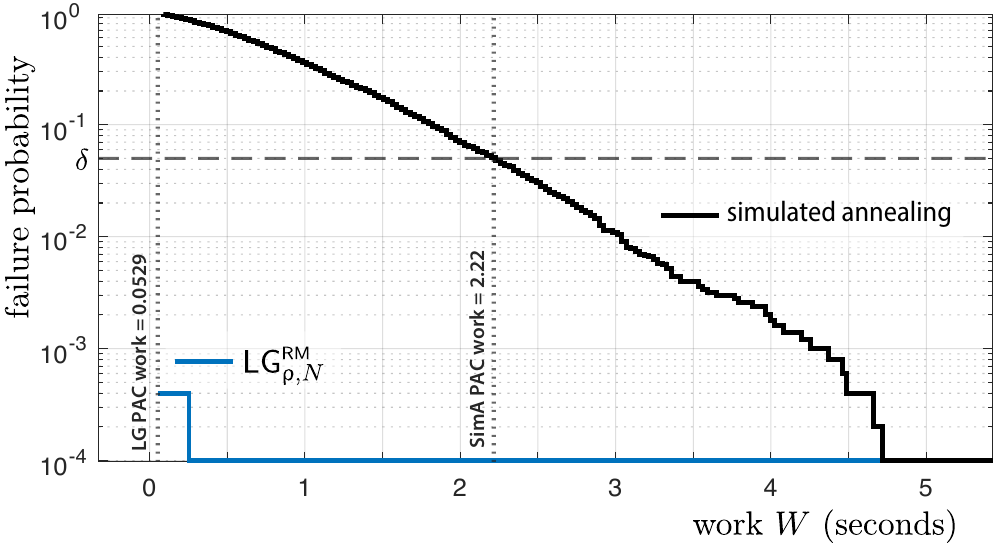}  
	\caption{Failure probability as a function of work for the six-hump camel objective 
	using $\LGRM{\temp,N}$  with    $\temp=24.5$  and $N=24$.}
			\label{f:PAC1}
\end{figure}  
\Cref{f:PAC1} shows empirical failure probability, based on 2000 independent trials, as a function of work for $\LGRM{\temp,N}$ using $\temp=24.5$ and a range of $\Tf^0$.     The resulting PAC work estimates are $0.0529$ seconds for $\LGRM{\temp,N}$ and $2.22$ seconds for simulated annealing.

This temperature value was selected as the best among a large grid.

\begin{figure}[h!]
	\centering
	\includegraphics[width=0.8\hsize]{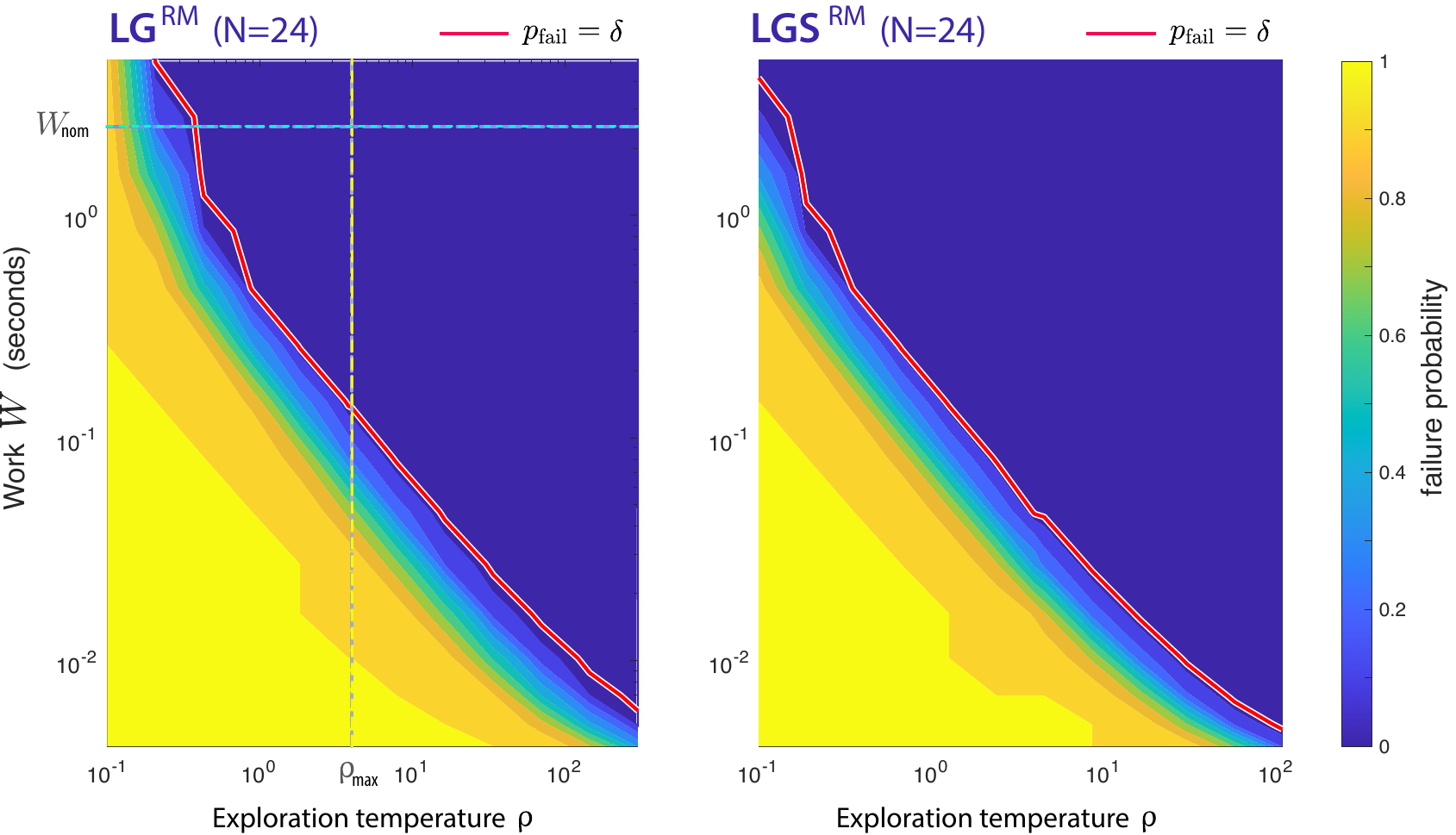} %WorkVsProbFailureCamel_N24.pdf}  
	\caption{Empirical performance of two best-state Langevin--gradient methods  for the modified six-hump camel objective.  The left hand side shows results for  $\LGRM{\temp,N}$ and the right hand side  $\LGSRM{\temp,N}$, each using $N=24$.}
		\label{f:WorkVsProbFailureCamel_N24}
\end{figure}

A larger study was done to investigate the sensitivity of performance to temperature in $\LGRM{\temp,N}$.  
With $N=24$ fixed,   
   based on a grid of $(W,\temp)$ values, and 2000 independent experiments, 
the empirical  failure probability   was obtained, denoted
\[
    \hapfail (\temp,W)
    =
    \widehat{\Pr}\{
        \Gamma(\widehat X)-\Gamma^\star>\epsy
    \} \, ,
\]
The contour plot  on the left hand side of \Cref{f:WorkVsProbFailureCamel_N24} shows its dependency on these two parameters.
The vertical dashed line marks
the low-temperature design value
$\tempLG\approx4.00$ obtained from~\eqref{e:nomTempAndTim}.

Recall that \(\LGSRM{\temp,N}\) is the version of best-state switching in which the   constant exploration
temperature $\temp$ is replaced by \eqref{e:LGStemp_eff}:
\begin{equation}
    \temp_{\rm eff}(x)
    =
    \temp\,[z+\Gamma(x)_+]\,,  \quad \textit{with $z> 0$.    
}
\end{equation}
The contour plot  on the right hand side of \Cref{f:WorkVsProbFailureCamel_N24} shows   results using this approach, with $z=1$.  
The minimal work required to achieve the PAC constraint is reduced slightly.

 \subsubsection{Spectral analysis}

\begin{figure}[h!]
	\centering
	\includegraphics[width=0.5\hsize]{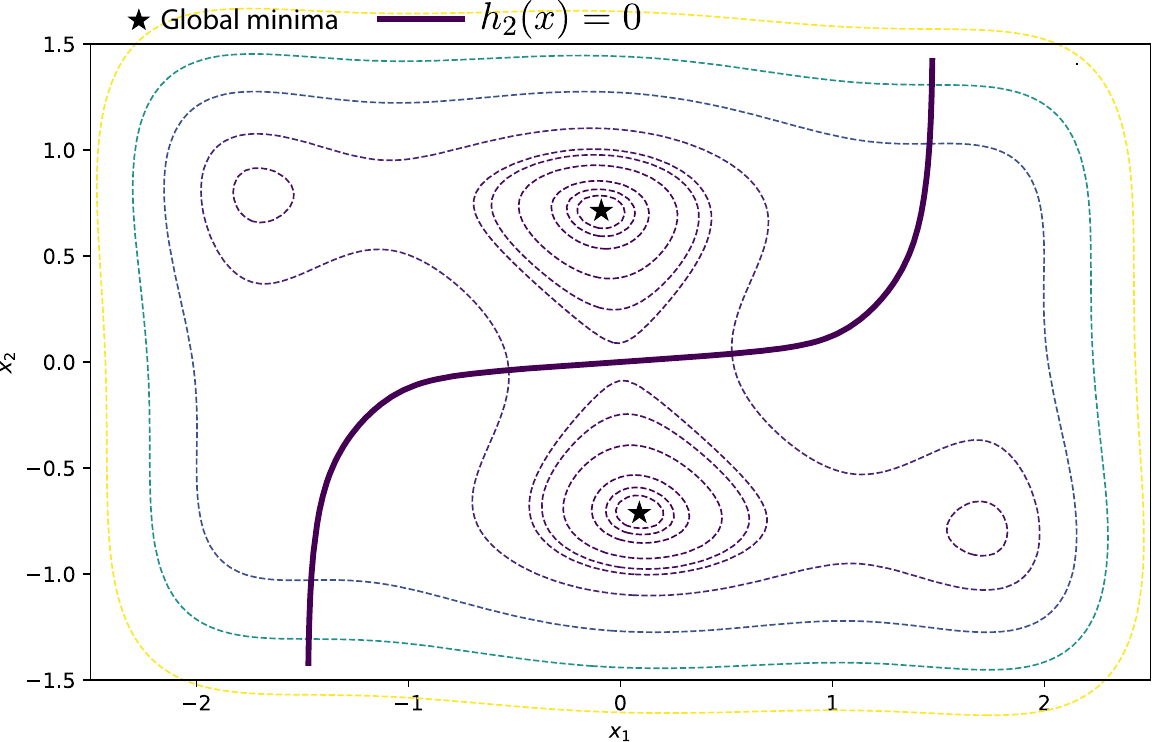}  
	\caption{Contour plots of the six hump camel function \eqref{e:sixHumps} and the curve $\{x:h_2(x) = 0\}$ with $h_2$ the second eigenfunction of the differential generator when    $\temp=0.2$.}
		\label{f:camel_h2_rho_0p2}
\end{figure}

Contour plots are shown in \Cref{f:camel_h2_rho_0p2},  along with the curve $\{x:h_2(x) = 0\}$ with $h_2$ the second eigenfunction of the differential generator with $\temp=0.2$.   Computations reveal that   the second eigenvalue $\lambda_2=-\gap_{\!\temp}$ is not repeated, so that $|\lambda_i | > \gap_{\!\temp}$ for $i \ge 3$.  The significance of this zero set and these observations regarding the spectrum 
 is the topic of \cite{huimeysch04a}, from which we obtain the following corollary.

Recall that $\gap_{\!\temp}$ denotes the spectral gap of the differential generator. 
Let $A$ be a closed set with non-empty interior  and denote  
\[
\tau_A = \min\{t\ge 0 :  X_t\in   A \}  \quad
\textit{and}
\quad
m_A( \vartheta \mid x) =     \Expect_x[ \exp( \vartheta \tau_A)]   \,, \ \vartheta>0 \, .
\]
 Let $H_+ =\{ x : h_2(x)>0 \}$ and  $H_- =\{ x : h_2(x)<0 \}$
\wham{Communication across a zero set.}
Suppose that   $A\subset H_+$.
Then for any $\temp>0$   we have in this example $ |\lambda_3| >\gap_{\!\temp}$, and
  \begin{equation}
m_A( \vartheta \mid x) 
\begin{cases}
<\infty & \text{$x\in H_+$ for any $\vartheta < |\lambda_3|$.}
\\
=\infty &  \text{$x\in H_-$ for any $\vartheta\ge   \gap_{\!\temp}$.}
\end{cases}
\label{e:MGFhuimeysch04a}
\end{equation}
In this sense, the zero set portrayed in \Cref{f:camel_h2_rho_0p2}
 is a stochastic barrier between $H_+$ and $H_-$.      
\clThm

We tested $\LGFE{\temp,N}$ with the stopping rule  \eqref{e:BasinStop} and level  $\RoA= \{x : \Obj(x ) \le \optObj + \eta\}$, using  $\eta = 0.5$ for which \Cref{a:jacquot}~(iii) is satisfied.   Results are summarized in two cases, where in each experiment simulations were conducted with common initial condition $X_0 = (2.0; 1.5 )$ (fixed for ease of comparison).
 Also, in each case the dominant eigenvalues and hence the spectral gap were approximated numerically.   

\begin{figure}[h!]
	\centering
	\includegraphics[width=\hsize]{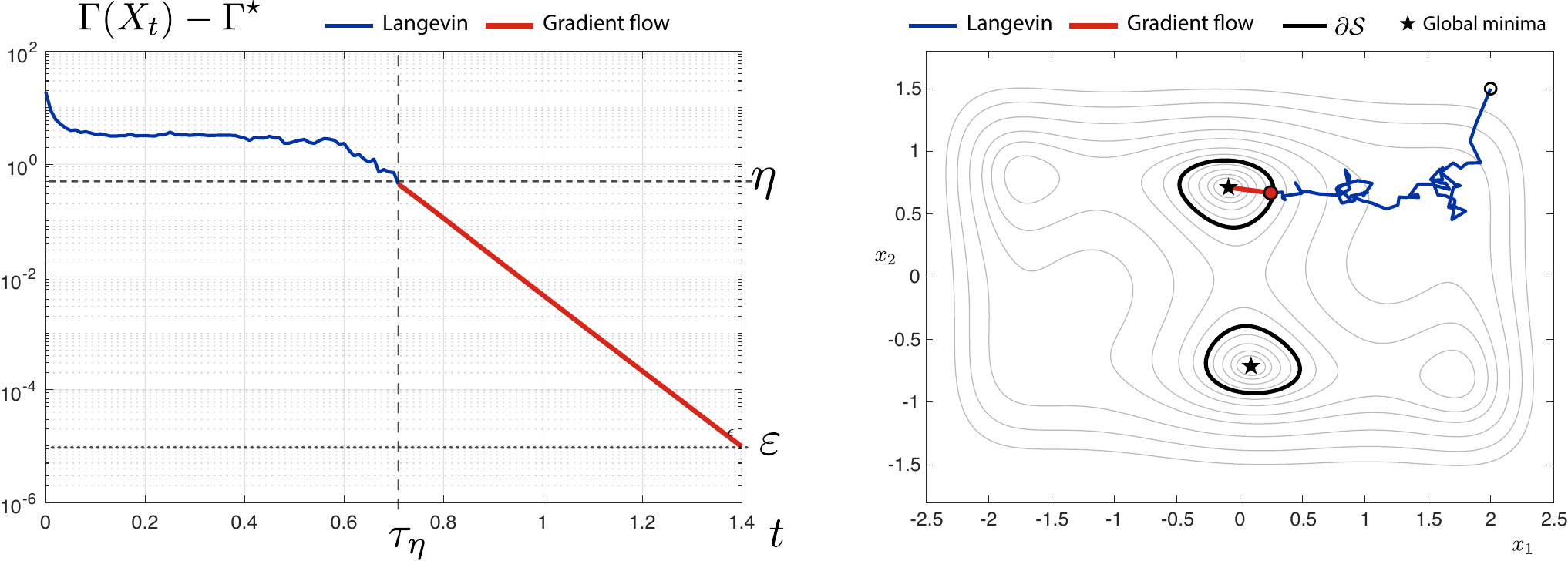}  
	\caption{$\LGFE{\temp,1}$
	for	the six hump camel function \eqref{e:sixHumps}.
Shown on the left is a plot of the excess cost with $\temp=0.2$.
	  The time $\tau_\eta$ indicated in the figure is    
	\eqref{e:BasinStop} using $\eta=0.5$.   
	%%%
	Shown on the right is the associated  trajectory of $X_t$,   	along with level sets of the objective:
	The outlines of two convex regions shown in the figure indicate the boundaries of the two connected components of $\RoA$.
}
		\label{f:Cost+PhaseCamelLG}
\end{figure} 

\wham{Low temperature:}  With $\temp=0.2$  the spectral gap is $\gap_{\!\temp} =|\lambda_2| \approx 0.01858$, and the
 first four eigenvalues are $\{\lambda_i \} = \{ 0, -0.01858,  -0.04184,  - 0.04217\}$. 
 Hence the ratio of the two limits for $r$ in    \eqref{e:MGFhuimeysch04a} is  
$\kappa = \lambda_3/\lambda_2  
\approx  2.25$.

%\begin{figure}[h]
%	\centering
%	\includegraphics[width=0.5\hsize]{CostCamelLG.pdf}  
%	\caption{Plot of the excess cost for a single run of $\LGFE{\temp,N}$  with $\temp=0.2$
%	 using the switching rule \eqref{e:BasinStop}
%for	the six hump camel function \eqref{e:sixHumps}.
%	  The time $\tau_\eta$ indicated in the figure is    
%	\eqref{e:BasinStop} using $\eta=0.5$.
%}
%		\label{f:CostCamelLG}
%\end{figure} 

Results from one run are illustrated in \Cref{f:Cost+PhaseCamelLG}.  The plot on the left hand side 	
		 shows the exponential convergence of $\Obj(X_t)$ to $\optObj$ following the switching time (denoted
		$\tau_\eta$ in the figure).    The trajectory on the right hand side reaches the set $\RoA$ at approximate $t=0.7$.		 
The simulation was stopped at  the time $T$ when the performance constraint  $\Obj(X_T)-\optObj = \epsy $ was reached, which was approximately $T=1.4$.

It was rare to obtain a trajectory converging to the second optimizer.
The spectral decomposition provides an explanation.  Choose the sign
of $h_2$ so that $X_0\in H_+$.  The upper optimizer also lies in
$H_+$, while the other optimizer lies in $H_-$.  Correspondingly,
write
\[
    \RoA=\RoA_+\cup\RoA_-,
\]
with $\overline{\RoA}_+\subset H_+$ and
$\overline{\RoA}_-\subset H_-$.

This asymmetry has a precise interpretation in terms of hitting times.
Let $A_+\subset\RoA_+$ and $A_-\subset\RoA_-$ be closed neighborhoods
of the respective optimizers, each with nonempty interior.  For any
$\gap_{\!\temp}<\bar r<|\lambda_3|$, (16) gives
\[
    \Expect_{X_0}\bigl[\exp(\bar r\tau_{A_+})\bigr]<\infty .
\]
On applying the same result to the eigenfunction $-h_2$, whose positive
and negative domains are interchanged, we obtain
\[
    \Expect_{X_0}\bigl[\exp(\gap_{\!\temp}\tau_{A_-})\bigr]=\infty .
\]
Thus the zero set of $h_2$ represents a genuine stochastic barrier:
from the specified initial condition, reaching the optimizer in the
opposite nodal domain occurs on the slow time scale associated with
the spectral gap.

The   behavior observed in these numerical examples is also
closely related to conclusions obtained for discrete Langevin algorithms \cite{tzeliarag18}.

\begin{figure}[h]
	\centering
	\includegraphics[width=\hsize]{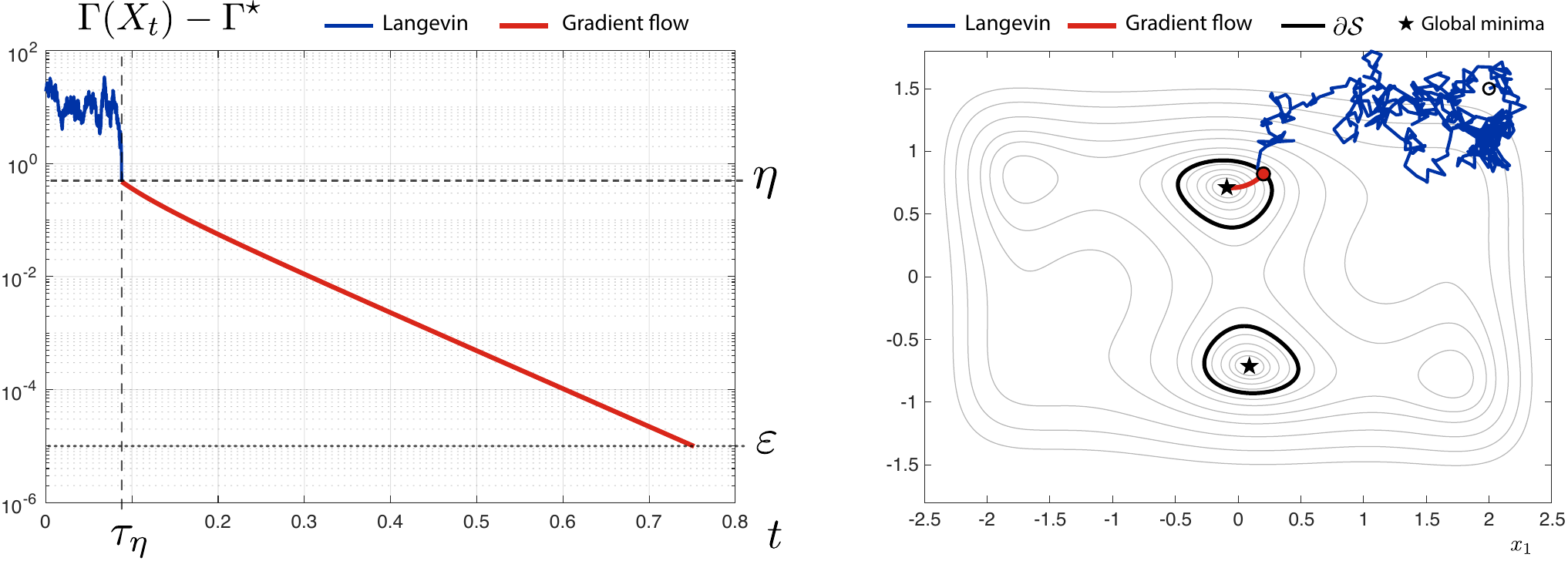}  
	\caption{A repetition of the experiment illustrated in
\Cref{f:Cost+PhaseCamelLG} using $\temp=10$, resulting in entry to $\RoA$ at approximately $\tau_\eta = 0.1$.} 
		\label{f:Cost+PhaseCamelLG10}
\end{figure} 

\wham{High temperature:}  With $\temp=10$  the spectral gap is $\gap_{\!\temp} =|\lambda_2| \approx 5.2908$ and 
first four eigenvalues are $\{\lambda_i \} = \{ 0, -5.2908, -12.6249, -  17.1669 \}$. 
We obtain a similar eigenvalue ratio as in the low temperature example:   
$\kappa = \lambda_3/\lambda_2  
\approx  2.39$.  

\begin{figure}[h!]
	\centering
	\includegraphics[width=\hsize]{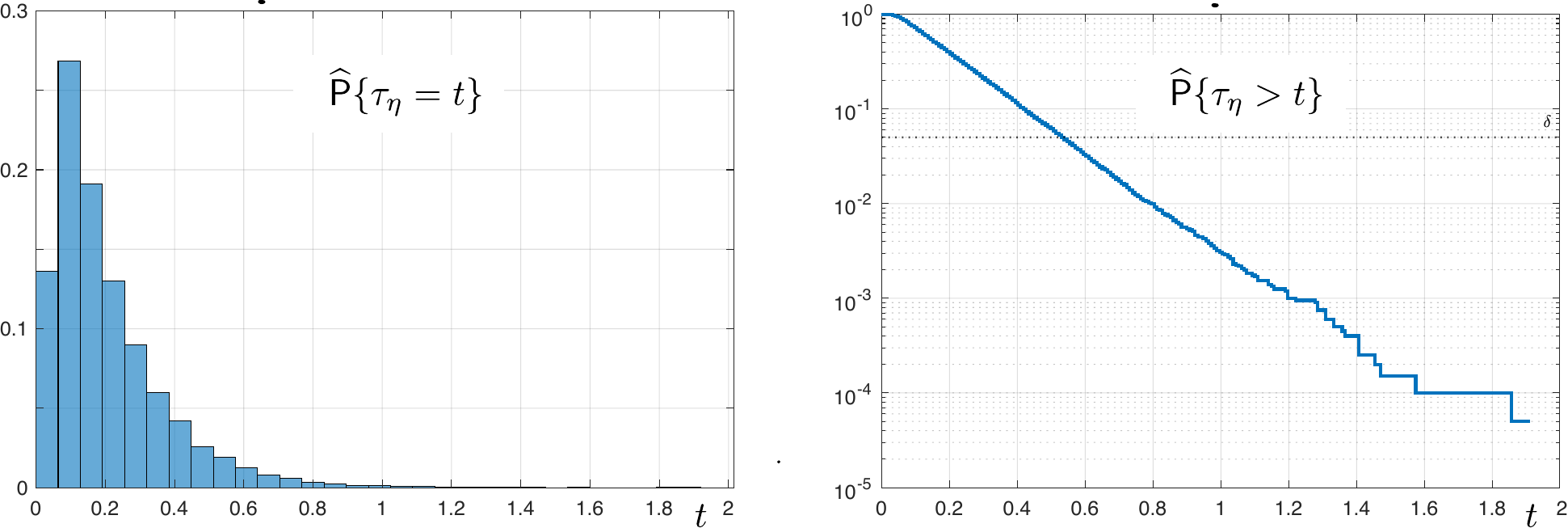}  
	\caption{Empirical distribution of $\tau_\eta$  with $\temp=10$. }
		\label{f:success_temp10}
\end{figure}

 On increasing $\temp=0.2$ to $\temp=10$, the switching time $\tau_\eta$  is  significantly reduced on average. 
 \Cref{f:Cost+PhaseCamelLG10} illustrates a successful run resulting in  $\tau_\eta < 0.1$.    	\Cref{f:success_temp10} 
shows an estimate of the distribution of  $\tau_\eta$ for $\temp=10$ based on 20,000 independent trials from the   single initial condition  $X_0 = (2.0; 1.5 )$.

%\begin{figure}[h]
%	\centering
%	\includegraphics[width=0.5\hsize]{CostCamelLG10.pdf}  
%	\caption{A repetition of the experiment illustrated in
%			\Cref{f:CostCamelLG} with temperature increased to $\temp=10$, resulting in entry to $\RoA$ at approximately $\tau_\eta = 0.1$.}
%		\label{f:CostCamelLG10}
%\end{figure} 
%
%
%
%
%\begin{figure}[h]
%	\centering
%	\includegraphics[width=0.5\hsize]{PhaseCamelLGtemp10.pdf}  
%	\caption{A repetition of the experiment illustrated in
%			\Cref{f:PhaseCamelLG} with temperature increased to $\temp=10$, resulting in entry to $\RoA$ at approximately $\tau_\eta = 0.1$. 
%	}
%		\label{f:PhaseCamelLG10}
%\end{figure}   

\smallskip 

It seems clear that in    this example, for this large set $\RoA$, \textit{there is no reason to use a small temperature}.

\begin{figure}[h!]
	\centering
	\includegraphics[width=0.6\hsize]{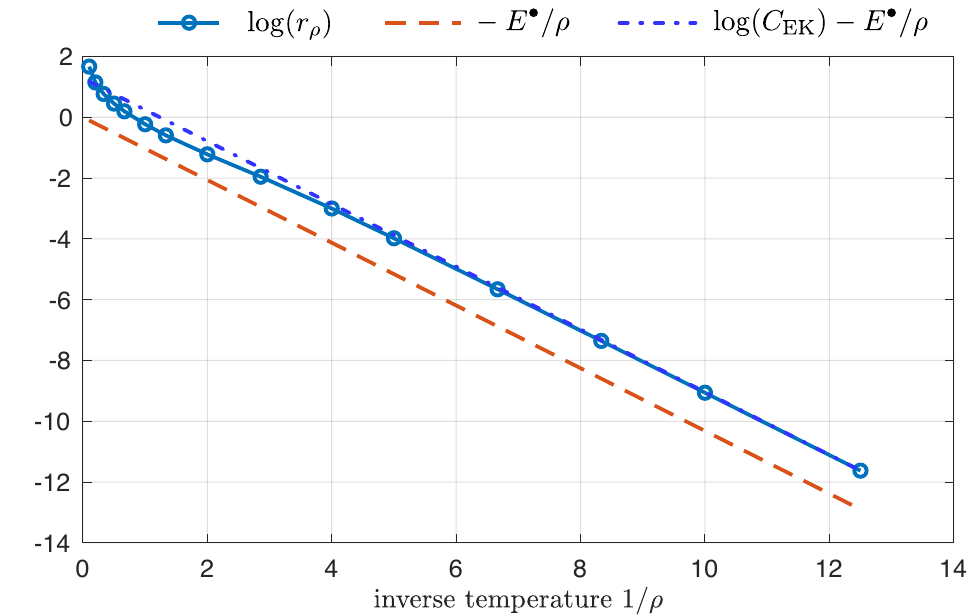}  
	\caption{The logarithm of the spectral gap for the six hump camel along with two approximations.}
		\label{f:SpectralGapComparisonsCamel}
\end{figure}

Theory in this paper focuses on small temperature only because useful
approximations are available.  In this example, these
approximations are informative over a surprisingly wide range of
temperatures.  In particular, \Cref{t:low-temp-basic}  below gives, 
\begin{equation}
    -\temp\log(\gap_{\!\temp})
    \longrightarrow \Dbar \,, 
    \qquad \temp\downarrow0  \,,
\label{e:low-temp-master_lambdaA}
\end{equation} 
where the dominant energy barrier $\Dbar$ is defined in \eqref{e:EnergyBarrier}.
We obtain in this example
$\Dbar\approx1.03162845$ and observe in
\Cref{f:SpectralGapComparisonsCamel} that
$\log(\gap_{\!\temp})$ is approximately affine in $1/\temp$, with slope close
to $-\Dbar$, over a substantial range of temperatures.
This is consistent with \eqref{e:low-temp-master_lambdaA}.

A sharper approximation may be obtained under additional assumptions
\[
    \log(\gap_{\!\temp})
    =
    \log(C_{\text{EK}})
    - \frac{\Dbar}{\temp}
    +o(1)\,, 
    \qquad \temp\downarrow0\,, 
  \]
which explains the second approximation appearing in \Cref{f:SpectralGapComparisonsCamel}.
The constant $C_{\text{EK}}$ is the Eyring--Kramers prefactor:  Since the six-hump camel objective
has two symmetry-related global minima, separated by a unique saddle $s\in\Re^2$,  we obtain from
\cite[Cor.~2.15]{MenzSchlichting2014} 
\begin{equation}
C_{\rm EK}
=
\frac{|\lambda_-(s)|}{\pi}
\sqrt{\frac{\det H_m}{|\det H_s|}} \,, 
\label{e:numCEKcamel}
\end{equation}
where $m$ is either global minimizer,  $H_m$ is the common Hessian of $\Obj$ at these locations,  
 and $H_s$ denotes the Hessians at the saddle.    
  
 \notes{
For the normalization used here, the optimal Poincar\'e constant
$\varrho_\temp$ and the spectral gap satisfy
$\gap_{\!\temp}=\temp\varrho_\temp$.  
}

\subsection{Objective with many local minima}  
\label{s:rast}

The Rastrigin objective $\Obj \colon\Re^d \to\Re_+$ is defined for any dimension $d$ by
\begin{equation}
\Obj(x) =  \sum _{i=1}^{d} g(x_i)\,, \quad \textit{with} \quad     g(z)=10+z^2-10\cos(2\pi z)\,, 
\label{e:Rastrigin}
\end{equation}
with unique minimizer $x^\star =0$,  giving   $\optObj=0$.  

In all numerical experiments in this subsection we maintained  
$\epsy=10^{-5}$,  $\delta=0.05$, and $\uptauGD=10^{-3}$, but increased the exploration stepsize by ten-fold, 
 $\uptau=10^{-3}$.

 \subsubsection{Spectral analysis}

The dominant energy barrier is independent of dimension.   The proof of the following is postponed to the Appendix.

\begin{proposition}
\label[proposition]{t:Rastrigin}  

The following hold for the Rastrigin objective, for any  $d\ge1$: 

\wham{(i)} 
The dominant energy barrier is
$
    \Dbar=19.2563\ldots$, 
independent of $d$.   Hence the approximation of the spectral gap \eqref{e:low-temp-master_lambdaA} is independent of dimension.  

%    \Dbar=19.2563139339 ...

\wham{(ii)}  $\RoA= \{x : \Obj(x ) \le  \eta\}$ satisfies the requirements of \Cref{a:jacquot} for any fixed 
$    0<\eta\le 0.99$.  

\wham{(iii)} 
For every $\temp>0$ the spectral
gap $\gap_{\!\temp}$ of the $d$-dimensional Langevin generator is exactly equal to the spectral gap with $d=1$.    
Moreover, the eigenvalue $\lambda_2=-\gap_{\!\temp}$ has multiplicity   $d$.
\clThm
\end{proposition}

\begin{figure}[h]
	\centering
 	\includegraphics[width=\hsize]{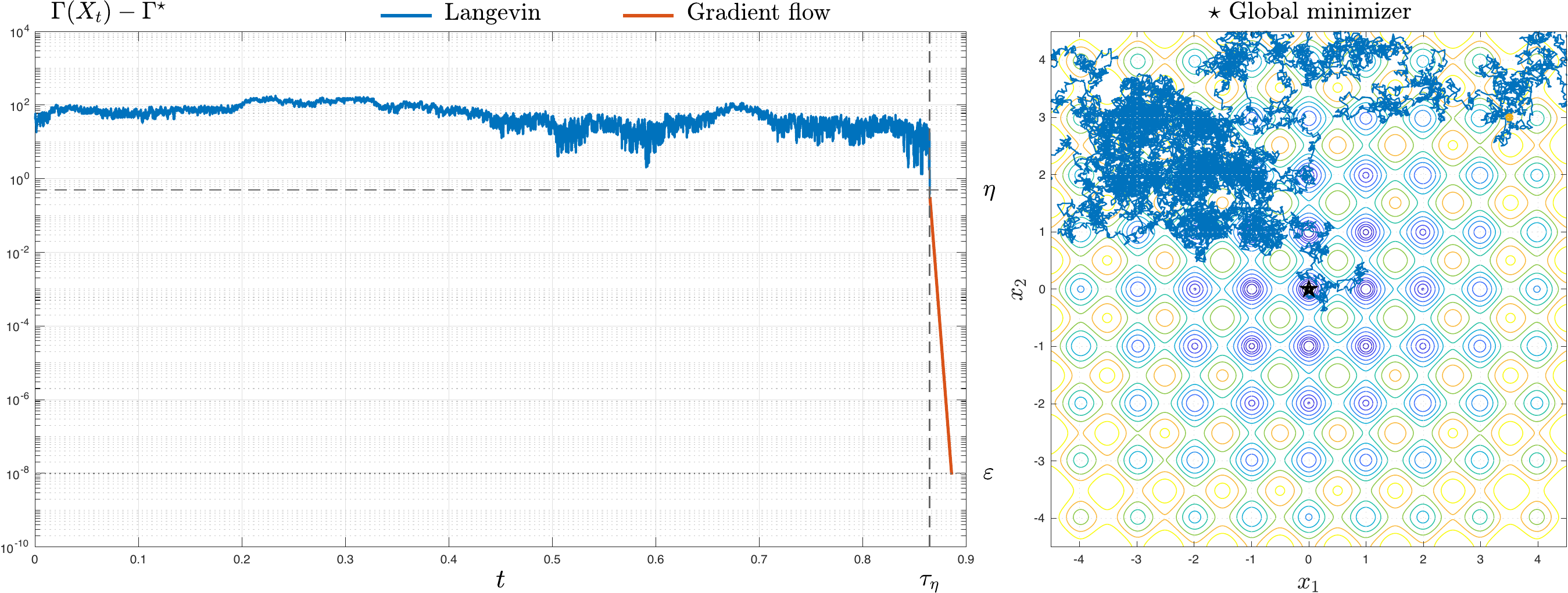}%rastaLG_temp25_eta0.5.pdf}  
\caption{A single run of $\LGFE{\temp,1}$  for the
Rastrigin objective in dimension 2. 
}
		\label{f:Rastrigin}
\end{figure}

\Cref{f:Rastrigin} shows results from a typical single  run of the Langevin--gradient method $\LGFE{\temp,1}$ in dimension $d=2$,
using  $\temp=50$ and $\eta=0.5$.
The plot on the left hand side shows the evolution of the excess cost:
At the entrance time $\tau_\eta$ the noise is removed and gradient
flow approaches $x^\star=0$ at an exponential rate.    
Shown on the right is the trajectory of $X_t$ in $\Re^2$.   
We observe aggressive global exploration during the  Langevin phase,
 traversing many local wells before
first entering the target set
$\RoA=\{x:\Gamma(x)-\Gamma^\star <\eta\}$.

\begin{figure}[h!]
	\centering
	\includegraphics[width=0.8\hsize]{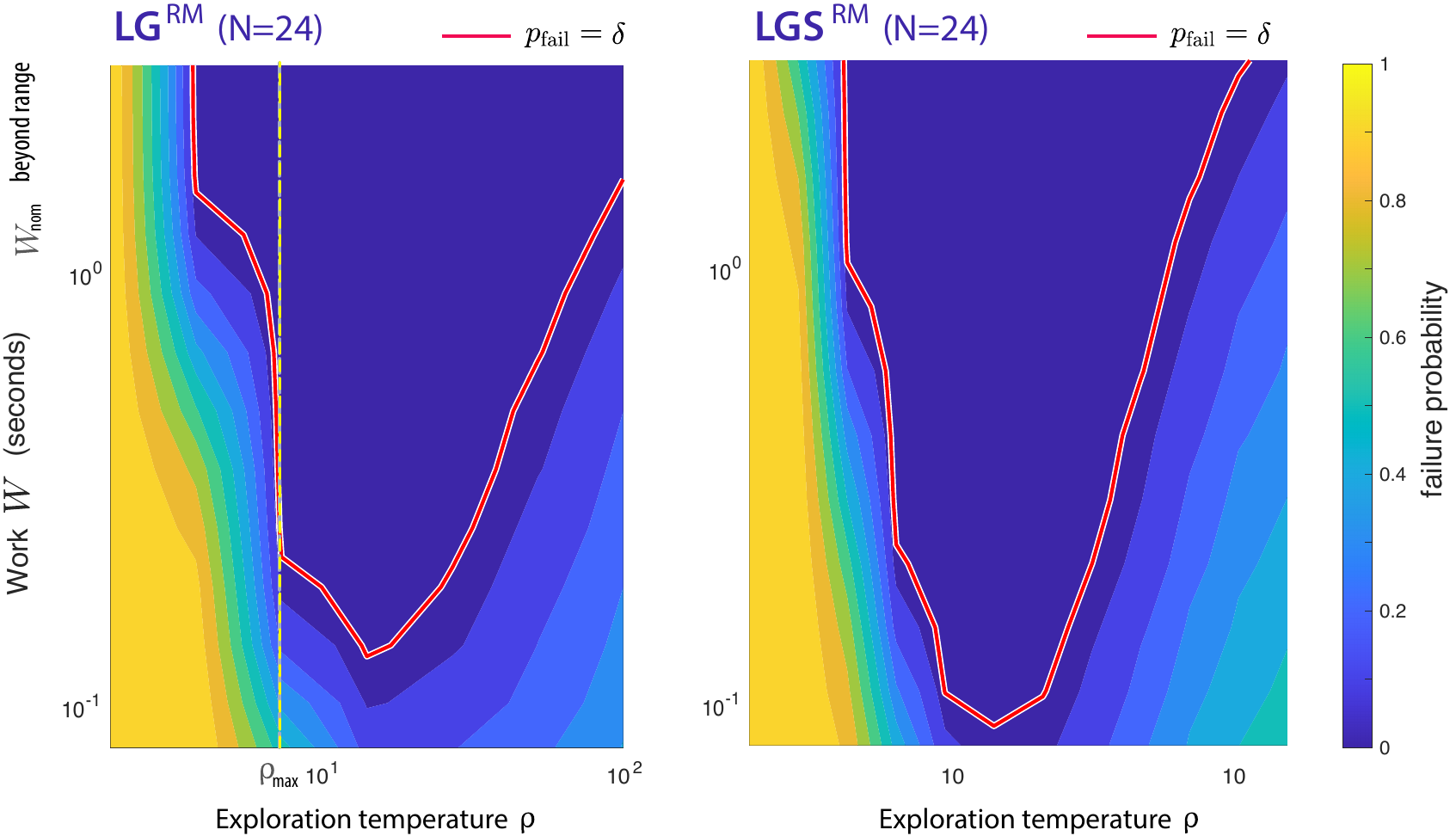}  
	\caption{Empirical performance of two best-state Langevin--gradient methods
	for the Rastrigin objective in dimension~2.  The left hand side shows results for  $\LGRM{\temp,N}$ and the right hand side  $\LGSRM{\temp,N}$, each using $N=24$.}
		\label{f:WorkVsProbFailureRastrigin_d2_N24}
\end{figure}

The performance of $\LGRM{\temp,N}$
and the state-dependent temperature vsion
 $\LGSRM{\temp,N}$
  are similar.
  Shown in  \Cref{f:WorkVsProbFailureRastrigin_d2_N24} are contour plots representing   the probability of failure for
 each of the two methods, with   $N=24$.   We see that the maximal temperature is a significant under-estimate of the best value of $\temp$ in this example, which explains the large choice of $\temp=50$ used in \Cref{f:Rastrigin}.    The contour plot for  $\LGSRM{\temp,N}$ shows a modest improvement in the best work observed.

\begin{figure}[h!]
	\centering
	\includegraphics[width=0.5\hsize]{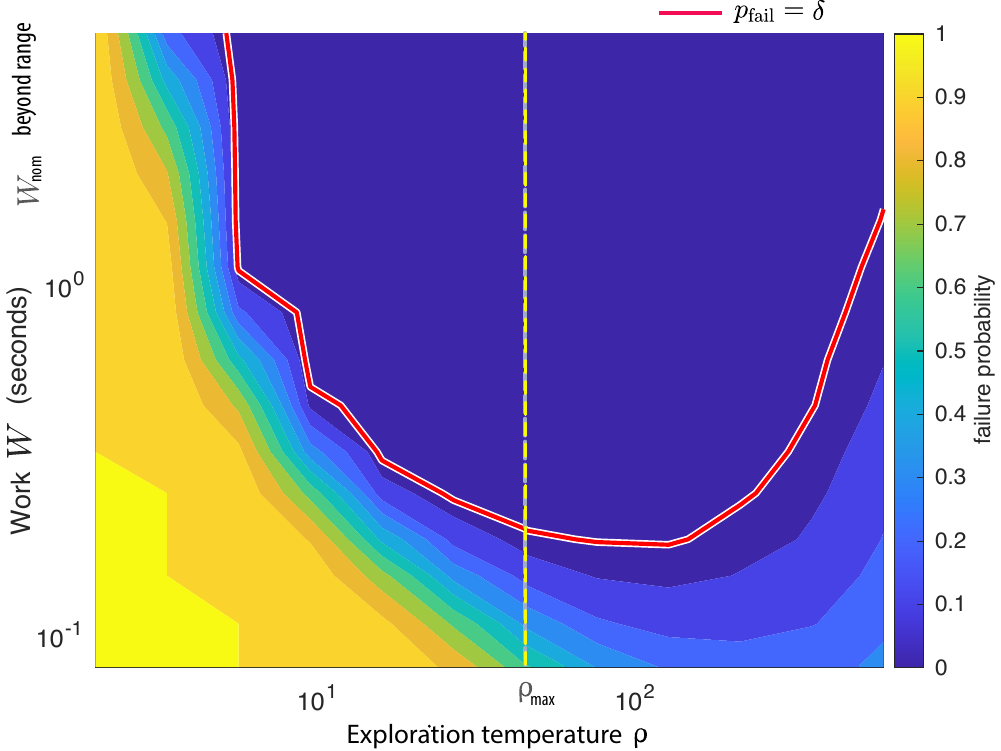}  
	\caption{Empirical performance of the best-state Langevin--gradient
scheme $\LGRM{\temp,N}$ for the Rastrigin objective in dimension 2, with  
$N=300$.}
		\label{f:WorkVsProbFailureRastrigin_d2_N300}
\end{figure}

Simulations were repeated for $N=300$ in $\LGRM{\temp,N}$  leading to the contour plot shown in 		\Cref{f:WorkVsProbFailureRastrigin_d2_N300}.
We find that the best performance observed is similar to what is obtained with $N=24$, but the range of acceptable temperature is improved.  

%\clearpage

\begin{figure}[h]
	\centering
 	\includegraphics[width=0.6\hsize]{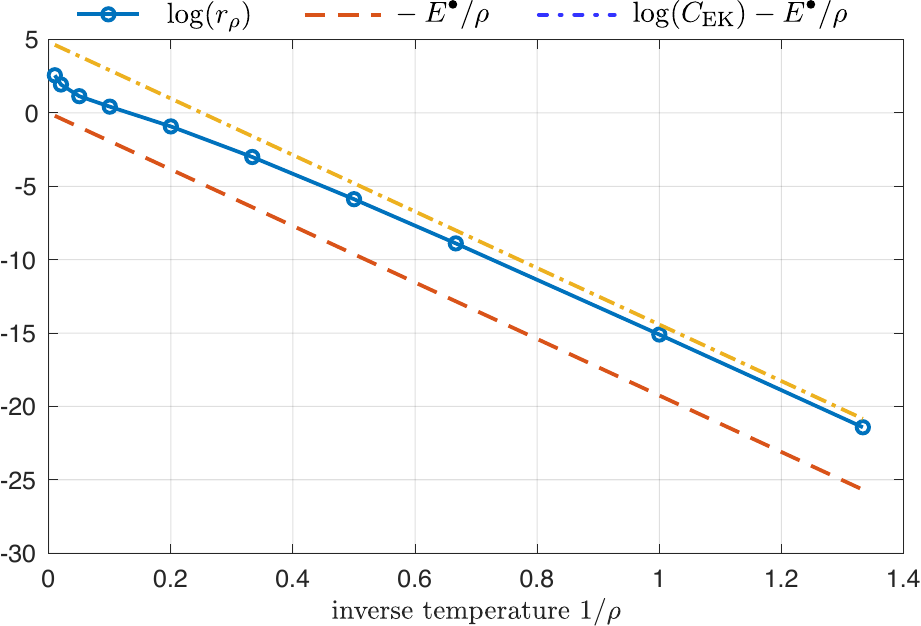}  
\caption{log spectral gap and approximations.}
		\label{f:LogGapVsInverseTemperature_rastrigin_v3}
\end{figure}

     \Cref{f:LogGapVsInverseTemperature_rastrigin_v3} provides plots of the spectral gap as a function of 
inverse temperature, showing similar behavior as in \Cref{f:SpectralGapComparisonsCamel}.   
 As anticipated by \Cref{t:Rastrigin},
in numerical approximation of the spectral gap, we observe 
 repeated eigenvalues at $\lambda_2=\lambda_3 = -\gap_{\!\temp}$ for any $\temp>0$.  
 There is an associated two dimensional eigenspace spanned by a pair of eigenfunctions
 $h_2, h_3$. 
  
Computation and approximation of the eigenvalues is simplified because  the spectral gap is independent of $d$ because
(see  \Cref{s:RastriginAppendix}.
for details) from which we obtain for any dimension, 
\[
    \gap_{\!\temp,d}
    \sim
    C_{\rm EK}
    \exp\left\{
       -\frac{g(s)-g(m)}{\temp}
    \right\},
    \qquad \temp\downarrow0,
\]
where the prefactor   $   C_{\rm EK}$  is derived from   the one-dimensional
Eyring--Kramers approximation.  Let $m>0$ denote the first
nonzero local minimum of $g$, and let $s\in(0,m)$ denote the
intervening local maximum.  Then,
\[
    C_{\rm EK}
    =
    \frac{1}{2}
    \frac{|g''(s)|}{\pi}
    \sqrt{\frac{g''(m)}{|g''(s)|}}
    =
    \frac{1}{2\pi}
    \sqrt{g''(m)|g''(s)|}.
\]
The division by 2 when comparing to \eqref{e:numCEKcamel}
is a consequence that here the optimizer is unique.

\begin{figure}[h!]
	\centering
 	\includegraphics[width=0.45\hsize]{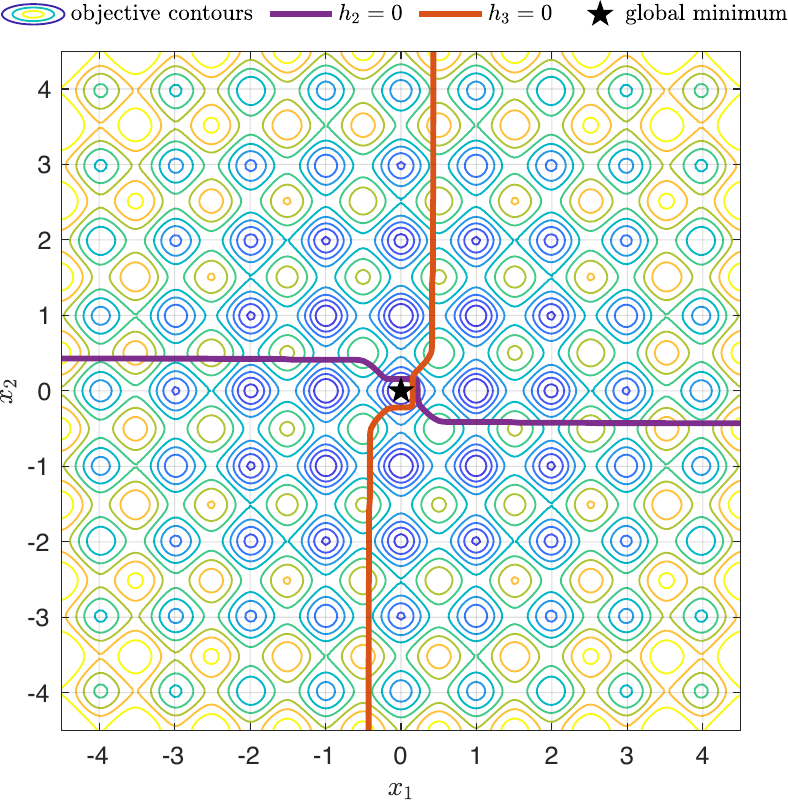}%rastaLG_temp25_eta0.5.pdf}  
	
%See the remarkable data in \verb+SpectralGapData_rastrigin_v3.csv+.    
% 
%
%{\tiny
% \[
%\begin{array}{llllll}
%\temp & -\lambda_2 & -\lambda_3 & -\lambda_4 & \text{crude gap Approx} &  \text{better gap Approx}\\ 
%0.75 & 4.96744515037313E-10 & 4.96774635201682E-10 & 7.64247214642751E-10 & 7.07053303322051E-12 & 8.88219107646807E-10\\ 
%1 & 2.76118811961289E-07 & 2.76118826968071E-07 & 4.83077739989776E-07 & 4.33599843432378E-09 & 5.4469962052335E-07\\ 
%1.5 & 0.000137771435447336 & 0.000137771435463143 & 0.000264829640358661 & 2.65904739205989E-06 & 0.000334036584040997\\ 
%2 & 0.00280286002946458 & 0.00280286002948095 & 0.00547201816088641 & 6.58482986441091E-05 & 0.00827203787704969\\ 
%3 & 0.0496450793426911 & 0.0496450793427681 & 0.0983096180627782 & 0.00163065857617709 & 0.204847654144813\\ 
%5 & 0.399134213402885 & 0.399134213402919 & 0.798268426805927 & 0.0212528816280528 & 2.66984334362045\\ 
%10 & 1.53514732231274 & 1.53514732231278 & 3.07029464462566 & 0.145783680938755 & 18.313732555169\\ 
%20 & 3.13630919813925 & 3.13630919813935 & 6.27261839627713 & 0.381816292133737 & 47.9647750304859\\ 
%50 & 6.94848533658342 & 6.94848533658349 & 13.8969706731692 & 0.680364715222955 & 85.4692195610113\\ 
%100 & 13.0233604527759 & 13.0233604527764 & 26.0467209055392 & 0.824842236565851 & 103.618868884375
%\end{array}
%\]
%}

\caption{Zero sets of two linearly independent eigenfunctions associated with  $\lambda_2=\lambda_3 = -\gap_{\!\temp}$. 
}
		\label{f:SecondEigenfunctions_rastrigin_rho_1_v5}
\end{figure}

\Cref{f:SecondEigenfunctions_rastrigin_rho_1_v5} shows the zero sets of two eigenfunctions chosen to be orthogonal in $L_2(\pi_\temp)$, for $\temp=1$.   They are related by
\[
h_3(x) = h_2(T x)
\]
where $T$ is a 90 degree rotation matrix.   

 We are not aware of an extension of \eqref{e:MGFhuimeysch04a}  when $\lambda_2$ is repeated.  

% \begin{figure}[h!]
%	\centering
%	\includegraphics[width=0.55\hsize]{WorkVsProbFailureRastrigin_d5_N24.pdf}  
%	\caption{Empirical performance of the best-state Langevin--gradient
%scheme $\LGRM{\temp,N}$ for the Rastrigin objective in dimension 5, with  
%$N=24$.}
%		\label{f:WorkVsProbFailureRastrigin_d5_N24}
%\end{figure} 

\begin{figure}[h!]
	\centering
	\includegraphics[width=0.7\hsize]{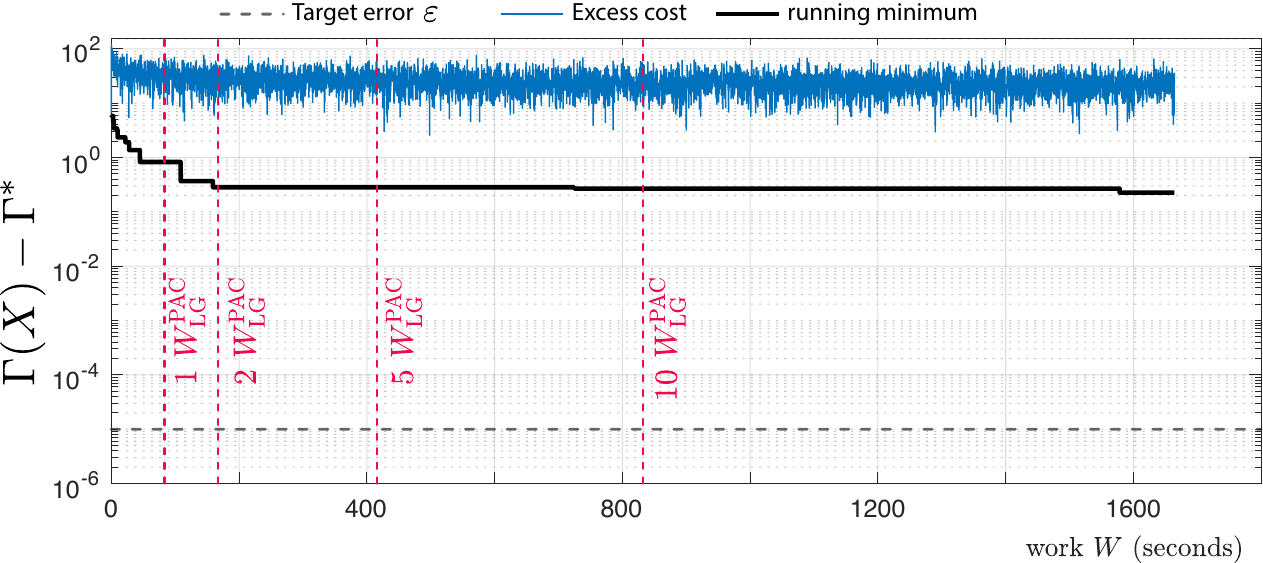}  
\caption{Representative sample path from simulated annealing for the Rastrigin objection with dimension $d=2$.}
		\label{f:Rastrigin_SimA_LongTrajectory}
\end{figure}

Simulated annealing experiments used 
the logarithmic cooling schedule \eqref{e:cooling} with   $E = 2\Dbar$ and   
and the same Euler step $\uptau=10^{-3}$.  
\Cref{f:Rastrigin_SimA_LongTrajectory}
shows results from a typical experiment.

It would take far too long to  estimate the substantially larger PAC work
required by classical simulated annealing.  	Instead, we compare its performance
at workload values
$K\times W_{\rm LG}^{\rm PAC}$ for a range of $K\ge 1$,  where
$W_{\rm LG}^{\rm PAC}=83.20$ seconds is the calibrated work at which
$\LGRM{\temp,N}$ attains the PAC target.

With 200 independent trials there were no   successes for the range of $K$ considered, which included $K=5$ corresponding to 
\[
    W=5W_{\rm LG}^{\rm PAC}=416.01\ {\rm sec},
\]
The observed success probability was zero.

The disparity is consistent with the different roles of exploration
and exploitation in the two schemes.  The logarithmic simulated-
annealing schedule begins at temperature
$\tempSimA(0)=E=38.51$ and cools only gradually.  At the work
budgets $W_{\rm LG}^{\rm PAC}$, $2W_{\rm LG}^{\rm PAC}$, and
$5W_{\rm LG}^{\rm PAC}$, the corresponding annealing temperatures
remain approximately $5.4$, $4.9$, and $4.4$, respectively.  Hence a
large fraction of the simulated-annealing computation is spent cooling
toward the temperature range that is already effective for
$\LGRM{\temp,N}$.  The latter instead performs exploration
directly at $\temp=4$, retains the best state encountered among the
restart trajectories, and then switches to deterministic gradient
descent exploitation phase.

\subsubsection{Dimension dependency and finite-temperature design}     

The challenge presented by dimension is clear in this example on noting that for $d=1$ the number of local minima is $63$,
so that for dimension $d=2$ considered in the previous discussion, this value increases to  $63^2 =3969 $.   
The next examples consider higher dimensions, resulting in
\[
\textit{Number of local minima:} \qquad
63^{5}\approx 10^9  \ \ \text{($d=5$),} 
\qquad
63^{10}\approx 10^{18}  \ \ \text{($d=10$)}
\]

\begin{figure}[h!]
	\centering
	\includegraphics[width=0.8\hsize]{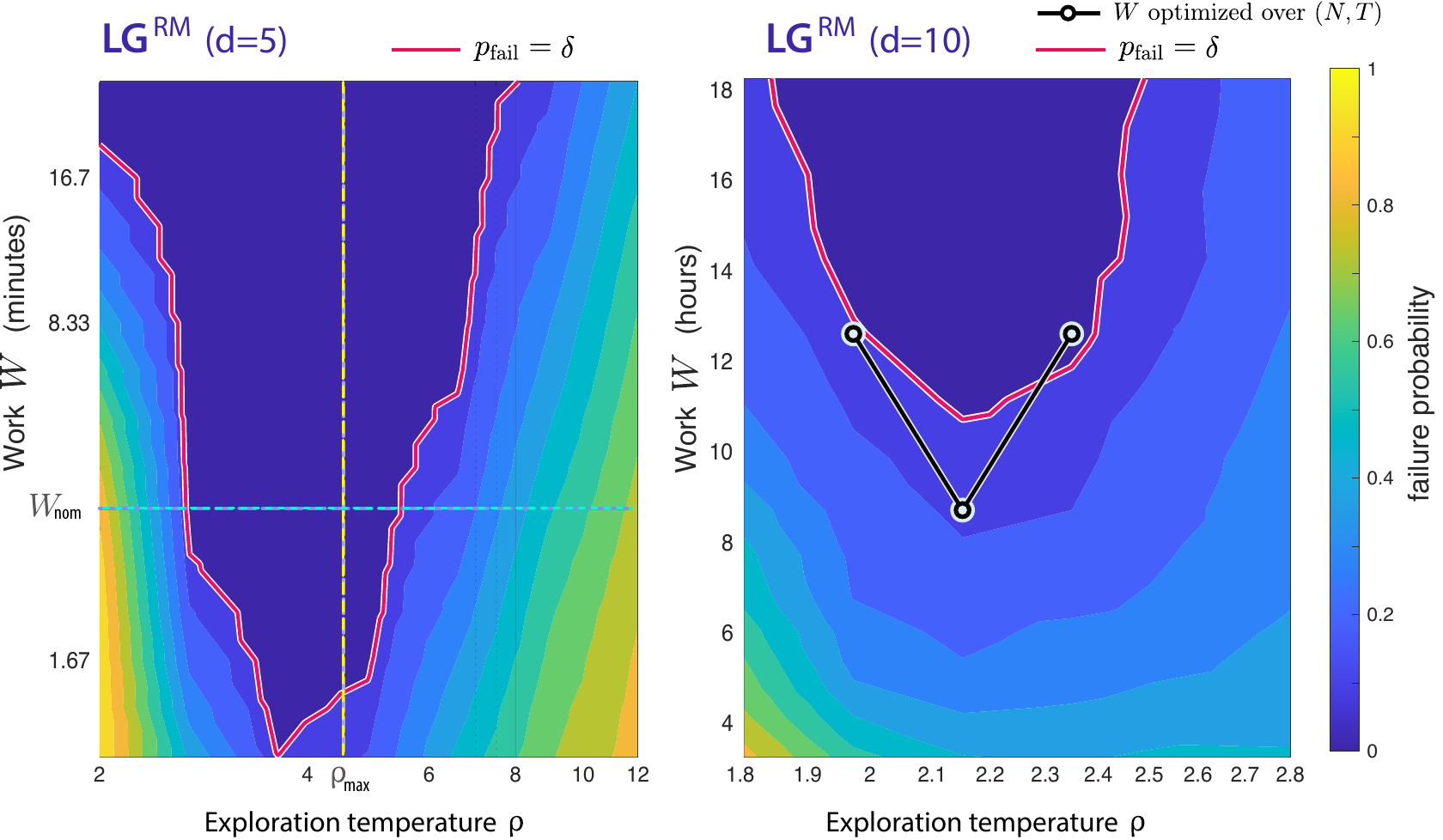}  
	\caption{Empirical performance of $\LGRM{\temp,N}$ 
	for the Rastrigin objective.  The left hand side shows results for $d=5$ and $N=24$.  The  
	   right hand side displays results for $d=10$ with  $N=300$, and in addition a plot of the realized work for a smaller range of temperature,
	   minimized  over a range of $(N,T)$ for each $\temp$ in the selected range \eqref{e:NTrangeForDim10}.}
		\label{f:WorkVsProbFailureRastrigin_d5_N24_d10_N300}
\end{figure}

With  dimension   $d=5$ we obtained results qualitatively consistent with $d=2$.  
A contour plot   showing failure probability as a function of work and temperature is shown on the left hand side of
\Cref{f:WorkVsProbFailureRastrigin_d5_N24_d10_N300} using $N=24$ as in previous experiments.  
Observe that work is measured in minutes (rather than the units of seconds used in \Cref{f:WorkVsProbFailureRastrigin_d2_N24}),
and the minimal work required to meet the PAC criterion rose from a fraction of one second to about one minute.

With dimension $d=10$, the low-temperature design rule
\eqref{e:nomTempAndTim} suggests increasing $N$ in order to reduce the
nominal switching time $\Tf^0$.  For example, taking $N=300$ gives
$\Tf^0\approx 2$ and increases the nominal temperature by a factor
$300/24$, giving $\tempLG\approx 23$.  This value is far outside the
low-temperature regime, and consequently the approximation underlying
\eqref{e:nomTempAndTim} should not be expected to provide useful
finite-temperature tuning in this example.

The contour plot shown on the right hand side of
\Cref{f:WorkVsProbFailureRastrigin_d5_N24_d10_N300}
illustrates the sensitivity of the probability of satisfying the PAC
constraint to temperature and total work.  Note that the units are now in hours.

The empirical PAC boundary
is strongly U-shaped, with its minimum near
$\temp\approx2.2$.  Hence, although increasing $N$ to $300$ gives the
nominal value $\tempLG\approx23$, direct finite-temperature evaluation
of $\LGRM{\temp,N}$ favors a temperature approximately one order of magnitude
smaller.  This discrepancy is not unexpected, since the nominal
parameter values are motivated by low-temperature asymptotics, whereas
the numerical experiment optimizes performance at finite temperature.

Note that all previous contour plots were based on 2000 independent trials for each pair $(\temp,\Tf^0)$.   
With dimension $d=10$   this value was reduced to 200, but the run time required to obtain the contour plot   was about 30 hours.  

While the nominal temperature prediction \eqref{e:temp_nom} is not predictive in this example, 
the theory leading to large $N$ was justified through a second set of experiments
 to determine whether the work could
be reduced by optimizing the number of Langevin replicas as well as
temperature and exploration time.  We considered
\begin{equation}
N\in\{100,150,200,250,300\},\qquad
    \temp\in\{1.8,2.0,2.2,2.4,2.6\},
\label{e:NTrangeForDim10}
\end{equation}    
using $100$ independent trials and several exploration horizons.
The results  also shown on the right hand side of
\Cref{f:WorkVsProbFailureRastrigin_d5_N24_d10_N300}
are consistent with the contour plot:  The smallest empirical work was again obtained at
$\temp=2.2$, now with
\[
    N=300,\qquad \Tf^0=2146.5 \  \text{secs.,}
\]
for which the observed failure probability was $0.04$ and the  calibrated work was $8.69$ hours.  

%Recall from  \Cref{f:Rastrigin_Transition_CRN_Theory} that a temperatures near  $\temp=2.2$ also maximized the single hop probability 
%to $x^\star$ from the nearby local minimizer $m_2$.

Reducing the number of replicas
provided little improvement: the best empirically feasible design with
$N=250$ required $9.26$ hours, while no tested design with $N\leq200$
attained an empirical failure probability below $\delta=0.05$.
The   $200$-trial experiment with $N=300$ gave a minimum of  
$11.12$ hours, also at $\temp=2.2$.    

These long run estimates are not surprising, given that 
we are faced with approximately $10^{18}$ local minima, with only one  global minimizer.

The main takeaways are that computational complexity is very sensitive to temperature, and increases dramatically with dimension
as predicted by discussion in \Cref{s:highd}.

 \begin{figure}[h!]
	\centering
	\includegraphics[width=0.85\hsize]{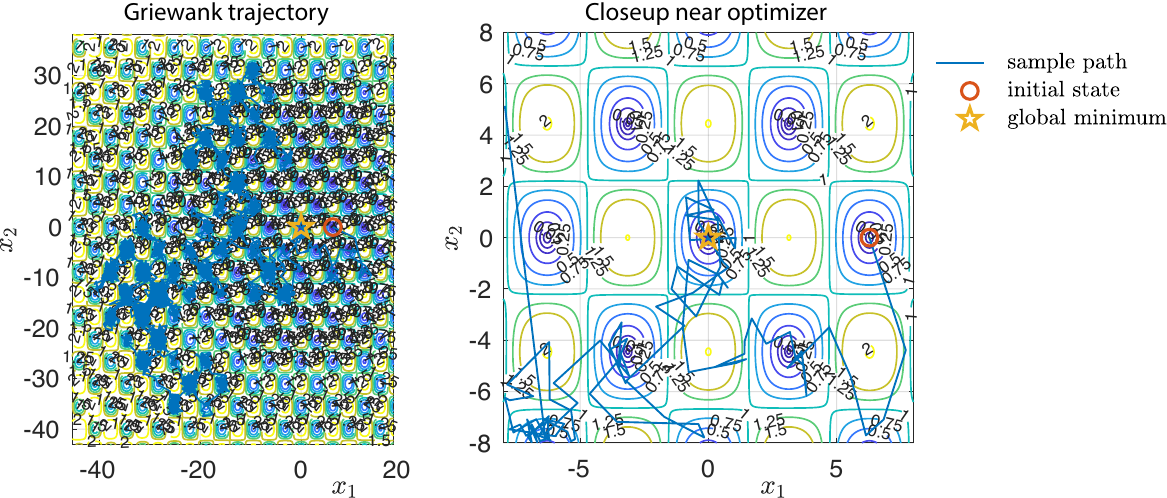}  
\caption{
Representative sample-path  of simulated annealing for the Griewank objective function in dimension $d= 2$.   }
		\label{f:Griewank_d2_SimA_TrajCloseup}
\end{figure}

\subsection{Griewank}
\label{s:SimA+LG_Griewank}

The Griewank function in two dimensions is
$$
 \Obj(x) = \frac{1}{4000} \|x\|^2 + 1-\cos(x_1)\cos ( {x_2}/{\sqrt2} )
 $$
with unique minimizer $x^\star=0$ and  $\Dbar\approx0.997$.  This objective is more challenging than the others considered because a suitable value of $\eta$ is very small.

The stepsize \(\uptau=10^{-2}\) was used in the following experiments.

In the following we compare $\LGRM{\temp,N}$ with simulated annealing.

\subsubsection{Comparisons}

 \begin{figure}[h!]
	\centering
	\includegraphics[width=0.55\hsize]{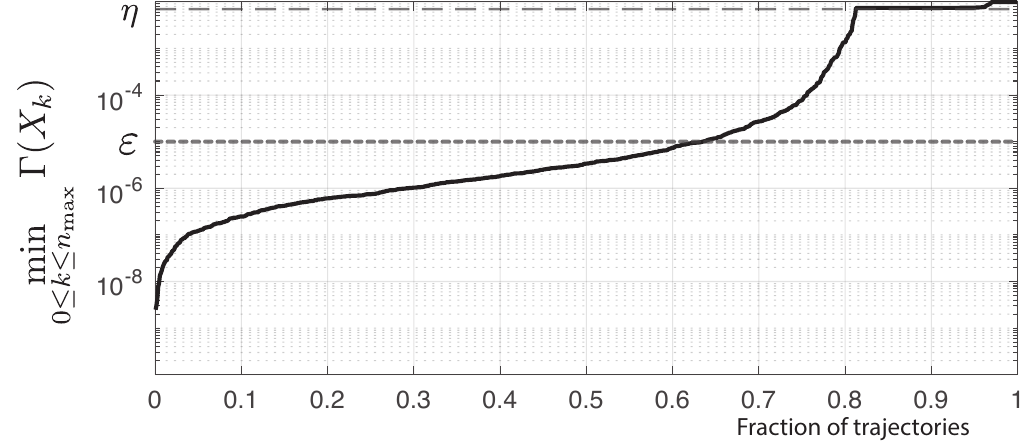}  
\caption{
Representative sample-path  of simulated annealing for the Griewank objective function in dimension $d= 2$.   }
		\label{f:Griewank_d2_SimA_EmpDist}
\end{figure}

\wham{Simulated Annealing} 

The cooling schedule \eqref{e:cooling} was used with 
$
E=2 \approx 2\Dbar  >\Dbar
$.
We fixed the time horizon of  $T=10^5$, so that with   \(\uptau=10^{-2}\) this corresponds to a computational budget of 
approximately $446$ seconds.  
\Cref{f:Griewank_d2_SimA_TrajCloseup} shows a sample trajectory.

 \Cref{f:Griewank_d2_SimA_EmpDist} provides a representation of the empirical distribution of   excess cost, 
$\displaystyle
\min_{0\leq k\leq n_{\max}}\Obj(X_k)
$, based on    $1000$ independent trials.
In particular,
\wham{$\circ$}
$812$ ($81.2\%$) reached the sublevel set $\{x: \Obj(x)\le \eta\}$
\wham{$\circ$}
 $635$ ($63.5\%$) reached the target set to satisfy the PAC requirement, $\{x: \Obj(x)\le \epsy\}$, with   $\epsy=10^{-5}$.  
 
The probability of missing the PAC   probability is thus $1-0.635=0.365$, which is well above the target $\delta = 0.05$.    Hence   the computational budget  required for simulated annealing to satisfy the PAC criterion is expected to be substantially more than  $446$ seconds.

We see next that the required computational work is reduced to less than 10 seconds when using \(\LGRM{\temp,N}\), for a wide range of $\temp $ and $N$.

 \wham{Langevin--gradient}

Recall that the energy barrier is moderate, $\Dbar\approx0.997$, 
    but to ensure $\RoA_\eta$  satisfies the required assumptions imposes    
$ \eta<0.0074$.  
In view of the guideline    \eqref{e:nomTempAndTim} for choice of $    \Tf^0  $, 
the small value of $\eta$ will motivate large $N$ in  $\LGRM{\temp,N}$.    
For example, 
with $N=500$  the expression \eqref{e:nomTempAndTim} gives
    $    \Tf^0 \approx 5$ and  $\tempLG \approx   
 N \times 0.0074 /L \approx 2.8 $ from \eqref{e:temp_nom}.

%exp(2*1*3/(500*0.0074))\approx 5
%\exp\left\{  2\frac{\Dbar L }{N\eta}  \right\} 
  
\begin{figure}[h!]
	\centering
	\includegraphics[width=0.8\hsize]{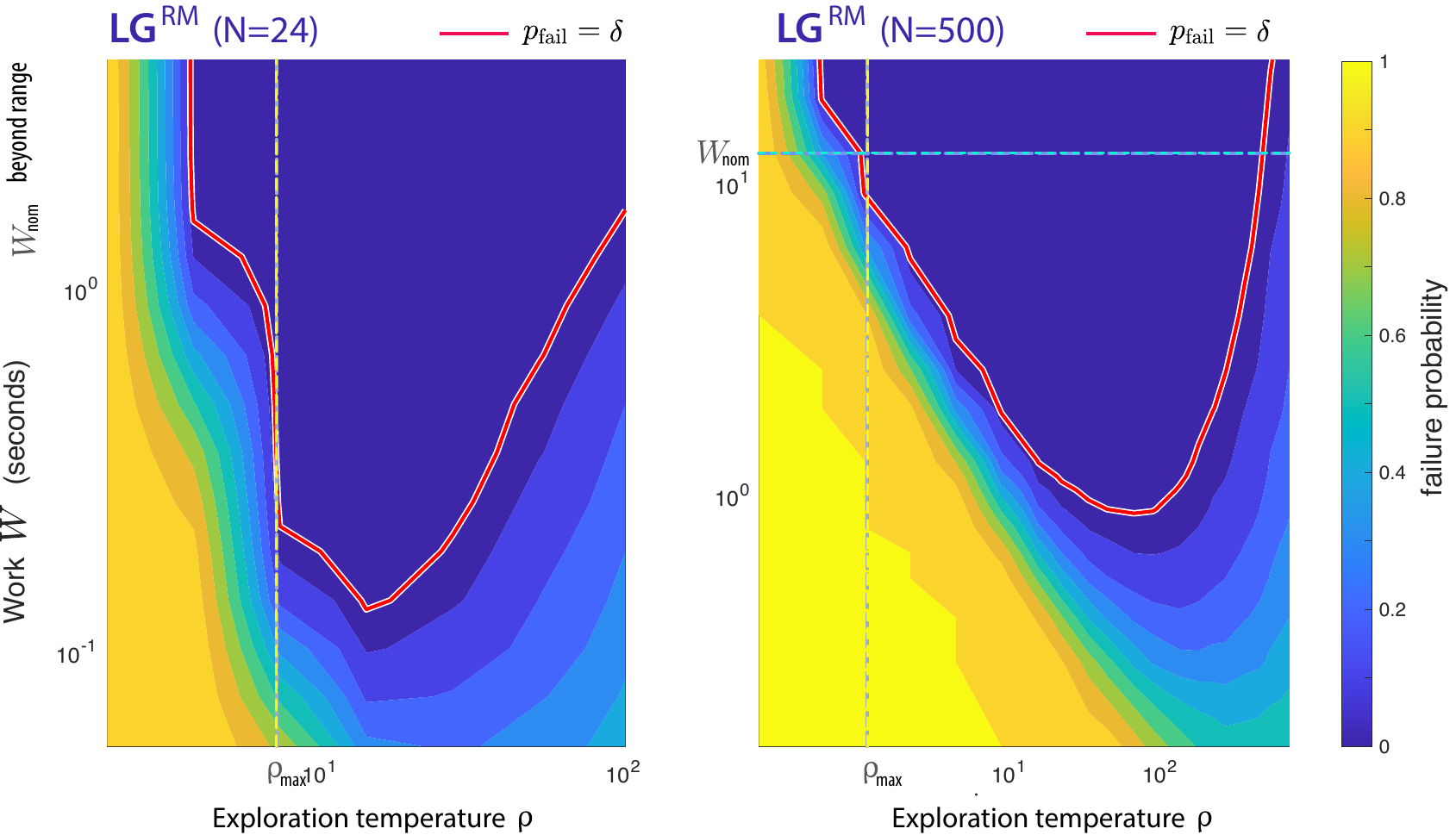}  
	\caption{Empirical performance of  $\LGRM{\temp,N}$ 
	for the Griewank objective in dimension~2. }	\label{f:WorkVsProbFailureGriewank_d2}
\end{figure}

 The  empirical failure probability for $\LGRM{\temp,N}$ using $N=500$ was estimated via Monte Carlo for a range of $(\temp,W)$.   The results   shown  in \Cref{f:WorkVsProbFailureGriewank_d2}  tell us that  
with a computational
work budget of $1$~second, the   PAC constraint is achieved over
the broad temperature range
$
        \temp \in [\temp_-, \temp_+]
$,
with $ \temp_- = \tempLG $ and $ \temp_+ =80$.

This is a striking improvement over the simulated-annealing experiment, in which the
empirical failure probability remains approximately $0.76$ with a computational
work budget of $89.2$~seconds.

\clearpage
 
%%%%%%%%%%%%%%%%%%%%%%%%%%%%%%%%%%%%%%%%%%%%%%
\section{Low Temperature Approximations}
\label{s:low}

\subsection{General theory}

 \wham{Notation for approximations}

 For positive functions $f$ and $g$, write  $    f\sim g$ if and only if $f(r)/g(r)\to 1$ as $r\to\infty$,  and
\[
\begin{aligned}
    f&\logprec g 
    \quad\Longleftrightarrow \   &
    \limsup_{r\to\infty}
    \frac1r\log\frac{f(r)}{g(r)}
   &  \leq0
\\
    f& \logsucc g 
    \quad\Longleftrightarrow \  &
    \liminf_{r\to\infty}
    \frac1r\log\frac{f(r)}{g(r)}
 &    \geq0
\\
    f&\logasyeq g 
    \quad\Longleftrightarrow \    &
       \lim_{r\to\infty}
    \frac1r\log\frac{f(r)}{g(r)}
  &   =0.
\end{aligned}
\] 
When considering vanishing temperature we typically have $r=1/\temp$, and   stress that $\temp\downarrow 0$.

\subsubsection{Spectral theory}

The differential generator is defined for $C^2$ functions $f\colon\Re^d\to\Re$ by
\begin{equation}
  \clL_\temp f
    =
    -\nabla\Obj\cdot\nabla f+\temp \Delta f 
\label{e:diffGen}
\end{equation}  
The corresponding  nonlinear generator of Fleming (e.g.~\cite{fle78a,fen99a,konmey05a}) is then defined for $C^2$ functions $V\colon\Re^d\to\Re$ 
\[
    \clH_\temp(V)
    \eqdef 
    e^{-V}\clL_\temp(e^V).
\]
The multiplicative drift condition {\rm(DV3)} introduced in \cite{konmey05a} is expressed
$$
    \clH_\temp(V)
    \leq
    -\delta W+b\,\ind_S
\eqno{(DV3)}
$$
where   $b<\infty$, and under the assumptions on the diffusion imposed here, it is convenient to assume that the set $S$ is compact and the two functions  $V,W$ are non-negative and coercive.  

We let   $\gap_{\!\temp}$ denote the spectral gap of the differential generator, which under the conditions of this paper defines the second eigenvalue 
$\lambda_2 = -\gap_{\!\temp}$.   The proof of this conclusion and further structure follows from existence of a solution to (DV3).    

However,  it is simplest to work with the extended generator, denoted $\clA_\temp$,
 whose domain allows functions that are not $C^2$    \cite{ethkur05}.
We write $\clA_\temp f = g$ for a function $f\colon\Re^d\to \Re$ if the stochastic process
\[
M_t =   f(X_0) - f(X_t) + \int_0^t  g(X_r)\, dr \,, \quad t\ge 0\,,
\]
is a local martingale with respect to the natural filtration $F_t = \sigma\{ X_r : r\le t\}$.

For a continuous function $v\colon\Re^d\to[1,\infty)$ let $L_\infty^v$ denote the Banach space of measurable functions   $f\colon\Re^d\to\Re$ with finite norm
\[
\| f\|_v = \sup_x \frac{|f(x)|}{v(x)}
\]

Part (i) of \Cref{t:DV3} follows from the  requirement $\vartheta_0^{-1} = \max(2, \temp)$ in the definition of  $V_\temp(x)$.   The remaining conclusions are from \cite{konmey05a,konmey17a}.

\begin{theorem}[Exponential Ergodicity]
\label[theorem]{t:DV3}
If \Cref{a:jacquot} holds then  (DV3) holds for any positive temperature:
  for   fixed $\vartheta\in (0, \vartheta_0)$
  with  $\vartheta_0^{-1} = \max(2, \temp)$, 
   take $V_\temp(x) = \vartheta [ \Obj(x)  - \optObj  ]$,
 $W =  ( \vartheta-\vartheta^2 \temp)   \|\nabla\Obj\|^2 $  for $x\in\Re^d$,  $b =  2 \temp\sup_x \| \Delta \Obj\, (x) \|$,   and $S = \{ x :  W(x) \le s \}$ for a sufficiently large constant $s$.  The following conclusions follow:    

\wham{(i)}
 $L_\infty^{v_\temp} \subset L_2(\pi_\temp)$, with $v_\temp = e^{V_\temp}$.

\wham{(ii)}  
The extended generator $\clA_\temp$ has a discrete spectrum $ \RoA_\temp$, whose elements are real and negative, with the exception of a single eigenvalue at the origin.  
The spectral gap may be expressed $ \gap_{\!\temp} = \min\{ -\lambda :  \lambda\neq 0\,, \ \lambda\in\RoA_\temp \}$.

\wham{(iii)}  
The diffusion is $v_\temp$-uniformly ergodic:  
there exist constants $b_\temp<\infty$ and $\kappa_\temp>0$ such that  for any $f \in L_\infty^v$,
\begin{equation}
\Big| \Expect_x [ f(X_t) ]  - \pi_\temp(f) ]  \Big|  \le  b_\temp  \| f\|_v  v_\temp(x)\exp(-\kappa_\temp t)
    \label{eq:v-uni}
\end{equation}\clThm
\end{theorem}

The theorem does not claim that $\kappa_\temp = \gap_{\!\temp}$ in     \eqref{eq:v-uni}.
This conclusion does hold on a different domain, and using a different bounding constant:

\begin{theorem}[Spectral Theory in $L_2$]
\label[theorem]{t:L2}
If \Cref{a:jacquot} holds then the semigroup  $\{ P_\temp^t : t\ge 0\}$  is a contraction on $L_2(\pi_\temp)$.  Moreover,  for $\temp>0$,

\wham{(i)}  
It   has a discrete spectrum in $L_2(\pi_\temp)$,  and for  $f\in L_2(\pi_\temp)$,
\begin{equation}
    \|P_\temp^t f-\pi_\temp(f)\|_{L_2(\pi_\temp)}
    \leq 
    \|f-\pi_\temp(f)\|_{L_2(\pi_\temp)}   \exp\{-\gap_{\!\temp} t\}  \,,    \qquad t\geq 0  
    \label{eq:L2}
\end{equation}
where 
the spectral gap $\gap_{\!\temp} >0 $ coincides with  the value in \Cref{t:DV3}.

\wham{(ii)}  
There is $C_\temp \colon\Re^d\to \Re_+$ such that for $A\in\clB$ (the Borel sigma field on $\Re^d$),  % [0,\infty]$, finite a.e.~$[\pi_\temp]$, 
 \begin{equation}
    \bigl|P_\temp^t(x,A)- \pi_\temp(A)\bigr|
    \leq     C_\temp(x)\exp\{-\gap_{\!\temp} t\}\,,    \qquad t\geq 0 \, .
    \label{eq:vbound-basic}
\end{equation}

\wham{(iii)}  
The function in (ii) satisfies the following bound, for all sufficiently
small $\temp>0$,
\[
    C_\temp(x)
    \leq
 \temp^{-d/4}
    \exp\left\{\temp^{-1} K(x)\right\}
\]
with $K\colon\Re^d\to\Re_+$ independent of $\temp$. Consequently,
\[
    \log C_\temp(x)=O(1/\temp)  \,, \quad \textit{and} \quad     \log\log\bigl(e+C_\temp(x)\bigr)
    =
    O(\log(1/\temp))
    =
    o(1/\temp).
\] 
 \end{theorem}

\def\Tr{\textsf{Tr}}

In the proof we  require reduction to consideration of densities.
For $s>0$, let $p_\temp^s(x,y)$ denote the density of
$P_\temp^s(x,\cdot)$ with respect to $\pi_\temp$.
Reversibility and the spectral-gap bound imply that the function
$C_\temp$ in \eqref{eq:vbound-basic} may be chosen so that
\[
 C_\temp(x)
 \le
 e^{\gap_{\!\temp} s}
 \left[
 1+\frac12
 \|p_\temp^s(x,\cdot)-1\|_{L_2(\pi_\temp)}
 \right].
\]

Reversibility and the semigroup property give
\[
 \|p_\temp^s(x,\cdot)-1\|_{L_2(\pi_\temp)}^2
 =
 p_\temp^{2s}(x,x)-1 \,, 
\]
and using the spectral expansion of the semigroup,
\[
\int
 \|p_\temp^s(x,\cdot)-1\|_{L_2(\pi_\temp)}^2
 \,\pi_\temp(dx)
 =
 \sum_{j\ge2} e^{2s\lambda_j} \, .
\]
Consequently, $C_\temp$ can be chosen in
$L_2(\pi_\temp)$ whenever the right hand side is finite.

 \smallskip
  \Cref{t:MenozziPesceZhang2021}  implies that for each $x$ and sufficiently small $t>0$, the probability measure $P_\temp^t(x,\cdot)$ has a density that lies in $L_\infty^{v_\temp}$, with $v_\temp$ defined in \Cref{t:DV3}.
The proof of \Cref{t:MenozziPesceZhang2021}
 follows from \cite[Theorem~1.2]{MenozziPesceZhang2021}.

\begin{lemma}
\label[lemma]{t:MenozziPesceZhang2021}
If \Cref{a:jacquot} holds then $P_\temp^t(x,\cdot)$
admits a transition density $q_\temp^t(x,\cdot)$ with respect to
Lebesgue measure, for each $t>0$. Moreover, there exist constants
$C_q<\infty$ and $c_q>0$, independent of $x$, $y$, and $\temp$, such that
for $0<t\le 1$,
\[
\begin{aligned}
q_\temp^t(x,y)
&\le
\frac{C_q}{(\temp t)^{d/2}}
\exp\left\{
    -c_q\frac{\|y-x_t \|^2}{\temp t}
\right\} 
\end{aligned}
\]
where $x_t$ is the solution to the gradient flow \eqref{e:GF} with $x_0=x$.
\clThm
\end{lemma}

\wham{Proof of \Cref{t:L2}}

For (i) first observe that it is shown in \cite{konmey12a} that there is a spectral gap in $L_2(\pi_\temp)$ whenever this is the case in $L_\infty^v$  for a function $v\colon\Re^d \to\Re$.  While stated for Markov chains in discrete time,  the result carries over to continuous time by considering the reslovent kernels as in 
\cite{konmey05a}.  The bound \eqref{eq:L2} is   \cite[Theorem~4.2.5]{bakgenled14}.

For (ii) we begin with \eqref{e:pi_temp}  which tells us that  
the density of $P_\temp^t(x,\cdot)$ with respect to $\pi_\temp$ is 
\[
    p_\temp^s(x,y)
    =
    Z_\temp
    \exp\left\{\temp^{-1}\Obj(y)\right\}
    q_\temp^s(x,y).
\]
Consequently,
\[
\begin{aligned}
    \|p_\temp^s(x,\cdot)\|_{L^2(\pi_\temp)}^2
    &=
    Z_\temp
    \int_{\Re^d}
    q_\temp^s(x,y)^2
    \exp\left\{\temp^{-1}\Obj(y)\right\}dy .
\end{aligned}
\]
Normalize $\optObj=0$ and choose any
$x^\star\in\argmin\Obj$.  Since $\nabla^2\Obj$ is bounded,
there is $M<\infty$ such that
\[
    \Obj(y)
    \leq
    \frac{M}{2}\|y-x^\star\|^2 \,,
\]
and by \Cref{t:MenozziPesceZhang2021},
\[
    q_\temp^s(x,y)^2
    \leq
    \frac{C_q^2}{(\temp s)^d}
    \exp\left\{
        -\frac{2c_q}{\temp s}\|y-x_s\|^2
    \right\}.
\]
Using
\[
    \|y-x_s\|^2
    \geq
    \frac12\|y-x^\star\|^2
    -
    \|x_s-x^\star\|^2,
\]
we obtain
\[
\begin{aligned}
q_\temp^s(x,y)^2
\exp\left\{\temp^{-1}\Obj(y)\right\}
&\leq
\frac{C_q^2}{(\temp s)^d}
\exp\left\{
    \frac{2c_q}{\temp s}
    \|x_s-x^\star\|^2
\right\}
\\
&\quad{}\times
\exp\left\{
    -\frac1{\temp}
    \left(
        \frac{c_q}{s}-\frac{M}{2}
    \right)
    \|y-x^\star\|^2
\right\}.
\end{aligned}
\]
Choose $s>0$ sufficiently small that
\[
    \frac{c_q}{s}>\frac{M}{2}.
\]
Integration then gives for a constant $k_s$,
\[
    \|p_\temp^s(x,\cdot)\|_{L^2(\pi_\temp)}^2
    \leq  w_\temp (x)\,, \qquad 
      w_\temp (x)\eqdef
    k_s Z_\temp\temp^{-d/2}
    \exp\left\{
        \frac{2c_q}{\temp s}
        \|x_s-x^\star\|^2
    \right\}.
\]

For $t\geq s$, the semigroup property and invariance of $\pi_\temp$
give
\[
\begin{aligned}
    P_\temp^t(x,A)-\pi_\temp(A)
    &=
    \int
    \bigl[p_\temp^s(x,y)-1\bigr]
    \bigl[
        P_\temp^{t-s}\ind_A(y)-\pi_\temp(A)
    \bigr]
    \pi_\temp(dy).
\end{aligned}
\]
By   
 Cauchy--Schwarz and    \eqref{eq:L2},
\[
\begin{aligned}
    \bigl|P_\temp^t(x,A)-\pi_\temp(A)\bigr|
    &\leq
    \|p_\temp^s(x,\cdot)-1\|_{L_2(\pi_\temp)}      \|P_\temp^{t-s} \ind_A-\pi_\temp(A)\|_{L_2(\pi_\temp)}
    \\[.5em]    
    &\leq
    \|p_\temp^s(x,\cdot)-1\|_{L_2(\pi_\temp)}  
    \| \ind_A-\pi_\temp(A)\|_{L_2(\pi_\temp)}       \exp\{-\gap_{\!\temp} (t-s)\} 
    \\[.5em]
    &\leq
       C_\temp^+(x) 
    \exp\{-\gap_{\!\temp} t\} \,,
 \\&
 \textit{where}  \ \         C_\temp^+(x) =     e^{\gap_{\!\temp} s}
    \|p_\temp^s(x,\cdot)-1\|_{L_2(\pi_\temp)}          \le  e^{\gap_{\!\temp} s} \Big[ 1 + \sqrt{ w_\temp(x)  } \Big]
\end{aligned}
\]
A bound is also obtained for $t<s$ using  
\[
    \bigl|P_\temp^t(x,A)-\pi_\temp(A)\bigr|
    \leq 1
    \leq e^{\gap_{\!\temp} s}\exp\{-\gap_{\!\temp} t\}.
\]
Part (ii)  therefore follows on defining
$
    C_\temp(x)
    =
    \max\left\{
            e^{\gap_{\!\temp} s} ,\,     C_\temp^+(x)   \right\}
$, and (iii) follows from the definition of $ C_\temp^+(x)$.
    \qed

Consider a measurable set $F$ for which  $\pi_\temp(F^c) \to 0$ as $\temp\downarrow 0$.  
Applying \Cref{eq:vbound-basic},   to achieve $ P_\temp^t(x,F^c) \le s$,  it is sufficient to take $\temp\in (0,1)$ sufficiently small so that $\pi_\temp(F^c)<s$ and then 
  take $t\ge     T_F(\temp;s)$ with      
\begin{equation}
    T_F(\temp;s)
    =
    \frac{1}{\gap_{\!\temp}}
    \log
    \left[
        \frac{C_\temp(x)}
        {s-\pi_\temp(F^c)}
    \right].
\label{e:T_F_temp}
\end{equation}
In applying this formula we must choose the set $F$, along with the  one-run failure probability $s$ to satisfy simultaneously $s>\pi_\temp(F^c)$ and 
\[
    s
    \le
    \begin{cases}
        \delta,
        & \text{for a single trajectory},\\[1mm]
        \delta^{1/N},
        & \text{for $N$ independent trajectories}.
    \end{cases}
\]

We consider next approximations of the spectral gap $\gap_{\!\temp}$ and the 
probability $\pi_\temp(F^c)$ for sets $F$ of interest.
 
%%%   Trace class operators 
%For a positive self-adjoint operator \(T\) on a Hilbert space, if \(\{\phi_j\}\) is any orthonormal basis, its trace is defined by
%
%$$ \Tr(T)=\sum_j \langle \phi_j,T\phi_j\rangle, $$
%
%provided the sum is finite. If \(T\) has discrete eigenvalues \(\{\mu_j\}\), this becomes
%
%$$ \Tr(T)=\sum_j\mu_j. $$
%
%In your setting, since
%
%$$ \Probtemp^t\phi_j=e^{t\lambda_j}\phi_j, \qquad 0=\lambda_1>\lambda_2\ge\lambda_3\ge\cdots, $$
%
%we have, when the sum is finite,
%
%$$ \Tr(\Probtemp^{2s}) = \sum_{j\ge1}e^{2s\lambda_j} = 1+\sum_{j\ge2}e^{2s\lambda_j}. $$
%
%The connection with your density is then a theorem/identity for an integral operator:
%
%$$ \Tr(\Probtemp^{2s}) = \int p_\temp^{2s}(x,x)\,\pi_\temp(dx), $$
%
%under the trace-class condition.
%
%But for your proof, there is an even cleaner route that avoids this terminology entirely. Start from
%
%$$ \|p_\temp^s(x,\cdot)-1\|_2^2 = p_\temp^{2s}(x,x)-1. $$
%
%Then integrate with respect to \(\pi_\temp(dx)\). Using the spectral expansion
%
%$$ p_\temp^s(x,y) = 1+\sum_{j\ge2}e^{s\lambda_j}\phi_j(x)\phi_j(y), $$
%

\subsubsection{Low temperature approximations}
\label{s:EK+}

Consideration of the choice of parameters to achieve  \eqref{e:performance-objective} concerns approximation of the probability $ \Prob_x\{X_T\in G\} $ for several different sets $G\in\clB$.    For this we apply exponential ergodicity,  which reduces the problem to approximation of terms  \eqref{eq:vbound-basic},
including the spectral gap $\gap_{\!\temp}$ and $\pi_\temp(G)$.

%%%%%
%One element in our analysis is   the following consequence of  the Laplace principle
%(see, e.g., \cite[Chap.~4]{deBruijn1970} or \cite{Hwang1980}),

The analysis requires two forms of the Laplace approximation.
The first is the logarithmic approximation used throughout the
low-temperature analysis: probabilities of sets are determined, at
exponential scale, by the minimum value of the objective over the set.
For the high-dimensional discussion in \Cref{s:highd} we also require a
sharper local approximation near the global minimizers.  When the
global minima are isolated and non-degenerate, this implies a
Gaussian approximation to the Gibbs measure $\pi_\temp$,
 and a 
$\chi_d^2$ approximation for   $ \Obj(X_\temp)-\optObj $ with
$X_\temp\sim\pi_\temp$.

The following proposition summarizes the results required in the following.
Part~(i) is the standard Laplace principle; the remaining conclusions
give refinements under nondegeneracy of the global minima.
See, for example, \cite[Chap.~4]{deBruijn1970} and \cite{Hwang1980}.

\begin{theorem}[Laplace approximations for Gibbs measures]
\label[theorem]{t:Gmin_low_temp}
Suppose that \Cref{a:jacquot} holds, and let
\[
    \pi_\temp(dx)
    =
    Z_\temp^{-1}
    \exp\{-\Obj(x)/\temp\}\,dx .
\]
\wham{(i)}
The normalizing constant satisfies
$
    \lim_{\temp\downarrow0}
    \temp\log Z_\temp
    =
    -\optObj $.
More generally, 
\[
\begin{aligned}
     \lim_{\temp\downarrow0}
    \temp\log
    \int_G
    \exp\{-\Obj(x)/\temp\}\,dx
   & =
    -\inf_{x\in G}\Obj(x)
    \\
        \lim_{\temp\downarrow0}
    \temp\log\pi_\temp(G)
&    =
    -
    \inf_{x\in G}
    [\Obj(x)-\optObj] \,, 
\end{aligned}
\]
for any  $G\in\clB$ satisfying $\displaystyle 
    \inf_{x\in G^\circ}\Obj(x)
    =
    \inf_{x\in G}\Obj(x)$.

\smallskip
 
In the remaining statements suppose that the set of global minimizers is finite,
$
    \argmin\Obj
    =
    \{x^\star_1,\ldots,x^\star_m\}$, 
and each minimum is non-degenerate:
\[
    H_j
    \eqdef
    \nabla^2\Obj(x^\star_j)>0,
    \qquad
    1\leq j\leq m.
\]

\wham{(ii)} 
The sharper Laplace approximation holds:
\[
\begin{aligned} 
    \lim_{\temp\downarrow0}
    \temp^{-d/2}e^{\optObj/\temp} &Z_\temp
    =
    (2\pi)^{d/2}K_\star \,, 
\\
&\textit{where} \ \ 
    K_\star
    \eqdef
    \sum_{j=1}^m
    \det(H_j)^{-1/2}
\end{aligned}
\]
That is,
$\displaystyle     Z_\temp  \sim  
    e^{-\optObj/\temp} 
    (2\pi\temp)^{d/2}K_\star$,   $    \temp\downarrow0$.

%Then the sharper Laplace approximation holds:
%\[
%    Z_\temp
%    \sim
%    e^{-\optObj/\temp}
%    (2\pi\temp)^{d/2}K_\star,
%    \qquad
%    \temp\downarrow0.
%\]

\wham{(iii)}
If $X_\temp\sim\pi_\temp$, then normalized excess cost converges in distribution:
\begin{equation}
    \frac{\Obj(X_\temp)-\optObj}{\temp}
    \ \darrow\
    \frac12\chi_d^2,
    \qquad
    \temp\downarrow0.
    \label{e:Gibbs-energy-chi-square}
\end{equation}
Consequently, for 
  $\eta(\temp)\downarrow0$ satisfying
$
    \sfrac{\eta(\temp)}{\temp}\longrightarrow a
    \in[0,\infty]$, 
we have
\begin{equation}
    \pi_\temp( \RoA_{\eta(\temp)})
    \longrightarrow
    \Prob\{\chi_d^2\leq2a\},
    \label{e:Gibbs-sublevel-chi-square}
\end{equation}

\wham{(iv)}
In the small-target regime
\[
    \eta(\temp)\downarrow0,
    \qquad
    \frac{\eta(\temp)}{\temp}\longrightarrow0,
\]
the approximation in part~(iii) can be sharpened to
\begin{equation}
    \pi_\temp(\RoA_{\eta(\temp)})
    \sim
    \frac{1}{\Gammafn(d/2+1)}
    \left(\frac{\eta(\temp)}{\temp}\right)^{d/2}.
    \label{e:Gibbs-small-target-joint}
\end{equation}
The leading constant is independent of the Hessians
$\{H_j\}$ and of the number of global minimizers.
 \clThm
\end{theorem}

%
%\begin{lemma}
%\label[lemma]{t:Gmin_low_temp}  
%The following hold under \Cref{a:jacquot}:
%\wham{(i)} $ \displaystyle 
%    \lim_{\temp \downarrow 0}
%    \temp\log Z_{\temp}
%    =
%    -\optObj $.
%    
%    \wham{(ii)} 
%Let $G\in\clB(\Re^d)$ satisfy
%$ \displaystyle 
%    \inf_{x\in G^\circ}\Obj(x)
%    =
%    \inf_{x\in G}\Obj(x)$.
%Then,
%\[
%    \lim_{\temp \downarrow 0}
%    \temp
%    \log\int_G
%    \exp\left\{-\frac{1}{\temp}\Obj(x)\right\}\,dx
%    =
%    -\inf_{x\in G}\Obj(x).
%\]
%Consequently,
%$ \displaystyle 
%    \lim_{\temp\downarrow0}
%    \temp\log\pi_\temp(G)
%    =
%    -\inf_{x\in G}[\Obj(x)-\optObj]
%$,
%or equivalently,
%\[
%    \pi_\temp(G)
%    \logasyeq
%    \exp\left\{
%        -\frac{1}{\temp}
%        \inf_{x\in G}[\Obj(x)-\optObj]
%    \right\}.
%\]
%\clThm
%\end{lemma}

   \Cref{t:Gmin_low_temp}~(i) implies the following,
\begin{equation}
\begin{aligned}
q_\mu(\temp) & \eqdef \pi_\mu(\clS_\mu^c)           \logasyeq
    \exp\{-\temp^{-1} \mu\},  \qquad \mu > 0
   \\
q_{\RoA}(\temp)   &    \eqdef     \pi_\temp(\RoA^c)     \logasyeq    \exp\{-\eta/\temp\} \,, 
\end{aligned}
\label{e:qApprox}
\end{equation}
    where  $
        \RoA_\mu
    =
    \{x:\Obj(x) \le \optObj+\mu\} $  
 and  $ \RoA$ appears in  \Cref{a:jacquot}~(iii).  These  
 approximations   will be combined with  \eqref{e:T_F_temp} to establish complexity bounds:

Approximations of   the spectral gap are based on $\Dbar$ defined in \eqref{e:EnergyBarrier}.
The proof of the following may be found in
Jacquot~\cite{Jacquot1992}; see also
Miclo~\cite{Miclo1992}.

\begin{theorem}[Logarithmic spectral gap asymptotics]
\label[theorem]{t:low-temp-basic}
Under \Cref{a:jacquot} the spectral gap admits the approximation,  
\begin{equation}
    \gap_{\!\temp}
    \logasyeq   
    \exp\{-\Dbar/\temp\} \,, \quad \temp \downarrow 0 
     \label{e:low-temp-master_lambda}
 \end{equation}
\clThm
 \end{theorem}

\subsection{Simulated annealing}

We consider
\begin{equation}
    dX_t=-\nabla\Obj(X_t)\,dt+\sqrt{2\temp(t)}\,dB_t ,
    \label{eq:SimAnn}
\end{equation}
with logarithmic cooling   \eqref{e:cooling}.   
 
In much of the literature approximations are based on not on $\Dbar$, but a potentially smaller value
\[
    \Dbar_{\rm SA}
    =
    \sup_{\Obj(x)>\optObj} D(x)
    =
    \sup_{\Obj(x)>\optObj}
    \inf_{\Obj(y)<\Obj(x)}
    \{H(x,y)-\Obj(x)\} \, .
\]
This is the maximum, over initial states $x$, of the
minimum increase in objective value required along a continuous path
from $x$ to some state $y$ satisfying $\Obj(y)<\Obj(x)$.   It is the quantity  
appearing in early analysis of simulated annealing \cite{haj88}.

However,  it seems that upper and lower bounds on the complexity in terms of total simulation time require $    \Dbar_{\rm SA} =     \Dbar$.    
Under the assumptions imposed in \cite{TangZhou2023}, the two
quantities coincide:
\[
    \Dbar_{\rm SA}=\Dbar=E^*,
\]
where $E^*$ is the \textit{critical depth} defined in that reference.   Also assumed there is uniqueness of the minimizer  $x^\star$ for the objective, and $\nabla^2 \Obj\, (x^\star) >0$.  Rather than list all of the assumptions here we defer to this recent published work.

The following result provides the approximation $T_1 \approx \exp( E L/\epsy)$ for  fixed $\epsy>0$ provided  $E > \Dbar+2\epsy$
 in the cooling schedule \eqref{e:cooling}.

\begin{proposition}
\label[proposition]{t:SimAnn}
Suppose that $\Obj$ satisfies the assumptions of \cite{TangZhou2023}. 
For $\epsy>0$ fixed, choose  $E$ in  \eqref{e:cooling}  to satisfy $E > \Dbar+2\epsy$.
Then, with $L = \log(1/\delta)$, the time $T_1$ to achieve     \eqref{e:performance-objective} satisfies
\[
    \lim_{\delta\downarrow0}
    \frac{1}
         {L}  \log T_1(\epsy,\delta;E)
    =
    \frac{E}{\epsy}.
\]
\end{proposition}

\PF  By \cite[Theorem~1]{TangZhou2023},
for every $\zeta>0$ there is a finite constant
$C=C(\epsy,\zeta,E)$ such that
\[
    \Prob\{\Obj(X_t)-\optObj>\epsy\}
    \leq
    C t^{-m(E,\epsy)+\zeta},
\]
where
\[
    m(E,\epsy)
    =
    \min\left\{
        \frac{\epsy}{E},
        \frac12\left(1-\frac{\Dbar}{E}\right)
    \right\} = \frac{\epsy}{E} \,,
\]
where the second equality follows from the assumed bound  $E > \Dbar+2\epsy$.
 
Hence, for any $0<\zeta<\epsy/E$, it is sufficient to take
\[
    T
    \geq
    \left( {C} /{\delta}\right)^{
        1/(\epsy/E-\zeta)} .
\]
It follows that
\[
    \limsup_{\delta\downarrow0}
    \frac{1}{L} \log T_1(\epsy,\delta;E)
    \leq
    \frac{1}{\epsy/E-\zeta}.
\]
Letting $\zeta\downarrow0$ gives
\[
    \limsup_{\delta\downarrow0}
    \frac{1}{L} \log T_1(\epsy,\delta;E)
    \leq
    \frac{E}{\epsy}.
\]
The proof will be completed after we establish 
 the reverse inequality.
 
 For this
  we apply a result of  \cite{Marquez1997}:  for
sufficiently small fixed $\epsy>0$,
\[
    \lim_{t\to\infty}
    \frac{1}{\log t}
    \log
    \Prob\{\Obj(X_t)-\optObj>\epsy\}
    =
    -\frac{\epsy}{E}.
\]
Hence, for every $\zeta>0$ and all sufficiently large $t$,
\[
    \Prob\{\Obj(X_t)-\optObj>\epsy\}
    \geq
    t^{-\epsy/E-\zeta}.
\]
Consequently, the requirement
$
    \Prob\{\Obj(X_T)-\optObj>\epsy\}\leq\delta
$
implies
\[
    \liminf_{\delta\downarrow0}
    \frac{1}{L} \log T_1(\epsy,\delta;E)
    \geq
    \frac{1}{\epsy/E+\zeta}.
\]
Since $\zeta >0$ is arbitrary we   obtain the required lower bound,
\[
    \lim_{\delta\downarrow0}
    \frac{1}{L} \log T_1(\epsy,\delta;E)
\ge
    \frac{E}{\epsy}.
\]
\qed

%%%%%%%%%%%%%%%%%%%%%%%%%%%%%
\subsection{Fixed-temperature Langevin}
%%%%%%%%%%%%%%%%%%%%%%%%%%%%%

In consideration of $\mathrm{L}_{\temp}$ we wish to choose $\temp$ in order to minimize the run time $T_2$ over $\temp$ while satisfying the PAC performance requirement     \eqref{e:performance-objective}.
    
     \begin{proposition}
\label[proposition]{t:2Ltemp}
Consider $\mathrm{L}_{\temp}$ under  \Cref{a:jacquot}, with $E>\Dbar$ fixed,  and with  $\temp_2$ given in 
\eqref{e:23temp}.   We then have,   
 \[
    T_2 \logasyeq  \exp\left\{ \frac{ EL + o(L)   }{ \epsy} \right\} 
 \,, \quad \delta\downarrow 0
\]
\end{proposition}

\PF
We apply \eqref{e:T_F_temp} with $F = \RoA_\epsy$, $s =\delta$, to obtain 
\[
T_2 = 
    T_F(\temp;\delta)
    =
    \frac{1}{\gap_{\!\temp}}
    \log
    \left[
        \frac{C_\temp(x)}
        {\delta - q_\epsy(\temp_2) }
    \right] \,,
\]

We have by assumption,
\[
    \frac{\epsy}{\temp_2}
    =
    \frac{E}{\Dbar}L .
\]
Applying \eqref{e:qApprox} we obtain  for fixed $\epsy>0$, as $\delta\downarrow0$,
\[
q_\epsy(\temp_2) \eqdef
   \pi_{\temp_2} (\RoA_\epsy^c) 
    \logasyeq
    \exp\left\{
        -\frac{E}{\Dbar}L 
    \right\}.
\]
Application of    \eqref{e:low-temp-master_lambda} then establishes the desired approximation:
\[
\begin{aligned}
T_2  \logasyeq      \exp\left\{ \frac{\Dbar}{\temp_2} \right\} 
    \log
    \left[
        \frac{C_\temp(x)}
        {\delta - q_\epsy(\temp_2) }
    \right] 
    &
    \logasyeq      \exp\left\{ \frac{ EL + o(L)   }{ \epsy} \right\} 
    \log
    \left[
        \frac{ C_\temp(x)}
        { (1- o(1)) \delta  }
    \right] 
    \\ 
    &
    \logasyeq      \exp\left\{ \frac{ EL + o(L)   }{ \epsy} \right\} 
\end{aligned}
\] 
\qed

\subsection{Parallel Langevin}

In 
$\mathrm{L}_{\temp,N}$ we obtain  $N$ independent copies of \eqref{eq:langevin}, all with the same $\temp$ and run length $T$, and choose
\[
    \haX_T
    \in \argmin_{1\leq i\leq N} \Obj(X_T^{(i)}).
\]
Denoting the one-run failure probability is $Q_T(\temp)$,   independence gives
\[
    \Prob\{\Obj(\haX_T)>\epsy\}
    =Q_T(\temp)^N.
\]
Our goal is to choose $T_3$ so that $Q_T(\temp)^N \le \delta$ when $T=T_3$.

     \begin{proposition}
\label[proposition]{t:3LtempN}

Consider $\mathrm{L}_{\temp,N}$ under  \Cref{a:jacquot}, with $E>\Dbar$ fixed,  and with  $\temp_3$ given in 
\eqref{e:23temp}.  Then,  for fixed $\epsy>0$, 
\begin{equation}
    T_3 \logasyeq
       \exp\left\{  \frac{EL}{N\epsy}   \right\} 
 \,, \quad \temp_3\downarrow 0 
\label{e:T3-E}
\end{equation}
% and   the total simulation time is  
%\[
%    \clC^{\LA{\temp,N}}
%    \logasyeq
%    N\exp\left\{
%        \frac{EL}{N\epsy}
%    \right\}.
%\]
\end{proposition}

\PF
We have
\[
    \frac{\epsy}{\temp_3}
    =
    \frac{E}{\Dbar}\frac{L}{N},
\]
so that from \eqref{e:qApprox} we obtain, as $    \temp_3\downarrow0$,
\[
q_\epsy(\temp_3) \eqdef
   \pi_{\temp_3} (\RoA_\epsy^c) 
\logasyeq
    \exp\left\{
        -\frac{E}{\Dbar}\frac{L}{N}
    \right\}     =
    o\left(\delta^{1/N}\right) 
\] 
Moreover,
\[
    \frac{\Dbar}{\temp_3}
    =
    \frac{EL}{N\epsy},
\]
so that \eqref{e:T3-E}
is obtained on applying    \eqref{e:low-temp-master_lambda}.
 \qed

\subsection{Parallel Langevin--gradient}

The goal of   $\LGFE{\temp,N}$ and  $\LGRM{\temp,N}$ is to  produce at least one terminal
point in $\RoA$; after the switching time $\Tf^0$, deterministic gradient
descent supplies the final accuracy.  The time is required to satisfy
\[
\Prob\{ X_T \not\in\RoA \} \le \delta^{1/N}  \,, \quad T = \Tf^0
\]
Recall $\tempLG$ was defined   in 
\eqref{e:temp_nom}.   In the following we choose $\temp_\delta = a \tempLG$ with $a=a_\delta \approx 1$ a function of $\delta$.
An example satisfying the required constraints is  $a_\delta = 1 - 1/L$.

     \begin{proposition}
\label[proposition]{t:4LGtempN}
Consider $\LGRM{\temp,N}$ under  \Cref{a:jacquot}, with $E>\Dbar$ and $N\ge 1$  fixed.
Choose the temperature  $\temp_\delta = a \tempLG$, with scaling satisfying  $a = a_\delta\uparrow 1$ 
and $\delta^a \to 0$ as $\delta\downarrow 0$. 
Then,
 as $\delta\downarrow 0$,
 \[
    \Tf^0 \logasyeq     
       \exp\left\{ \frac{EL}{N \eta}    \right\}  
\]
\end{proposition}

\PF
This is again similar to the previous proof, beginning with an application of   \eqref{e:T_F_temp} with $s =\delta^{1/N}$.
However, we now use   $F =\RoA$, giving 
\[
\Tf =  
    \frac{1}{\gap_{\!\temp}}
    \log
    \left[
        \frac{C_\temp(x)}
        {\delta^{1/N} -    \pi_\temp(\RoA^c)  }        \right]    
\]
Applying  \eqref{e:qApprox} and the relation $    \sfrac{\eta}{\temp}
    =  L/(aN)$
gives  
\[
 \pi_\temp(\RoA^c)   \logasyeq    \exp\{-\eta/\temp\}   =   \exp\left\{
        -\frac{1}{a}\frac{L}{N}
    \right\}     =
    o\left(\delta^{1/N}\right) 
\]
giving
\[
\Tf \logasyeq \frac{1}{\gap_{\!\temp}}      \log
    \left[
        \frac{C_\temp(x)}
        {\delta^{1/N}   }        \right]       
\]
The logarithmic term is of order $O(L)$ for fixed $N$ while the spectral gap $\gap_{\!\temp}$  grows at an exponential rate with $L$.     Consequently,  only the spectral gap is significant in a logarithmic scale approximation of the required exploration time:
\[
\Tf \logasyeq \frac{1}{\gap_{\!\temp}}            \logasyeq \frac{1}{\gap_{\!\tempLG}}  \logasyeq   \exp\left\{ \frac{\Dbar}{\tempLG}\right\} \,,
\]
where the final approximation follows from 
    \eqref{e:low-temp-master_lambda}.
    This completes the proof on substituting the expression for $\tempLG$ in \eqref{e:temp_nom}.   
 \qed

\section{Conclusions}

We have compared Langevin-based approaches to
$(\epsy,\delta)$-PAC global optimization, with emphasis on   total simulation work.  The main conclusion is that
independent restart can greatly reduce the cost of high-confidence
global exploration, and that a further improvement is obtained by
separating exploration from exploitation.  In the
Langevin--gradient schemes, stochastic dynamics are used only to
locate a region of attraction of a global minimizer, after which
gradient flow supplies the final accuracy at cost $O(\log (1/\epsy))$. 

The analysis therefore gives a concrete answer to the exploration–exploitation question posed at the outset. Exploration should continue only long enough—and at sufficiently high temperature—to place at least one trajectory in a region from which deterministic dynamics are reliable. Beyond that point, additional stochastic exploration is wasteful: terminal accuracy is obtained far more efficiently by exploitation.

The analysis is based primarily on low-temperature approximations,
which make explicit the roles of the dominant energy barrier
$\Dbar$, the attraction margin $\eta$, and the confidence parameter
$L=\log(1/\delta)$.   It will be of interest in future work to exploit more fully moderate-temperature theory,
as it is well known that low-temperature carries an additional computational cost (see e.g.\ \cite{tanwuzho24}).    

More broadly, improved exploration need not come only from temperature
selection.  Other considerations include correlated excitation  \cite{orvkerprobacluc22,laumey25b} or
deterministic   excitation  \cite{laumey22e}.

\addcontentsline{toc}{section}{References}

\bibliographystyle{abbrv}
%\bibliography{strings,markov,q,SimAnnealing} 
 \def\cprime{$'$}\def\cprime{$'$}

 \clearpage
 
 \appendix

 \section{Rastrigin objective}
\label{s:RastriginAppendix}

Recall that
\[
    \Obj(x)=\sum_{i=1}^d g(x_i),
    \qquad
    g(z)=10+z^2-10\cos(2\pi z).
\]
The function $g$ is even, nonnegative, and has its unique global
minimum at the origin.

\wham{Proof of \Cref{t:Rastrigin}}
We first consider the one-dimensional objective $g$, whose critical
points satisfy
\[
    g'(z)=2z+20\pi\sin(2\pi z)=0.
\]
Observe that $|20\pi\sin(2\pi z)|\leq20\pi$, and hence  the equation $g'(z)=0$ has no
solutions for $|z|>10\pi$.

Let $s$ denote the first positive local maximum and $m$ the first
positive nonglobal local minimum.   These scalars and their evaluations can be computed numerically:
\[
\begin{aligned}
    s&=0.5025460366  &  \quad    m &=0.9949586377 
\\
    g(s)&=20.2512729910&  \quad  
        g(m) &=0.9949590571
 \end{aligned}
\]

The communication height in one dimension is particularly simple: since every continuous path from $a$ to $b$ must traverse the
interval $[a,b]$, we have for any $a<b$,
\[
    H_1(a,b)
    =
    \max_{z\in[a,b]} g(z),
\]
 Consequently, specializing the definition
\eqref{e:EnergyBarrier} to the one-dimensional objective $g$, and
using $\min g=0$,  
\begin{equation}
    \Dbar_1
    =
    \sup_{a,b\in\Re}
    \bigl\{
        H_1(a,b)-g(a)-g(b)
    \bigr\}.
    \label{e:RastriginBarrier1Ddef}
\end{equation}

\begin{figure}[h!]
	\centering
	\includegraphics[width=0.5\hsize]{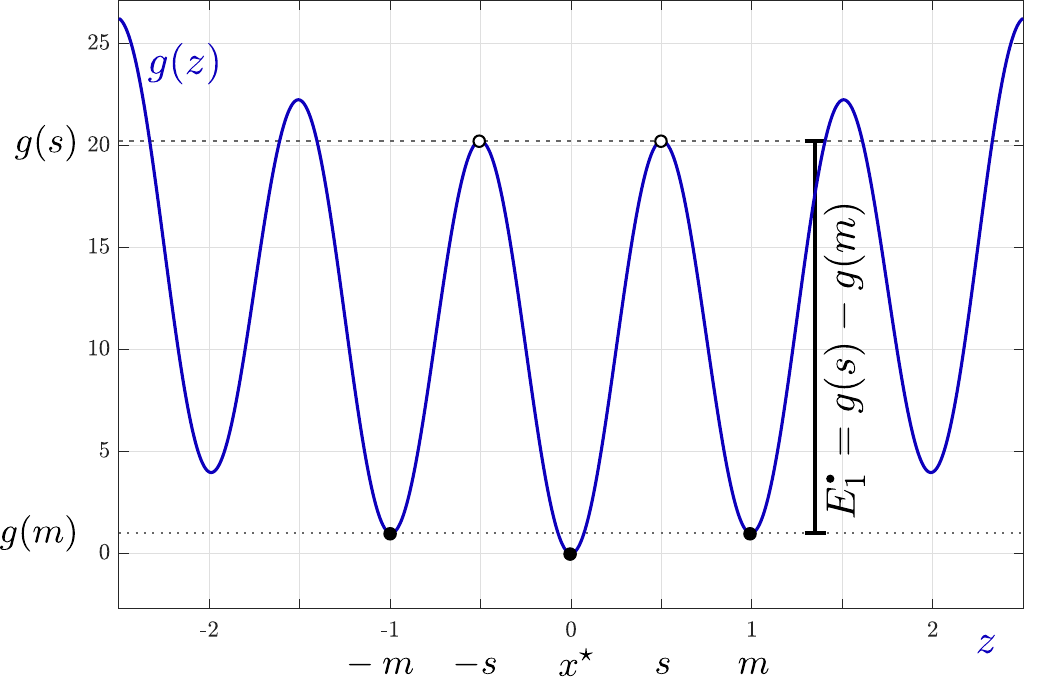}  
 \caption{One-dimensional Rastrigin potential
$g(z)=10+z^2-10\cos(2\pi z)$, with  dominant one-dimensional barrier
is $\Dbar_1=g(s)-g(m)$.  The global minimum is at $z=0$, the
first neighboring local minima are at $z=\pm m$, and the intervening
local maxima are at $z=\pm s$.   }

		\label{f:Rastrigin1D}
\end{figure}

The one-dimensional geometry is illustrated in
\Cref{f:Rastrigin1D}.  For a pair of wells, the expression in
\eqref{e:RastriginBarrier1Ddef} is largest when the endpoints are
taken at the bottoms of the wells.  The first neighboring local
minima $\pm m$ are the deepest nonglobal minima, and a path from
either of these minima to the global minimum at $0$ must cross the
intervening local maximum at $\pm s$.  Inspection of the remaining
local minima shows that none gives a larger value in
\eqref{e:RastriginBarrier1Ddef}.  Hence the supremum is attained for
the pair $(0,m)$, or symmetrically $(0,-m)$, and
\begin{equation}
    \Dbar_1
    =
    H_1(0,m)-g(0)-g(m)
    =
    g(s)-g(m)
    =
    19.2563.
    \label{e:RastriginBarrier1D}
\end{equation}

We next show that the same barrier is obtained in every dimension.
Let $\Dbar_d$ denote the dominant energy barrier for $\Obj$ on
$\Re^d$.  For $x,y\in\Re^d$, construct a path from $x$ to $y$ by
changing the coordinates one at a time.  During the change of
coordinate $i$, the maximal objective is bounded by
\[
    H_1(x_i,y_i)
    +
    \sum_{j<i}g(y_j)
    +
    \sum_{j>i}g(x_j).
\]
Since
\[
    H_1(x_i,y_i)
    \leq
    g(x_i)+g(y_i)+\Dbar_1
\]
and $g\geq0$, each of these maxima is at most
\[
    \Obj(x)+\Obj(y)+\Dbar_1.
\]
It follows that
\[
    \Dbar_d\leq\Dbar_1.
\]
The reverse inequality follows by considering points that differ in
only one coordinate, with all remaining coordinates equal to zero.
Hence
 $   \Dbar_d=\Dbar_1 
$
for every $d\geq1$.

It remains to verify the claim concerning the level set  $   \RoA_\eta=\{x:\Obj(x)\le\eta\}$.  

The function $g$ is strongly convex on the closed interval $ \{ z : g(z) \le \eta \}$  provided $0<\eta <g(m) =0.9949590571$.
This is clear from the figure and can be verified from the second derivative formula 
$  g''(z)
    =
    2+40\pi^2\cos(2\pi z)$.

Fix a value of  $\eta$ satisfying this upper bound,   and let $a_\eta$ denote the unique positive solution to $ g(a_\eta)=\eta$.  
We then have the inclusion $
    \RoA_\eta\subset(-a_\eta,a_\eta)^d$.
 For $\eta\leq0.99$ we have
\[
    a_\eta\leq a_{0.99}=0.0712307811,
\]
and on this interval
$
    g''(z) \ge   g''(a_{0.99} )  \ge \underline\mu \eqdef    357    $, 
%    2+40*pi^2*cos(2*pi*0.0712307811)   =  357.9008
and hence for $x\in\RoA_\eta$,
\[
    \nabla^2\Obj(x)\succeq\underline\mu I  
\]
Strong convexity on this set gives the Polyak--Lojasiewicz bound
\[
    \|\nabla\Obj(x)\|^2
    \geq
    2\underline\mu\,\Obj(x),
    \qquad x\in\RoA_\eta.
\]
Since $\Obj$ is nonincreasing along gradient flow,
$\RoA_\eta$ is forward invariant, and therefore
\[
    \Obj(x_t)
    \leq
    \Obj(x_0)e^{-2\underline\mu t}
    \leq
    \eta e^{-2\underline\mu t}.
\]
This verifies the requirements of \Cref{a:jacquot}(iii).

We now turn to (iii).   The positive spectral gap follows from   \Cref{t:L2}.    More structure regarding the eigenvalues follows from 
the special form of the objective, $    \Obj(x)=\sum_{i=1}^d g(x_i)$.
It follows that the  Gibbs distribution is product form, 
$
    \pi_{\temp,d}
    =
    \pi_{\temp,1}^{\otimes d}
$,
and the generator may be expressed
\[
    \clL_{\temp,d}
    =
    \sum_{i=1}^d \clL_{\temp}^{(i)},
\]
where each $\clL_{\temp}^{(i)}$ is a copy of the one-dimensional
generator acting on coordinate $i$.

Let
$
    0=\upgamma_1>\upgamma_2>\upgamma_3>\cdots
$
denote the one-dimensional eigenvalues, with corresponding
eigenfunctions $\{\phi_k\}$.  These eigenvalues are simple by  
one-dimensional Sturm--Liouville theory; see, for example,
\cite[Chap.~9]{teschl14}.   The eigenfunctions in dimension $d$ are  
products
\[
    \prod_{i=1}^d\phi_{k_i}(x_i),
\quad \textit{
with eigenvalues} 
\quad
    \sum_{i=1}^d\upgamma_{k_i}.
\]
It follows immediately that the spectral gap for the $d$-dimensional generator is 
$\gap_{\!\temp} = -\upgamma_2$, independent of $d$, and that it has multiplicity $d$:
the associated eigenspace is
\[
\clE_\temp = 
    \spn\{
        \phi_2(x_1),\ldots,\phi_2(x_d)
    \}.
\]
Thus the spectral gap   is independent of dimension.

For $d=2$ this gives
$
    \lambda_2=\lambda_3=-\gap_{\!\temp}$ and $\clE_\temp
    =
    \spn\{\phi_2(x_1),\phi_2(x_2)\}$.  
Since $g$ is even, $\phi_2$ may be chosen odd.  If $T$ denotes rotation
through $90^\circ$, then $\clE_\temp$ is invariant under
$h\mapsto h\circ T$, and an orthogonal pair may be chosen so that
\[
    h_3(x)=h_2(Tx).
\]
This explains the symmetry observed in \Cref{f:SecondEigenfunctions_rastrigin_rho_1_v5}.
\qed

\medskip

There is an important distinction from the six-hump camel example.
When the spectral gap is simple, the associated eigenfunction is
unique up to multiplication by a nonzero scalar.  Consequently, its
zero set, and hence the pair of nodal domains $H_+$ and $H_-$, is
uniquely determined up to interchange of the two domains (made possible because $h_2$ and $-h_2$ are each eigenfunctions).

For Rastrigin in dimension $d>1$, the eigenvalue associated with the
spectral gap has multiplicity $d$.  The corresponding eigenspace
$\clE_\temp$ is uniquely determined by the generator, but there is no
distinguished eigenfunction within this space.  Different linear
combinations of eigenfunctions in $\clE_\temp$ generally have
different zero sets.  Consequently, there is no uniquely determined
two-way partition of the state space associated with the spectral gap.

\section{Impact of dimension}

\wham{Proof of \Cref{t:running-min-temp-scaling}}

Since $d>2$, the choice \eqref{e:rho-epsilon-scaling} gives
\[
    \temp_{\epsy}\downarrow0
\quad \textit{
and}
\quad
    \frac{\epsy}{\temp_{\epsy}}
    =
    \frac{1}{\kappa}\epsy^{2/d}
    \longrightarrow0.
\]
Hence the small-target conclusion of
\Cref{t:Gmin_low_temp} applies with
$r_{\temp}=\epsy$, giving
\[
    \pi_{\temp_{\epsy}}
    \{\RoA_\epsy\}
    \sim
    \frac{1}{\Gammafn(d/2+1)}
    \left(
        \frac{\epsy}{\temp_{\epsy}}
    \right)^{d/2}.
\]
Substituting $
    \temp_{\epsy}
    =
    \kappa\epsy^{1-2/d}$
gives
\[
    \left(
         {\epsy}/{\temp_{\epsy}}
    \right)^{d/2}
    =
    \kappa^{-d/2}\epsy,
\]
and therefore recalling $K_d$ defined in \eqref{e:Kd_kappa}, 
\[
    \pi_{\temp_{\epsy}}
    \{\RoA_\epsy\}
    \sim
    \frac{\epsy}
    {\Gammafn(d/2+1)\kappa^{d/2}} =     \frac{\epsy}{K_d(\kappa)}
\]
which
 establishes 
    \eqref{e:Gibbs-target-temp-scaled}.

Write
$
    p_{\epsy}
    \eqdef
    \pi_{\temp_{\epsy}}
    \{\RoA_\epsy\}$ so that  $
    \Prob\{M_n>\epsy\}
    =
    (1-p_{\epsy})^n$.  
Then, from the definition \eqref{e:nstarIdeal},
\[
    n^\star(\epsy,\delta)
    =
    \left\lceil
        \frac{L}{-\log(1-p_{\epsy})}
    \right\rceil.
\]
Since $p_{\epsy}\to0$,
\[
    -\log(1-p_{\epsy})
    \sim
    p_{\epsy}
    \sim
    \frac{\epsy} {
K_d(\kappa)}
\]
It follows that
$    n^\star(\epsy,\delta)
    \sim
    K_d(\kappa) L/\epsy$, 
which proves \eqref{e:n-epsilon-delta}.
\qed

\section{Calibrated computational-work model}
\label{s:WorkModel}

The primitive costs used in the current simulations are
\[
\begin{array}{c|c}
\text{primitive operation} & \text{measured cost (seconds)}\\ \hline
\text{gradient evaluation}
    &c_\nabla=4.4198\times10^{-6}\\
\text{objective evaluation}
    &c_\Obj=4.1655\times10^{-6}\\
\text{joint objective--gradient evaluation}
    &c_{\Obj,\nabla}=4.3551\times10^{-6}\\
\text{one \(d\)-dimensional Gaussian vector}
    &c_Z=1.8629\times10^{-6}\\
\text{one scalar comparison}
    &c_{\rm cmp}=1.3082\times10^{-8}\\
\text{one state-vector copy}
    &c_{\rm copy}=2.1174\times10^{-8}
\end{array}
\]
These coefficients are obtained from separate microbenchmarks on the
machine used for the experiments.  They are not algorithmic tuning
parameters.

The state-copy cost is tracked only as a diagnostic and is not included
in the plotted work.  Therefore the work is a deterministic function of
the discrete run length.

\subsection{One Langevin exploration path}

The work for one exploratory Langevin trajectory at design
\((\temp,a)\) is
\begin{equation}
    w_{\temp,a}^0
    \eqdef
    c_{\Obj,\nabla}(n_{\temp,a}+1)
    +(c_Z+c_{\rm cmp})n_{\temp,a}.
    \label{eq:w0}
\end{equation}

The initial \(c_{\Obj,\nabla}\) term accounts for evaluating
\((\Obj,\nabla\Obj)\) at the initial state; every Euler step then
requires one joint objective--gradient evaluation, one Gaussian vector,
and one comparison with the running minimum.

\subsection{Gradient descent}

The computational work of the single exploitation phase is
\begin{equation}
    \wGD 
    \eqdef
    c_\nabla \nGD .
    \label{eq:wGD}
\end{equation}

\subsection{Total Langevin--gradient work}

The exact work charged by the code is
\begin{equation}
    \WLGRM
    \eqdef
    Nw_{\temp,a}^0
    +c_{\rm cmp}(N-1)
    +\wGD .
    \label{eq:WLGexact}
\end{equation}
The extra \(c_{\rm cmp}(N-1)\) accounts for selecting the best retained
state among the \(N\) exploratory trajectories.  Since this term is
negligible in the present experiments, the conceptual approximation is
\begin{equation}
    \WLGRM
    \approx
    Nw_{\temp,a}^0+\wGD .
    \label{eq:WLGapprox}
\end{equation}

Under the fixed stochastic step size, fixing \(a\) essentially fixes
\(\WLGRM\); temperature affects performance primarily
through the probability of success rather than through the work itself.

\end{document}